\documentclass[11pt]{article}
\usepackage{makecell}
\usepackage[utf8]{inputenc} 

\usepackage[T1]{fontenc}

\usepackage{amssymb}
\usepackage{bm}
\usepackage{latexsym,amsmath,url,epsfig}
\usepackage{textcomp}
\usepackage{amsfonts,euscript}
\usepackage{graphicx}
\usepackage{amsmath}
\usepackage{amssymb}
\usepackage{amsfonts}
\usepackage{amsthm}
\usepackage{comment}
\usepackage{mathrsfs}
\usepackage[all]{xy}
\usepackage{epstopdf}
\usepackage{subfigure}
\usepackage{rotating}
\usepackage{appendix}
\usepackage{multirow,rotating}
\usepackage{url}
\usepackage{subfigure}
\usepackage{setspace}
\usepackage{listings}
\usepackage{color}
\usepackage{tikz}
\usepackage{graphics}
\usepackage[breaklinks=true]{hyperref}
\usepackage[english]{babel}
\usepackage{datetime}
\usepackage{subfigure}
\usepackage{booktabs}
\usepackage{rotfloat}
\usepackage[numbers]{natbib}
\usepackage{algorithm}
\usepackage{algorithmicx}
\usepackage{algpseudocode}
\usepackage{framed}
\usepackage{xcolor}
\usepackage{csquotes}
\usepackage{arydshln}
\usepackage{tabularx}
\usepackage{dsfont}
\usepackage{comment}
\graphicspath{{pic/}}
\newtheorem{theorem}{Theorem}
\newtheorem{corollary}{Corollary}

\newtheorem{lemma}{Lemma}

\newtheorem{assumption}{Assumption}
\theoremstyle{definition}

\newtheorem{remark}{Remark}

\newcommand{\R}{\mathbb{R}}

\newcommand{\var}{\text{Var}}

\newcommand{\cP}{\mathcal{P}}

\newcommand{\E}{\mathbb{E}}

\algdef{SE}[DOWHILE]{Do}{doWhile}{\algorithmicdo}[1]{\algorithmicwhile\ #1}

\allowdisplaybreaks

\usepackage{tikz}
\usepackage{pgfplots}
\pgfplotsset{compat=1.18} %

\begin{document}
\title{Bayesian Distributionally Robust Merton Problem with Nonlinear Wasserstein Projections}
\author{ 
   Jose Blanchet\thanks{\scriptsize Stanford University, CA 94305, US. Email: \texttt{jose.blanchet@stanford.edu}}\and
   Jiayi Cheng\thanks{\scriptsize New York University, NY 10003, US. Email: \texttt{jiayicheng@nyu.edu}}
    \and Hao Liu\thanks{\scriptsize Stanford University, CA 94305, US. Email: \texttt{haoliu20@stanford.edu}} 
    \and 
    Yang Liu\thanks{\scriptsize The Chinese University of Hong Kong (Shenzhen), Guangdong 518172, China. Email: \texttt{yangliu16@cuhk.edu.cn}}
        }
\date{}%
\maketitle

\begin{abstract}
We revisit Merton's continuous-time portfolio selection through a data-driven {lens of distributionally robust optimization (DRO)}. Our aim is to tap the benefits of frequent trading over short horizons while acknowledging that drift is hard to pin down, whereas volatility can be screened using realized or implied measures for appropriately selected assets. Rather than time-rectangular distributionally robust control (DRC), which replenishes adversarial power at every instant and induces over-pessimism, we place a single ambiguity set on the drift prior within a Bayesian Merton model. This prior-level ambiguity preserves learning and tractability: a minimax swap reduces the robust control to optimizing a nonlinear functional of the prior, enabling a Karatzas and Zhao \cite{KZ98}-type closed-form evaluation for each candidate prior. We then characterize small-radius worst-case priors under Wasserstein uncertainty via an explicit asymptotically optimal pushforward of the nominal prior, and we calibrate the ambiguity radius through a nonlinear Wasserstein projection tailored to the Merton functional. Synthetic and real-data studies demonstrate reduced pessimism relative to DRC and improved performance over myopic DRO–Markowitz under frequent rebalancing.

\vspace{0.5em}

\vspace{0.5em}
\end{abstract}

\section{Introduction}\label{intro}

This paper revisits Merton’s classical continuous-time portfolio selection model through a data-driven,
distributionally robust lens. Our goal is to tap the benefits of frequent trading over short horizons (days to weeks) while acknowledging that drift is difficult to pin down over such horizons. In contrast, volatility can often be stabilized by screening via realized or implied measures for appropriately chosen assets\footnote{We acknowledge that this screening can limit the investing universe; we discuss ways to incorporate volatility uncertainty in the conclusion.}. We therefore focus on the robustness of the drift and explicitly separate uncertainty modeling: we treat volatility as pre-estimated (e.g., realized/implied volatility) within the Bayesian filtering setup, and we robustify only the drift by placing a single ambiguity set on the drift prior.

Our starting point is the observation that, while data-driven robust portfolio selection performs competitively, it is static and myopic; it does not leverage the value of frequent rebalancing that Merton's continuous-time model affords. In particular, for example, myopic DRO--Markowitz policies do not capitalize on information revealed over time except through rolling re-estimation. At first sight, it may be surprising that DRO--Markowitz often outperforms dynamic investment strategies in practice (as illustrated in \citep{BlanchetChenZhou2021}). As explained there, sophisticated dynamic decisions rely heavily on model assumptions; when these are violated, errors compound over time, undermining dynamic strategies, especially under unpredictable non-stationarities. By contrast, the Merton framework provides a principled way to exploit frequent trading, provided we handle model uncertainty in a way that does not induce excessive pessimism.

The robust control literature (DRC/DRMDP{; DR Control/Markov Decision Processes}) places time-rectangular ambiguity on the data-generating process and derives policies via dynamic programming; see, e.g., \cite{HansenSargent2001,hansen2008robustness,RB1,RB2} and subsequent developments in DRMDP/DRRL {(Reinforcement Learning)} \cite{NianSi,wang2022policy,Wang2023,Liu2022,Zhou2021,lu2024drrl}. Rectangularity replenishes the adversary's power at every time step, often yielding over-conservative allocations when applied to portfolio choices with short horizons. Intuitively, in a one-dimensional drift-shift toy example, a rectangular adversary can depress the drift at each instant, compounding pessimism over time; by contrast, a prior-level ambiguity perturbs the drift distribution once. Nevertheless, rectangularity is widely used because it preserves time consistency and the dynamic-programming structure (Bellman equations), which confers strong tractability and algorithmic scalability, the very benefit delivered by the “replenishing” mechanism. Our design choice is different: instead of time-rectangular uncertainty, we adopt a \,\emph{prior-level} ambiguity in a Bayesian Merton model. We place a single ambiguity set around the drift prior (primarily a Wasserstein ball%
). {However, unlike rectangular formulations, this choice does not build Bellman time consistency into the robust layer of the problem. Accordingly, the robust solution studied in this paper should be interpreted in an ex-ante, or pre-committed, sense: we optimize the initial robust criterion at time zero and then implement the resulting observation-adapted \textit{feedback} policy over the investment horizon. Importantly, the resulting robust portfolio policy itself retains non-trivial feedback features adapted to the stock-price filtration, as in the underlying Bayesian Merton problem. The philosophy is that the portfolio manager faces an environment in which the drift is roughly constant over a fixed time horizon, say $T$, but drawn from a distribution that can be calibrated based on historical data. However, this is a coarse approximation of the more realistic situation in which the drift, while uncertain, is not exactly constant and that the prior distribution is difficult to calibrate.}
The volatility is assumed to be constant and therefore easy to estimate in the continuous-time setting. This \emph{distributionally robust Bayesian control} (DRBC) design reduces pessimism and preserves the learning structure of the Bayesian Merton formulation.  

Technically, an important tractability lever is a minimax swap (Sion-type) that holds for broad prior-level ambiguity sets (including Wasserstein%
). This swap allows us to evaluate, for any fixed prior in the ambiguity set, the optimal Bayesian Merton value and policy in closed form using the formulas of \cite{KZ98}. As a result, the DRBC game reduces to a constrained distributional optimization over the drift prior: we optimize a nonlinear functional of the prior that arises from Karatzas and Zhao \cite{KZ98}'s expression. The reduction holds under mild conditions standard in the Bayesian Merton literature. \footnote{Even if the minimax swap is difficult to justify, one may start with the formulation in which the adversary moves first. While this is not the most natural formulation (because the adversary typically models an environment that occurs after the manager makes its decision), still, it may be a pragmatic way to induce robustness while mitigating overconservative policies.}

Optimizing over a Wasserstein ball is subtle here because the resulting objective is highly nonlinear in the prior; standard Wasserstein DRO tools do not apply off-the-shelf. We derive small-radius asymptotics for non-linear functionals and apply these results to the worst-case prior. We obtain a \emph{constructive} approximation: an explicit asymptotically optimal pushforward perturbation of the nominal prior that realizes the first-order effect. Beyond optimization, we also address calibration: we select the ambiguity radius by general \emph{nonlinear Wasserstein projections}. Then, we tailor these general results to the Merton functional, extending linear RWPI-style projection ideas to this nonlinear setting. We highlight that this nonlinear projection perspective may be of independent interest, given the broad and growing use of Wasserstein projections \citep{SiMurthyBlanchetNguyen2021,Blanchet2021WassersteinDRO}.

We complement the theory with evidence in synthetic and real-data settings. In synthetic experiments, we simulate an environment consistent with Merton's assumptions: volatilities are known (or well-estimated), while asset drifts are unknown but deterministic and time-varying, generated from sinusoidal bases spanning a wide range of oscillatory frequencies. We compare two empirical strategies to construct the nominal drift prior from data: (i) batched, disjoint time windows (e.g., days or weeks) that form an empirical prior from window-level average returns; and (ii) day-of-week aggregation within larger windows (e.g., averages of ``Mondays'', ``Tuesdays'', etc.). Perhaps surprisingly, the batched-window prior performs slightly better even under periodic drifts, and we adopt it in our real-data study. Across both synthetic and real data, DRBC exhibits reduced pessimism compared to DRC under matched radii and improves performance over myopic DRO--Markowitz when frequent rebalancing is possible.

\paragraph{Contributions.} Our main contributions are as follows.
\begin{itemize}
  \item \textbf{Duality for prior-level ambiguity.} We prove a minimax swap for Wasserstein drift-prior ambiguity sets,  %
  reducing DRBC--Merton to optimizing a nonlinear functional of the prior while preserving closed-form evaluation via \cite{KZ98}.
  \item \textbf{Constructive worst-case prior sensitivity under Wasserstein and non-linear functionals.} For small ambiguity radii, we derive a first-order expansion of the robust objective and give an explicit asymptotically optimal push-forward perturbation of the nominal prior. These results are of independent interest since they are derived for the evaluation of the worst-case non-linear functions of probabilities.
  \item \textbf{Calibration via nonlinear Wasserstein projections.} We introduce a projection-based, data-driven method to select the ambiguity radius tailored to the Merton functional, generalizing linear RWPI-style projections to a nonlinear manifold. Again, the results involve general nonlinear projections in Wasserstein geometry which are of independent interest.
  \item \textbf{Empirical validation and reduced pessimism.} Synthetic and real-data experiments demonstrate reduced over-conservatism relative to DRC and improved performance over myopic DRO--Markowitz under frequent rebalancing.
\end{itemize}

{Our portfolio strategy can still be interpreted within a Bayesian lens. What we offer is a systematic approach to infuse both robustness and tractability in the choice of the prior, situating our framework within the scope of contemporary tools of distributionally robust decision making. From the Bayesian perspective, a statistician used to study Bayesian estimators from a frequentist standpoint may object that the prior's choice typically induces a bias that is relatively small as time increases. But this is not the environment we have in mind. In our setting, the investment horizon (which we denote by $T$) is fixed. This time horizon is long enough so that the manager may take advantage of multi-stage, even frequent, decisions but short enough that the volatility and drift are roughly constant. With this in mind, the volatility can be treated as fixed and the drift is unknown, so a Bayesian setting is natural but with a prior that requires robust calibration. This is precisely the mindset that motivates our development.}

\paragraph{Related work.} Distributionally robust optimization (DRO) has been extensively studied in statistics and machine learning, including Wasserstein DRO \cite{Blanchet2021WassersteinDRO,Blanchet2024DRO} and surveys \cite{rahimian2019distributionally,Bayraksan2015,ksw_2024_dro}. In control, DRC and DRMDP formulations typically adopt a time-rectangular ambiguity that preserves Bellman dynamic programming and tractability, but can be conservative because the adversary constantly replenishes its power  \cite{HansenSargent2001,hansen2008robustness,wkr_2013_rmdp,RB1,RB2,NianSi,wang2022policy,Wang2023,Liu2022,Zhou2021,lu2024drrl}. We instead place a single ambiguity set on the drift prior in a Bayesian Merton model, enabling closed-form evaluation while tempering pessimism. Our approach relates to Bayesian DRO in static settings \cite{doi:10.1137/21M1465548} and differs from DRBO \cite{pmlr-v108-kirschner20a}: beyond being online and discrete-time, DRBO does not robustify within a stochastic control framework (it is closer to rolling-horizon risk minimization) and provides limited guidance for selecting the ambiguity size.

On the sensitivity-analysis side, our work builds on the Wasserstein-DRO expansions of \cite{BartlDrapeauOblojWiesel2021}, who develop first-order asymptotics and optimal perturbations under \text{Wasserstein} balls, and to subsequent statistical analyses of \text{Wasserstein} estimators such as \cite{BKW19, Blanchet2021WassersteinDRO, BMZ21}. A parallel line of work studies divergence-based robustness, most notably the KL- and $\phi$-divergence sensitivity framework of \cite{Lam2016} and the Rényi-divergence bounds of \cite{AtarChowdharyDupuis2015}. In contrast to these approaches, which focus on linear or convex performance measures, we analyze the sensitivity of a highly non-linear Merton value functional in \text{Wasserstein} geometry and use its first-order expansion to construct problem-specific worst-case priors in continuous-time Bayesian control models. A broader robustness and sensitivity literature, ranging from ambiguity-averse portfolio selection \cite{PflugWozabal2007} to distributionally robust SAA {(Sample Average Approximation)}, stochastic programming, and adversarial training and hedging \cite{AndersonPhilpott2019, Dupacova1990, BonnansShapiro2013, ArmacostFiacco1974, AraujoEtAl2019, GaoKleywegt2016, sauldubois_touzi_2024}, provides a useful context but tackles problems of a different structural form.

On the projection side, most existing results focus on linear or convex functionals and thus differ from the non-linear, Merton-specific projection problem studied here. Linear optimal transport projections have been extensively analyzed in the optimal-transport DRO literature, including recent developments on small-sample behavior \cite{LinBlanchetGlynnNguyen2024}, unifying OT-based DRO reformulations \cite{BlanchetKuhnLiTaskesen2023}, and stability evaluations via distributional perturbations \cite{BlanchetCuiLiLiu2024}. Related projection methodologies also arise in confidence-region construction \cite{BlanchetMurthySi2022}, fairness testing through OT projections \cite{SiMurthyBlanchetNguyen2021, TaskesenBlanchetKuhnNguyen2021}, and martingale projections under adapted Wasserstein distances \cite{BlanchetWieselZhangZhang2024}. Unlike these works, which center on linear expectation functionals or convex risk measures, we analyze a fully non-linear projection defined by the Merton value functional and derive a problem-specific first-order expansion that yields constructive worst-case priors tailored to the continuous-time Bayesian control.

\paragraph{Organization.} Section~\ref{Merton} reviews the classical and Bayesian Merton formulations. Section~\ref{finmath} introduces the DRBC model, establishes the minimax swap enabling closed-form evaluation, and gives the small-radius Wasserstein approximations for worst-case priors. Section~\ref{projection} presents our data-driven construction of the drift prior and develops a nonlinear Wasserstein projection for ambiguity calibration. Sections~\ref{syn} and~\ref{real} report synthetic and real-data experiments, respectively. Additional experimental details are deferred to the Appendix. Due to space constraints and to ease the exposition, we focus on compactly supported adversarial priors. We provide an online supplementary material dealing with more general priors in the important power-utility setting.

\section{Preliminaries}\label{Merton}
In this section, we review the classical and Bayesian formulations of Merton's portfolio selection problem and set the notation used throughout. We denote by $W = \left(W_1, \ldots, W_d\right)^\top$ an $\mathbb{R}^d$-valued Brownian motion under a complete filtered probability space $(\Omega, \mathcal{F}, \{\mathcal{F}(t)\}_{t \in [0,T]}, \mathbb{P})$. The risk-free asset is given by $S_{0}(0) = s_0 > 0$ and 
$$
d S_{0}(t) = r S_{0}(t) d t, \;\; 0 \leq t \leq T,
$$
with risk-free rate $r> 0$, and $d$ risky assets are represented by the vector $S = \left(S_1, \ldots, S_d\right)^\top$. The dynamics of the risky assets follow the geometric Brownian motion: for $i = 1, \ldots,d$, $S_{i}(0) > 0$ and 
$$
 d S_{i}(t) = S_{i}(t)\left[b_idt + \sum_{j=1}^d\sigma_{ij}dW_{j}(t)\right],\;\; 0 \leq t \leq T {,}  
$$ {where $\sigma = \left(\sigma_{ij}\right)_{1\leq i,j\leq d}$ is assumed to be a known $d \times d$ full-rank matrix. While assuming that the volatility is known and fixed presents a limitation, %
we adopt this assumption in the same spirit as in the non-robust Bayesian formulation of \cite{KZ98}. The motivation is that this parameter is fully recovered instantly in continuous time. The way in which this requirement translates in practice is that assets should be screened so that high-frequency volatility estimates should remain relatively constant over the trading time horizon $T$.}

A portfolio (or control, or policy) is a stochastic process $\pi = \{\pi(t)\}_{t \in [0,T]}$ such that for a fixed $t \in [0,T]$, $\pi(t) = \left(\pi_{1}(t), \ldots, \pi_{d}(t)\right)^\top$ and $\pi_{i}(t)$ represents the amount of money invested in the $i$th stock at time $t$. This induces the dynamics of a controlled wealth process with $X(0) = x_0$ (we simplify the notation so that $X^{\pi}$ is written as $X$): 
\begin{equation}\label{wealthpm}
		d X(t) =  rX(t)dt + \pi(t)^\top\left(b-r\mathbf{1}\right)dt + \pi(t)^\top \sigma dW(t),
\end{equation}where $\mathbf{1} = \left(1, \ldots, 1\right)^\top$ is the vector of all ones.
We call an $\mathcal{F}$-progressively measurable (under $\mathbb{P}$ up to $\mathbb{P}$-null sets) stochastic process (control) $\pi= \{\pi(t)\}_{t \in [0,T]}$ admissible if $X(0) = x_0$, $\int_0^T \left \Vert \pi(t) \right\Vert_2^2dt < \infty$ $\mathbb{P}$-almost surely, and 
     Equation (\ref{wealthpm}) admits a unique strong solution with $X(t) > 0$ for any $t \in [0,T]$. The collection of all admissible controls is denoted by $\mathcal{A}(x_0)$.

The objective function of Merton's problem is $V(x_0) = \sup_{\pi \in \mathcal{A}(x_0)}\mathbb{E}_{\mathbb{P}}\left[U(X(T))\right],$ 
where $U$ is the utility function. We will specify the utility in Assumption \ref{ass:utility} later. %
Using techniques from dynamic programming, Merton's problem admits a closed-form formula of the optimal strategy $\pi$. In practice, the parameters $b$ and $\sigma$ are estimated from the market data with statistical techniques (e.g., maximum likelihood estimation), and then these estimates are plugged into the closed-form solutions.

However, in practice, estimating $b$ is difficult at the horizons of interest. To address this, \cite{KZ98} introduces a partially observed Bayesian variant: we keep the same probability space but model the instantaneous return as an unobservable random vector $B: \Omega \to \mathbb{R}^d$, independent of $W$ under $\mathbb{P}$, with prior distribution $\mu$. The price dynamics become the same SDE with $b$ replaced by $B$:
\begin{equation}\label{stock}
 d S_{i}(t) = S_{i}(t)\left[B_idt + \sum_{j=1}^d\sigma_{ij}dW_{j}(t)\right],\;\; 0 \leq t \leq T.  
\end{equation}
Equivalently, under $\mathbb P$, $B\sim\mu$ (some probability measure) and $B$ is independent of $W$; we refer to $\mu$ as the prior. We write $\mathcal N(m,\Sigma)$ for a Gaussian law with mean $m$ and covariance $\Sigma$, and $\varphi_s$ for the density of $\mathcal N(0, s I_d)$.

We denote the natural filtration of the process $S$ by $\mathcal{F}^S$ and denote the $\mathbb{P}$-augmentation of $\mathcal{F}^S$ by $\mathcal{F}^S_{\mathbb{P}}$ and this right-continuous completion is the observation filtration. Decisions are based only on information from stock prices. That is, in this case, the admissible controls are restricted to those that are $\mathcal{F}^S_{\mathbb{P}}$-progressively measurable and satisfy the same assumptions of integrability and SDE as before.

The full information of the Brownian motion and the random vector is encoded in the $\mathbb{P}$-augmentation of the enlarged filtration $\mathcal{G}^{B,W} = \{\mathcal{G}^{B,W}(t)\}_{t \geq 0}$ with 
$$\mathcal{G}^{B,W}(t) = \sigma\left(B, W(s), 0 \leq s \leq t\right) = \sigma(B)\vee \mathcal{F}^W(t),$$
where we denote this augmentation by $\mathcal{G}$. 
Therefore, for each $t \geq 0, \mathcal{F}^S_{\mathbb{P}}(t) \subset \mathcal{G}(t)$, where the inclusion can be strict. The wealth dynamics are given by the equation \begin{equation}\label{wealthp}
		d X(t) =  rX(t)dt + \pi(t)^\top\left(B-r\mathbf{1}\right)dt + \pi(t)^\top \sigma dW(t).
\end{equation}
We call an $\mathcal{F}^S_{\mathbb{P}}$-progressively measurable (under $\mathbb{P}$ up to $\mathbb{P}$-null sets) stochastic process $\pi= \{\pi(t)\}_{t \in [0,T]}$ { an admissible control} if $X(0) = x_0$, $\int_0^T \left \Vert \pi(t) \right\Vert_2^2dt < \infty$ $\mathbb{P}$-almost surely, and Equation (\ref{wealthp}) admits a unique strong solution with $X(t) > 0$ for any $t \in [0,T]$. The collection of all admissible controls is denoted by $\mathcal{A}_{\text{Bayesian}}(x_0)$.
The Bayesian diffusion control problem is defined as
\begin{equation}\label{OriginalKaratzas}{
  V_{\text{Bayesian}}(x_0) = \sup_{\pi \in \mathcal{A}_{\text{Bayesian}}(x_0)}\mathbb{E}_{\mathbb{P}}\left[U(X(T))\right]. }
\end{equation}

We will study a general utility family shown in Assumption \ref{ass:utility}. 

\begin{assumption}[Utility Function; refinement of Assumption 3.1 of \cite{KZ98}]
\label{ass:utility}
(1) The utility function $U: (0,\infty)\to\mathbb{R}$ is twice continuously differentiable,
strictly concave and strictly increasing, satisfies the Inada conditions $\lim_{x\downarrow 0}U'(x)=\infty,$ $\lim_{x\uparrow\infty}U'(x)=0,$
and satisfies a polynomial growth condition:
there exist constants $C_U>0$ and $p_U>0$ such that
{\[
  |U(x)| \;\le\; C_U\Bigl(1+x+x^{-p_U}\Bigr),
  \qquad x>0.
\]}
(2) Define the inverse marginal utility
$I=(U')^{-1}:(0,\infty)\to(0,\infty)$, which is strictly decreasing. 
We assume that $I$ and its derivative satisfy a two-sided polynomial growth bound, and that 
$I$ does not decay too fast at infinity: there exist constants $C_I>0$, $c_I>0$, $p_I>0$ {and $q_I > 0$} such that
{\[
  |I(y)| + |I'(y)|
  \;\le\; C_I\Bigl(1+y^{-p_I}\Bigr),
  \qquad y>0,
\]}
and
\[
  I(y) \;\ge\; c_I\,y^{{-q_I}},
  \qquad y>0.
\]
\end{assumption}

\begin{remark}
The growth conditions for $U$ and for $I=(U')^{-1}$ are stated separately because
bounds on $U$ do not in general imply corresponding bounds on $U'$ or on~$I$.
Nevertheless, for standard utility functions used in applications
(such as power/CRRA and logarithmic utility), both parts of
Assumption~\ref{ass:utility} are satisfied. 
\end{remark}

We now introduce the objects needed to state the Bayesian solution. Define 
$F(t,y) := F_{\mu}(t,y) := \int_{\mathbb{R}^d}  L_t(b,y) \mu(d b) = \E_{\mu}[L_t(B, y)]$   
with
\begin{equation}\label{eq:L_func}
L_t(b,y) = 
\exp\left( \langle \sigma^{-1}(b - r\mathbf{1}), y \rangle - \frac{1}{2} \|\sigma^{-1}(b - r\mathbf{1})\|^2 t \right), \;\; t \in (0, \infty), \; b \in \mathbb{R}^d, \; y \in \mathbb{R}^d. %
\end{equation}

{
For a fixed $b$, $L_t(b,y)$ is the complete-information likelihood-ratio kernel associated with the observation state $y$ under drift $b$. In particular, along the realized observation path $Y(t)$ (see Theorem \ref{thm:karatzas_zhao_general_compact}), one has the state-price density under complete information $\Lambda(t):=L_t(B,Y(t))^{-1}$ (notation of \cite{KZ98}). 
Averaging over the prior $\mu$ gives $F(t,y)$. So, along the realized observation path $Y(t)$, the quantity $F(t,Y(t))$ is the filtered (partially observed) likelihood-ratio process, or equivalently, the Bayesian evidence term for the hidden drift. Similarly, we have the filtered state-price density under incomplete information $\widehat\Lambda(t) :=F(t,Y(t))^{-1}$.
}

For $k>0$, $s\in[0,T]$, and $y\in\mathbb{R}^d$, set {the budget functional appearing in the duality step of the martingale method in \cite{KZ98}:}
\begin{equation}\label{eq:LKdef}
L(k;s,y):=
\begin{cases}
e^{-rs}\displaystyle\int_{\mathbb{R}^d} I\!\left( \frac{k e^{-rT}}{F(T, y+z)} \right) \varphi_s(z) \, d z, & s>0,\\[1.2em]
I\!\left( \frac{k e^{-rT}}{F(T, y)} \right), & s=0.
\end{cases}
\end{equation}

To simplify tedious technical discussions, we assume that the drift variable $B$ is compactly supported (Assumption \ref{ass:compact_B}). We provide more technical discussions under specific utilities with relaxed concentration assumptions of $B$ in supplementary materials. 

\begin{assumption}[Compact Support for $B$]
\label{ass:compact_B}
The random variable $B$ is compactly supported in Merton's model.  Specifically, %
there exists a compact set $K \subset \mathbb{R}^d$ such that $\mu(B \in K) = 1$ for any prior $\mu$. %
\end{assumption}

Under this setting, a standard Gaussian integral computation yields the following lemma.
\begin{lemma}\label{lem:Lregularity}
Under Assumptions~\ref{ass:utility} and \ref{ass:compact_B}, the map $L$ in \eqref{eq:LKdef} is finite for all $(k,s,y)\in(0,\infty)\times[0,T]\times\mathbb{R}^d$, continuously differentiable in $(k,s,y)$ on $(0,\infty)\times(0,T]\times\mathbb{R}^d$, and for each fixed \(s\in[0,T]\), the map \((k,y)\mapsto L(k;s,y)\) is twice continuously differentiable on \((0,\infty)\times\mathbb R^d\).
\end{lemma}

According to \cite{KZ98}, under Assumption \ref{ass:utility}, the strictly decreasing function
\begin{equation}\label{eq:Lagr}
\begin{aligned}
    k \mapsto e^{-rT} \int_{\mathbb{R}^d} I\left( \frac{ke^{-r T}}{F(T, z)} \right) \varphi_T(z) d z = L(k; T, 0)
\end{aligned}
\end{equation}
is continuous and maps $(0, \infty)$ onto itself. Thus, the equation $L(k; T, 0) = x_0$ is satisfied for a unique constant $\mathcal{K}(x_0) \in (0, \infty)$.  

The optimal solution is summarized in the following two results, first for the value function and later for the optimal strategy (i.e., the portfolio weights at every point in time); these results restate the characterization in \cite{KZ98} in our notation and compact form, and will be used repeatedly in our development.
\begin{theorem}[Karatzas--Zhao \cite{KZ98}'s Solution]
\label{thm:karatzas_zhao_general_compact}
Suppose that Assumptions \ref{ass:utility} and \ref{ass:compact_B} hold. 

(1) The optimal value function of Problem (\ref{OriginalKaratzas}) is given by
$$
V_{\text{Bayesian}}(x_0) = \int_{\mathbb{R}^d} F(T, z) U\left( I\left( \frac{\mathcal{K}(x_0) e^{-rT}}{F(T,z)} \right) \right) \varphi_T(z) dz,
$$
where $\mathcal{K}(x_0) > 0$ is the unique solution to the budget constraint:
\begin{equation}\label{eq:budget}
x_0 e^{rT} = \int_{\mathbb{R}^d} I\left( \frac{\mathcal{K}(x_0) e^{-rT}}{F(T,z)} \right) \varphi_T(z) dz.
\end{equation}

In particular, for power utility $U(x) = x^\alpha / \alpha$, this reduces to the classical Karatzas--Zhao \cite{KZ98}'s formula:
$$
V_{\text{Bayesian}}(x_0) = \frac{(x_0 e^{rT})^\alpha}{\alpha} \left( \int_{\mathbb{R}^d} F(T,z)^{\frac{1}{1-\alpha}} \varphi_T(z) dz \right)^{1-\alpha}.
$$
(2) The optimal fractions invested in each stock at time $t \in [0,T]$ are given by the vector
$$
\frac{\pi^*(t)}{X^*(t)} = {\left(\sigma^\top\right)}^{-1} (-\mathcal{K}(x_0)) e^{-rT}\cdot \frac{
\int_{\mathbb{R}^d} \frac{\nabla_z F\left(T,z+Y(t)\right)}{F\left(T,z+Y(t)\right)^2} \cdot I'\left( \frac{\mathcal{K}(x_0) e^{-rT}}{F(T,z+Y(t))} \right) \varphi_{T-t}(z) dz
}{
\int_{\mathbb{R}^d} I\left( \frac{\mathcal{K}(x_0) e^{-rT}}{F(T,z+Y(t))} \right) \varphi_{T-t}(z) dz
},
$$
where $Y(t) = \sigma^{-1} (B - r\mathbf{1}) t + W(t)$, $t \in [0, T]$. Note that the filtration generated by $\{Y(t)\}_{t \in [0,T]}$ is equal to $\mathcal{F}^S$ under $\mathbb{P}$.

In particular, for power utility $U(x) = x^\alpha / \alpha$, this reduces to the classical formula:
$$
\frac{\pi^*(t)}{X^*(t)} = {\left(\sigma^\top\right)}^{-1} \cdot \frac{
\int_{\mathbb{R}^d} \nabla_z F\left(T,z+Y(t)\right) \cdot F\left(T,z+Y(t)\right)^{\frac{\alpha}{1-\alpha}} \varphi_{T-t}(z) dz
}{
(1-\alpha) \int_{\mathbb{R}^d} F\left(T,z+Y(t)\right)^{\frac{1}{1-\alpha}} \varphi_{T-t}(z) dz
}.
$$
\end{theorem}

In practice, the prior distribution is chosen by experts and other available information, and the fraction of investment into risky assets is computed via the formula provided by Theorem \ref{thm:karatzas_zhao_general_compact} plugging in the observations of stock prices, normalizing the weight to maintain self-financing. %

\section{Formulation and Main Structural Results}\label{finmath}
In this section, we introduce our formulation of distributionally robust Bayesian control (DRBC) of Merton's problem, compare it with the classical DRC methods, and show the tractability of the DRBC formulation. We focus on optimal-transport-based uncertainty sets and, as mentioned in the Introduction, we will revisit $\phi$-divergence extensions in a later section. The work of \cite{Blanchet2025Duality} considers duality results for DRBC with $\phi$-divergence uncertainty. %
{In this setting, the admissible controls are still defined with respect to the prior distribution, since all other probability measures in the uncertainty set are absolutely continuous to the nominal prior, and thus the formulation is still well-defined.} 
However, this is not immediate in the Wasserstein uncertainty case that we study here, so some care is needed to handle this. %

We use $D_c$ to denote the optimal transport discrepancy generated by $c$ (equivalent to the Wasserstein distance when $c$ is a metric), as we now describe. Precisely, let
$\mathcal{P}(\mathbb{R}^{d}\times\mathbb{R}^{d})$ be the space of Borel
probability measures supported on $\mathbb{R}^{d}\times\mathbb{R}^{d}$. A
given element $\upsilon\in\mathcal{P}(\mathbb{R}^{d}\times\mathbb{R}^{d})$ is
associated to a random vector $\left(U,V\right)  $, where $U$ and $V$ take values in $\mathbb{R}^{d}$, in the following way: $\upsilon
_{U}\left(A\right) = \upsilon\left(  A\times\mathbb{R}^{d}\right)$ and $\upsilon
_{V}\left(  A\right) = \upsilon\left(  \mathbb{R}^{d}\times A\right)$ for every
Borel set $A\subset\mathbb{R}^{d}$, where $\upsilon_{U}$ and $\upsilon_{V}$ 
respectively denote the marginal distributions of $U$ and $V$ under $\upsilon$.

To define $D_c$, we need to introduce a cost function $c:\mathbb{R}^{d}\times\mathbb{R}%
^{d}\rightarrow\lbrack0,\infty]$, where we shall assume that $c$ is lower
semicontinuous and $c\left(  u,u\right)  =0$ for any $u\in
\mathbb{R}^{d}$. Finally, given two probability distributions ${P}_1$ and ${P}_2$ supported on
$\mathbb{R}^{d}$ and a cost function $c$, we define
\begin{equation}\label{ot}
D_{c}({P}_1, {P}_2):=\inf\{\mathbb{E}_{\upsilon}[c(U,V)]:\upsilon\in\mathcal{P}(\mathbb{R}^{d}%
\times\mathbb{R}^{d}),\upsilon_{U}={P}_1,\upsilon_{V}={P}_2\},
\end{equation}
which can be interpreted as the optimal (minimal) transportation cost of
moving the mass from ${P}_1$ into the mass of ${P}_2$ under a cost $c\left(
x,y\right)  $ per unit of mass transported from $x$ to $y$. If for a given
$\rho\ge 1$, $c^{1/\rho} $ is a metric, then $D_{c}^{1/\rho}$ defines a metric on probability measures (the Wasserstein distance of order $\rho$); see \cite{villani2008optimal}.
Throughout the rest of the paper, we will choose the cost
function $c(u,v)=||v-u||_{2}^{2}$ when $B$ is compactly supported. Other choices of cost function are discussed in supplementary materials.

{To rigorously define the DRBC formulation, we need to make sure only the distribution of $B$ is changed and all other problem structures (e.g., the adaptedness of the controls, the distribution of the other sources of randomness, independence structures, and integrability conditions) are kept the same. 
To preserve the structure of the original partial observation problem, we keep track only of the joint law of $(B,W)$ and maintain their independence. 

We place the model on the canonical product space $\Omega=\mathbb{R}^d\times C([0,T];\mathbb{R}^d)$ equipped with the Borel $\sigma$-algebra $\mathcal{B}(\Omega)$, with canonical coordinates $(B,W)$ and write $S=S(B,W)$ for the price process. We let $\mathbb{P}_W$ be the Wiener measure on the space $C\left([0,T]; \mathbb{R}^d\right)$ and denote the prior distribution of $B$ as $\mathbb{P}_0$ (satisfying Assumption \ref{ass:compact_B}). The observation filtration is fixed once and for all as
$$\mathcal F^S(t):=\sigma\big(S(u):0\le u\le t\big)\vee\mathcal N,$$
where $\mathcal N$ collects the null sets, and it does not depend on the choice of the nominal prior $\mathbb P_0$ (nor on $\mu$ in the ambiguity set). %
The observation filtration is completed using the common null sets of the observation laws.
Indeed, for each fixed \(b\in K\),
\[
Y^b(t)=\sigma^{-1}(b-r\mathbf 1)t+W(t)
\]
has a law on \(C([0,T];\mathbb R^d)\) equivalent to Wiener measure. Hence all mixtures
over priors supported on \(K\) induce the same null sets on the observation path space.
The completed observation filtration is therefore prior-independent. A detailed construction
is provided in the supplementary material.
}

Now we define ${\mathbb{Q}_0} := \mathbb{P}_0 \otimes \mathbb{P}_W$. 
{For any reference prior \(\nu\in\mathcal P(K)\), define
\[
\mathcal{U}_{\delta,B}^{OT}(\nu)
:=
\{\mu\in\mathcal P(K):D_c(\mu,\nu)\le\delta\}.
\]
Specifically, define 
\[
\mathcal U_{\delta,B}^{\text{OT}}(\mathbb P_0)
:=
\{\mu\in\mathcal P(K):D_c(\mu,\mathbb P_0)\le \delta\}.
\]
For each \(\mu\in\mathcal U_{\delta,B}^{\text{OT}}(\mathbb P_0)\), write $
Q^\mu:=\mu\otimes \mathbb{P}_W$. 
Then the Wasserstein uncertainty set with a radius $\delta > 0$ is defined as:
\[
\mathcal U_\delta^{\text{OT}}(\mathbb P_0)
=
\{{Q}^\mu:\mu\in\mathcal U_{\delta,B}^{\text{OT}}(\mathbb P_0)\},
\]
where we assume that the prior measure and all adversarial measures have a compact support. 
Next,}
{
a control \(\pi\) is admissible if it is \(\mathcal F^S\)-progressively measurable and,
for every \(\mu\in\mathcal P(K)\), with \(Q^\mu:=\mu\otimes \mathbb P_W\),
\[
\int_0^T \|\pi(t)\|_2^2\,dt<\infty\qquad Q^\mu\text{-a.s.},
\]
the wealth equation admits a unique solution, and \(X(t)>0\) for all \(t\in[0,T]\),
\(Q^\mu\)-a.s. The collection of such controls with \(X(0)=x_0\) is denoted by
\(\mathcal A(x_0)\). }
The DRBC problem is defined by 
\begin{equation}\label{DRO}
  V(x_0) = \sup_{\pi \in \mathcal{A}(x_0)} 
\inf_{{Q} \in \mathcal{U}^{\text{OT}}_{\delta}(\mathbb{P}_0)} \mathbb{E}_{{Q}}\left[U(X^{\pi}(T))\right] 
\end{equation}
and we denote the optimal solution of Problem \eqref{DRO} as $\pi_{\text{DRBC}}$.

As discussed in \cite{Blanchet2025Duality}, Problem \eqref{DRO} is typically not directly tractable in its original form. To solve the problem, \cite{Blanchet2025Duality} discusses several techniques, such as discretizing the prior distribution in the context of KL-uncertainty sets or applying simulation techniques that may be time-consuming. However, exploiting the special structure of the Merton problem, under reasonable assumptions, we can derive the following theorem, which will significantly simplify approximating the solution to Problem \eqref{DRO}.

\begin{theorem}[Min-Max Equality]
\label{thm:minmax}
Under Assumptions~\ref{ass:utility} and \ref{ass:compact_B}, for any initial wealth $x_0 > 0$, the following min-max equality holds:
$$
\sup_{\pi \in \mathcal{A}(x_0)} \inf_{{Q} \in \mathcal{U}^{\text{OT}}_{\delta}(\mathbb{P}_0)} \mathbb{E}_{Q}[U(X^\pi(T))] = \inf_{{Q} \in \mathcal{U}^{\text{OT}}_{\delta}(\mathbb{P}_0)} \sup_{\pi \in \mathcal{A}(x_0)} \mathbb{E}_{Q}[U(X^\pi(T))].
$$
\end{theorem}

\begin{proof}
We prove the theorem in five steps, leveraging the compact support of $B$ and the structure of the utility function.

\textbf{Step 1: Weak duality.} The inequality $\sup_\pi \inf_{Q} \leq \inf_{Q} \sup_\pi$ holds by definition (weak duality). {For $\pi\in \mathcal{A}(x_0)$, write $\pi^b(t,w):=\pi\bigl(t;(b,w)\bigr)$ and $X^{\pi,b}(T,w):=X^\pi\bigl(T;(b,w)\bigr).$ 
We proceed to prove the reverse inequality. 
Fix constants \(r_+>1\), \(r_->p_U\) and constants %
\(M_0,M_+,M_->0\). %
Let 
\(
\mathcal{A}'_{r_+,r_-,M_0,M_+,M_-}(x_0)
\) 
denote the set of all \(\pi\in \mathcal{A}(x_0)\) satisfying the following conditions: %
\begin{itemize}
    \item[(R1)] there exists a set $N_\pi\subset C([0,T];\mathbb R^d)$ with
    $\mathbb P_W(N_\pi)=0$ such that, for every $w\notin N_\pi$, the map $b\longmapsto \pi^b(\cdot,w)$
    is continuous from $K$ into $L^2([0,T];\mathbb R^d)$;
    \item[(R2)] the random variable $G_\pi(w):=\int_0^T \sup_{b\in K}\|\pi^b(t,w)\|^2\,dt$
    belongs to $L^1(\mathbb P_W)$ and $\mathbb E_{\mathbb P_W}[G_\pi]\le M_0$; %
    \item[(R3)] %
    \[
    \sup_{b\in K}\mathbb E_{\mathbb P_W}\!\left[(X^{\pi,b}(T))^{r_+}\right]\le M_+,
    \qquad
    \sup_{b\in K}\mathbb E_{\mathbb P_W}\!\left[(X^{\pi,b}(T))^{-r_-}\right]\le M_-.
    \]
\end{itemize}
{For notational simplicity, we write \(\mathcal{A}'(x_0)\) for this fixed class.} 
We equip $\mathcal{A}'(x_0)$ with the metric
\[
d(\pi,\pi') := \sup_{{Q} \in \mathcal{U}^{\text{OT}}_{\delta}(\mathbb{P}_0)}\left( \mathbb{E}_{Q} \left[  \int_0^T \|\pi(t) - \pi'(t)\|^2 dt \right]\right)^{1/2}
\] and the topology induced by the $L^2$-norm.
By construction, $\mathcal{A}'(x_0)$ is convex. 
{Indeed, if \(\pi^1,\pi^2\in A'(x_0)\) and \(\lambda\in[0,1]\), then
\(
\pi^\lambda:=\lambda\pi^1+(1-\lambda)\pi^2
\)
satisfies (R1) by linearity. Moreover,
\(
G_{\pi^\lambda}\le \lambda G_{\pi^1}+(1-\lambda)G_{\pi^2},
\) 
so (R2) holds. Since the wealth equation is affine in \(\pi\), 
\(
X^{\pi^\lambda}
=
\lambda X^{\pi^1}+(1-\lambda)X^{\pi^2},
\) 
and positivity and the moment bounds in (R3) follow from convexity of
\(x\mapsto x^{r_+}\) and \(x\mapsto x^{-r_-}\).} 
Moreover, if ${Q}^\mu =\mu\otimes\mathbb P_W\in \mathcal{U}^{\text{OT}}_{\delta}  (\mathbb P_0)$, then $\mu(K)=1$.  Hence (R3) implies
\[
\mathbb E_{{Q}^\mu}\!\left[(X^\pi(T))^{r_+}\right]
=
\int_K \mathbb E_{\mathbb P_W}\!\left[(X^{\pi,b}(T))^{r_+}\right]\mu(db)
\le M_+. 
\]
Similarly, 
\[
\mathbb E_{{Q}^\mu}\!\left[(X^\pi(T))^{-r_-}\right]\le M_-.
\]
Thus the sliced bounds in (R3) automatically imply the corresponding moment bounds under
every ${Q}\in \mathcal{U}^{\text{OT}}_{\delta} (\mathbb P_0)$.}
{Although the metric \(d\) depends on the fixed ambiguity radius \(\delta\), the classes
\(\mathcal{A}(x_0)\) and \(\mathcal{A}'(x_0)\) themselves do not.}
As $\mathcal{A}'(x_0)\subset\mathcal{A}(x_0)$, we have 
$$
\sup_{\pi \in \mathcal{A}(x_0)} \inf_{{Q}\in\mathcal{U}^{\text{OT}}_{\delta}(\mathbb{P}_0)} \mathbb{E}_{Q}[U(X^\pi(T))]
\;\ge\; \sup_{\pi \in \mathcal{A}'(x_0)} \inf_{{Q}\in\mathcal{U}^{\text{OT}}_{\delta}(\mathbb{P}_0)} \mathbb{E}_{Q}[U(X^\pi(T))].
$$
In the next steps, we establish the corresponding reverse inequality for $\mathcal{A}'(x_0)$. The final step (discussed later) will verify that the optimizer obtained at the end of the analysis belongs to $\mathcal{A}'(x_0)$, thereby closing the gap.

\textbf{Step 2: Variation of constants and continuity in the model.} 
Let $K\subset\mathbb{R}^d$ be the compact support guaranteed by Assumption~\ref{ass:compact_B}, fix $\pi\in\mathcal{A}'(x_0)$, and note that the linear wealth dynamics yield the variation-of-constants formula
\begin{equation}\label{Eq_VIC1}
{}
{X^{\pi,b}(T)} = x_0e^{rT} + e^{rT}\int_0^T e^{-rs}{\pi^b(s,W)}^\top(b-r\mathbf{1})\,ds + e^{rT}\int_0^T e^{-rs}{\pi^b(s,W)}^\top\sigma\,dW(s).
\end{equation}
{We define the map
\[
g_\pi(b):=\mathbb E_{\mathbb P_W}\!\left[U\!\left(X^{\pi,b}(T)\right)\right], \qquad b\in K.
\]
Let $b_n\to b$ in $K$. By (R1),
\[
\int_0^T \|\pi^{b_n}(t,W)-\pi^b(t,W)\|^2\,dt \to 0
\qquad \mathbb P_W\text{-a.s.}
\]
Also,
\[
\int_0^T \|\pi^{b_n}(t,W)-\pi^b(t,W)\|^2\,dt
\le
2\int_0^T \sup_{\tilde b\in K}\|\pi^{\tilde b}(t,W)\|^2\,dt
=
2G_\pi(W),
\]
and $G_\pi\in L^1(\mathbb P_W)$ by (R2). Hence, the dominated convergence theorem gives
\[
\mathbb E_{\mathbb P_W}\!\left[\int_0^T
\|\pi^{b_n}(t,W)-\pi^b(t,W)\|^2\,dt\right]\to 0.
\]
Using the variation-of-constants formula, Jensen's inequality, and It\^o's isometry, 
it follows that
\[
X^{\pi,b_n}(T) \to X^{\pi,b}(T)
\qquad\text{in }L^2(\mathbb P_W).
\]
In particular, after passing to a subsequence if necessary,
\[
X^{\pi,b_n}(T) \to X^{\pi,b}(T)
\qquad \mathbb P_W\text{-a.s.}
\]
Now the growth envelope on $U$ and the strict inequalities
$r_+>1$ and $r_->p_U$ imply that $\left\{U\!\left(X^{\pi,\widetilde b}(T)\right):
\widetilde b\in K\right\}$
is uniformly integrable under $\mathbb{P}_W$. Indeed, choose
$\eta>0$ such that $1+\eta<r_+$ and $(1+\eta)p_U<r_-.$
By Assumption~\ref{ass:utility}, there exists $C_\eta>0$ such that
\[
|U(x)|^{1+\eta}
\le
C_\eta\left(
1+x^{1+\eta}+x^{-(1+\eta)p_U}
\right),
\qquad x>0.
\]
Therefore, by (R3),
\[
\sup_{\widetilde b\in K}
\mathbb{E}_{\mathbb{P}_W}
\left[
\left|U\!\left(X^{\pi,\widetilde b}(T)\right)\right|^{1+\eta}
\right]
<\infty.
\]
Hence the family is uniformly integrable. Consequently,
$g_\pi(b_n)\to g_\pi(b)$.} 
Thus $g_\pi$ is continuous on the compact set $K$.
{Hence, for fixed \(\pi\) and \({Q}^\mu=\mu\otimes \mathbb{P}^W\),
\[
\mathbb{E}_{{Q}^\mu}[U(X^\pi(T))]
=
\int_K g_\pi(b)\,\mu(db),
\]
and the right-hand side depends continuously on \(\mu\) under weak convergence because
\(g_\pi\) is continuous and bounded on \(K\).} 
This proves the desired continuity in the model variable.

\textbf{Step 3: Concavity and continuity in the control.}
For fixed ${Q} \in \mathcal{U}^{\text{OT}}_{\delta}(\mathbb{P}_0)$, the map
$\pi \mapsto \mathbb{E}_{Q}[U(X^\pi(T))]$ is concave because $U$ is concave
and the state equation is affine in $\pi$.

To show continuity on $(\mathcal{A}'(x_0),d)$, let $\pi_n \to \pi$ in the topology
induced by $d$, and set $\Delta\pi(t) = \pi_n(t)-\pi(t)$. 
{Using the linear wealth dynamics, we may write the unsliced version of
representations~\eqref{Eq_VIC1}: 
\[
X^{\pi_n}(T)-X^\pi(T)
=
e^{rT}\int_0^T e^{-rs}\Delta\pi(s)^\top(B-r\mathbf{1})\,ds
+
e^{rT}\int_0^T e^{-rs}\Delta\pi(s)^\top\sigma\,dW(s).
\]} Taking quadratic moments,
and using $(a_1+a_2)^2 \le 2a_1^2 + 2a_2^2$, we obtain
\[
\mathbb{E}_{Q}\!\big[|X^{\pi_n}(T) - X^\pi(T)|^2\big]
\;\le\; C \,\mathbb{E}_{Q}\!\left[\int_0^T \|\Delta\pi(s)\|^2 ds\right],
\]
where $C = 2 e^{2rT}\Big( T \sup_{b\in K} \|b - r\mathbf{1}\|^2
\;+\; \|\sigma \sigma^\top\|_{\mathrm{op}}\Big)$
is independent of ${Q}$. Since $d(\pi_n,\pi) \to 0$, the right-hand side
converges to zero uniformly over ${Q}\in\mathcal{U}^{\text{OT}}_{\delta}(\mathbb{P}_0)$.
Hence $X^{\pi_n}(T) \to X^\pi(T) \text{ in } L^2({Q})$ (and thus in probability)
uniformly over ${Q}\in\mathcal{U}^{\text{OT}}_{\delta}(\mathbb{P}_0)$.

Next we establish uniform integrability
{ directly from (R3).
Since $\pi_n,\pi \in \mathcal{A}'(x_0)$ and the class $\mathcal{A}'(x_0)$ is defined with fixed constants
$M_+,M_->0$, Step~1 implies that for every fixed
${Q}\in\mathcal{U}^{\mathrm{OT}}_{\delta}(\mathbb{P}_0)$,
\[
\sup_{n\ge 1}\mathbb{E}_{{Q}}\!\left[(X^{\pi_n}(T))^{r_+}\right]\le M_+,
\qquad
\sup_{n\ge 1}\mathbb{E}_{{Q}}\!\left[(X^{\pi_n}(T))^{-r_-}\right]\le M_-,
\]
and the same bounds hold for $X^\pi(T)$.
By Assumption~\ref{ass:utility},
\[
|U(x)| \le C_U\bigl(1+x+x^{-p_U}\bigr),\qquad x>0.
\]
Choose $\eta>0$ such that
\[
1+\eta<r_+,
\qquad
(1+\eta)p_U<r_-.
\]
Then there exists a constant $C_\eta>0$ such that
\[
|U(x)|^{1+\eta}
\le
C_\eta\Bigl(1+x^{1+\eta}+x^{-(1+\eta)p_U}\Bigr),
\qquad x>0.
\]
Therefore,
\[
\sup_{n\ge 1}\mathbb{E}_{{Q}}\!\left[|U(X^{\pi_n}(T))|^{1+\eta}\right]
\le
C_\eta\left(
1+
\sup_{n\ge 1}\mathbb{E}_{{Q}}\!\left[(X^{\pi_n}(T))^{1+\eta}\right]
+
\sup_{n\ge 1}\mathbb{E}_{{Q}}\!\left[(X^{\pi_n}(T))^{-(1+\eta)p_U}\right]
\right)
<\infty,
\]
and similarly
\[
\mathbb{E}_{{Q}}\!\left[|U(X^\pi(T))|^{1+\eta}\right]<\infty.
\]
Hence, the family
\[
\{U(X^{\pi_n}(T))\}_{n\ge 1}\cup\{U(X^\pi(T))\}
\]
is uniformly integrable under ${Q}$.
Since $X^{\pi_n}(T)\to X^\pi(T)$ in probability under ${Q}$ and $U$ is continuous, we also have
\[
U(X^{\pi_n}(T))\to U(X^\pi(T))
\text{ in probability under }{Q}.
\]
Uniform integrability then yields
\[
\mathbb{E}_{{Q}}\!\left[\left|U(X^{\pi_n}(T))-U(X^\pi(T))\right|\right]\to 0,
\]
and in particular
\[
\mathbb{E}_{{Q}}\!\left[U(X^{\pi_n}(T))\right]
\longrightarrow
\mathbb{E}_{{Q}}\!\left[U(X^\pi(T))\right].
\]}
Thus, for each fixed ${Q}$, the map
$\pi\mapsto \mathbb{E}_{Q}[U(X^\pi(T))]$ is continuous on
$(\mathcal{A}'(x_0),d)$.

\textbf{Step 4: Apply Sion’s min-max theorem on $\mathcal{A}'(x_0)$.} The space $\mathcal{A}'(x_0)$ is convex, and $\mathcal{U}^{\text{OT}}_{\delta}(\mathbb{P}_0)$ is weakly compact since each admissible $B$-marginal is supported in the compact set $K$ and the Wasserstein %
balls are weakly closed (Prokhorov’s theorem).
Moreover, the map $(\pi, {Q}) \mapsto \mathbb{E}_{Q}[U(X^\pi(T))]$ satisfies:
\begin{itemize}
    \item \textbf{Concavity in $\pi$}: For fixed ${Q}$, $\pi \mapsto \mathbb{E}_{Q}[U(X^\pi(T))]$ is concave by the previous step.
    \item \textbf{Convexity in ${Q}$}: For fixed $\pi$, ${Q} \mapsto \mathbb{E}_{Q}[U(X^\pi(T))]$ is linear (hence convex).
    \item \textbf{Continuity in $\pi$}: For fixed ${Q}$, $\pi \mapsto \mathbb{E}_{Q}[U(X^\pi(T))]$ is continuous on $(\mathcal{A}'(x_0),d)$ (Step 3).
    \item \textbf{Continuity in ${Q}$}: For fixed $\pi$, ${Q} \mapsto \mathbb{E}_{Q}[U(X^\pi(T))]$ is continuous under the weak topology on $\mathcal{U}^{\text{OT}}_{\delta}(\mathbb{P}_0)$ by the argument in Step 2.
\end{itemize}
By Sion’s min-max theorem applied to $\mathcal{A}'(x_0)\times\mathcal{U}^{\text{OT}}_{\delta}(\mathbb{P}_0)$, we conclude:
$$
\sup_{\pi \in \mathcal{A}'(x_0)} \inf_{{Q} \in \mathcal{U}^{\text{OT}}_{\delta}(\mathbb{P}_0)} \mathbb{E}_{Q}[U(X^\pi(T))] = \inf_{{Q} \in \mathcal{U}^{\text{OT}}_{\delta}(\mathbb{P}_0)} \sup_{\pi \in \mathcal{A}'(x_0)} \mathbb{E}_{Q}[U(X^\pi(T))].
$$

\textbf{Step 5: Return to $\mathcal{A}(x_0)$.}
{Let $\mu$ be any admissible prior satisfying Assumption~\ref{ass:compact_B}, and consider the
singleton ambiguity set $\{{Q}^\mu\}$, where ${Q}^\mu:=\mu\otimes \mathbb{P}_W.$
Under the refined prior-independent observation filtration, the DRBC problem with uncertainty set
$\{{Q}^\mu\}$ reduces exactly to the fixed-prior Bayesian Merton problem. Hence
Theorem~\ref{thm:karatzas_zhao_general_compact} furnishes an optimal feedback control, which we denote by $\pi^{*,\mu}$. We now verify that the fixed-prior optimizer $\pi^{*,\mu}$ belongs to the same fixed class $\mathcal{A}'(x_0)$ with {constants $r_+,r_-$ and constants \(M_0,M_+,M_-\) chosen large enough and independent of $\mu$. }
For clarity, write
\[
F_\mu(T,y):=\int_K L_T(b,y)\,\mu(db),
\]
and define
\[
\Theta_\mu(t,y):=\int_{\mathbb{R}^d}
I\!\left(\frac{\mathcal{K}_\mu(x_0)e^{-rT}}{F_\mu(T,z+y)}\right)\varphi_{T-t}(z)\,dz,
\]
\[
\Xi_\mu(t,y):=\int_{\mathbb{R}^d}
\frac{\nabla_u F_\mu(T,z+y)}{F_\mu(T,z+y)^2}
\,I'\!\left(\frac{\mathcal{K}_\mu(x_0)e^{-rT}}{F_\mu(T,z+y)}\right)\varphi_{T-t}(z)\,dz,
\]
where $\mathcal{K}_\mu(x_0)>0$ is the unique multiplier in the budget equation under prior $\mu$.
By Theorem~\ref{thm:karatzas_zhao_general_compact},
\[
\frac{\pi^{*,\mu}(t)}{X^{*,\mu}(t)}
=
(\sigma^\top)^{-1}\bigl(-\mathcal{K}_\mu(x_0)\bigr)e^{-rT}
\frac{\Xi_\mu(t,Y(t))}{\Theta_\mu(t,Y(t))}.
\]
On the other hand, by the definition of the budget functional,
\[
X^{*,\mu}(t)
=
L_\mu\!\bigl(\mathcal{K}_\mu(x_0);\,T-t,\,Y(t)\bigr)
=
e^{-r(T-t)}\Theta_\mu(t,Y(t)).
\]
Therefore, 
\[
\pi^{*,\mu}(t)
=
-(\sigma^\top)^{-1}\mathcal{K}_\mu(x_0)e^{-r(2T-t)}\,\Xi_\mu(t,Y(t)).
\]
Set $A_K:=\sup_{b\in K}\|\sigma^{-1}(b-r\mathbf{1})\|<\infty$. 
Then for every $b\in K$ and $u\in\mathbb{R}^d$,
\[
L_T(b,u)
=
\exp\!\left(
\left\langle \sigma^{-1}(b-r\mathbf{1}),u\right\rangle
-\frac{T}{2}\|\sigma^{-1}(b-r\mathbf{1})\|^2
\right). 
\]
Hence, 
\[
e^{-A_K\|u\|-\frac{T}{2}A_K^2}
\le
L_T(b,u)
\le
e^{A_K\|u\|}.
\]
Integrating against $\mu$ gives constants $c_1,c_2>0$, depending only on
$K,T,r,\sigma$, such that for every admissible prior $\mu$ and every $u\in\mathbb{R}^d$,
\[
c_1 e^{-A_K\|u\|}
\le
F_\mu(T,u)
\le
c_2 e^{A_K\|u\|}.
\]
Moreover,
\[
\nabla_u F_\mu(T,u)
=
\int_K \sigma^{-\top}(b-r\mathbf{1})L_T(b,u)\,\mu(db),
\]
so there exists a constant $C_K>0$, independent of $\mu$, such that
\[
\frac{\|\nabla_u F_\mu(T,u)\|}{F_\mu(T,u)}\le C_K,
\qquad u\in\mathbb{R}^d.
\]
By Theorem~\ref{thm:karatzas_zhao_general_compact}, for each $\mu$ there is a unique multiplier $\mathcal{K}_\mu(x_0)>0$ solving the budget equation.
Using the uniform bounds on $F_\mu$ above, together with the monotonicity of $I$ and the same
argument used in the proof of Theorem~\ref{thm:karatzas_zhao_general_compact} (see \cite{KZ98}), there exist constants
$0<\underline{k}\le \overline{k}<\infty$
depending only on $x_0,T,r,\sigma,K$ and Assumption~\ref{ass:utility}, such that for every admissible $\mu$, $\underline{k}\le \mathcal{K}_\mu(x_0)\le \overline{k}$.
Using the bound on $\nabla_u F_\mu/F_\mu$, the upper-growth bound on $I'$ from
Assumption~\ref{ass:utility},
\[
|I'(y)|\le C_I\bigl(1+y^{-p_I}\bigr),
\qquad y>0,
\]
and the bounds on $F_\mu$ and $\mathcal{K}_\mu(x_0)$, there exist constants $C,c>0$, independent of
$\mu,t,y,z$, such that
\[
\left\|
\frac{\nabla_u F_\mu(T,z+y)}{F_\mu(T,z+y)^2}
I'\!\left(\frac{\mathcal{K}_\mu(x_0)e^{-rT}}{F_\mu(T,z+y)}\right)
\right\|
\le
C e^{c\|z+y\|}.
\]
Therefore, 
by a standard Gaussian integral and a dominated convergence argument, 
the map $(t,y)\mapsto \Xi_\mu(t,y)$ is continuous on $[0,T]\times\mathbb{R}^d$.}

{
Now we verify the (R1) condition.
Fix $w\in C([0,T];\mathbb{R}^d)$. For $b\in K$, define
\[
Y^b(t,w):=\sigma^{-1}(b-r\mathbf{1})t+w(t),\qquad 0\le t\le T.
\]
Then $b\longmapsto Y^b(\cdot,w)$
is continuous from $K$ into $C([0,T];\mathbb{R}^d)$ under the uniform norm. By the formula above,
\[
\pi^{*,\mu,b}(t,w)
=
-(\sigma^\top)^{-1}\mathcal{K}_\mu(x_0)e^{-r(2T-t)}\,
\Xi_\mu\!\bigl(t,Y^b(t,w)\bigr).
\]
Since $(t,y)\mapsto \Xi_\mu(t,y)$ is continuous, the map $(t,b)\longmapsto \pi_t^{*,\mu,b}(w)$
is continuous on the compact set $[0,T]\times K$. Hence, if $b_n\to b$ in $K$, then
\[
\sup_{0\le t\le T}
\bigl\|\pi^{*,\mu,b_n}(t,w)-\pi^{*,\mu,b}(t,w)\bigr\|
\to 0,
\]
and therefore
\[
\int_0^T
\bigl\|\pi^{*,\mu,b_n}(t,w)-\pi^{*,\mu,b}(t,w)\bigr\|^2\,dt
\to 0.
\]
Thus (R1) holds.}

{Next, we verify (R2). From the representation of $\pi^{*,\mu}$ and the bound on $\Xi_\mu$,
\[
\|\pi^{*,\mu,b}(t,w)\|
\le
C e^{c\|Y^b(t,w)\|},
\qquad 0\le t\le T,\ b\in K,
\]
with constants independent of $\mu$. Hence
\[
G_{\pi^{*,\mu}}(w)
:=
\int_0^T \sup_{b\in K}\|\pi^{*,\mu,b}(t,w)\|^2\,dt
\le
T C^2
\exp\!\left(
C_0+2c\|w\|_\infty
\right),
\]
where
\[
\|w\|_\infty:=\sup_{0\le t\le T}\|w(t)\|,
\qquad
C_0:=2cT\sup_{b\in K}\|\sigma^{-1}(b-r\mathbf{1})\|.
\]
Since the Brownian supremum on a finite horizon has finite exponential moments, there exists
$M_0<\infty$, independent of $\mu$, such that
\[
\mathbb{E}_{\mathbb{P}_W}\!\left[G_{\pi^{*,\mu}}\right]\le M_0.
\]
Thus (R2) holds.}

{Next, we verify (R3).
At terminal time,
\[
X^{*,\mu,b}(T,w)
=
I\!\left(
\frac{\mathcal{K}_\mu(x_0)e^{-rT}}{F_\mu(T,Y^b(T,w))}
\right).
\]
Using the upper-growth bound on $I$ from Assumption~\ref{ass:utility},
\[
|I(y)|\le C_I\bigl(1+y^{-p_I}\bigr),
\qquad y>0,
\]
together with the bounds on $F_\mu$ and $\mathcal{K}_\mu(x_0)$, there exist constants $C,c>0$,
independent of $\mu$ and $b$, such that
\[
X^{*,\mu,b}(T,w)\le C e^{c\|Y^b(T,w)\|}.
\]
Hence, for any chosen exponent $r_+>1$,
\[
\bigl(X^{*,\mu,b}(T,w)\bigr)^{r_+}
\le
C e^{c r_+\|Y^b(T,w)\|}.
\]
For the inverse moments, use the lower bound from Assumption~\ref{ass:utility}: 
\[
I(y)\ge c_I y^{-q_I},
\qquad y>0.
\]
Combining this with the bounds on $F_\mu$ and $\mathcal{K}_\mu(x_0)$ yield constants $C,c>0$,
independent of $\mu$ and $b$, such that
\[
X^{*,\mu,b}(T,w)\ge C e^{-c\|Y^b(T,w)\|}.
\]
Therefore, for any chosen exponent $r_->p_U$,
\[
\bigl(X^{*,\mu,b}(T,w)\bigr)^{-r_-}
\le
C e^{c r_-\|Y^b(T,w)\|}.
\]
Since
\[
Y^b(T) = W(T) + T\sigma^{-1}(b-r\mathbf{1}),
\]
and $K$ is compact, we have
\[
\sup_{b\in K}
\mathbb{E}_{\mathbb{P}_W}\!\left[e^{\lambda \|Y^b(T)\|}\right]
<\infty
\qquad\text{for every }\lambda>0.
\]
Consequently, there exist constants $M_+,M_->0$, independent of $\mu$, such that
\[
\sup_{b\in K}
\mathbb{E}_{\mathbb{P}_W}\!\left[
\bigl(X^{*,\mu,b}(T)\bigr)^{r_+}
\right]
\le M_+,
\qquad
\sup_{b\in K}
\mathbb{E}_{\mathbb{P}_W}\!\left[
\bigl(X^{*,\mu,b}(T)\bigr)^{-r_-}
\right]
\le M_-.
\]
Thus (R3) holds.
We have proved that $\pi^{*,\mu}$ satisfies (R1)--(R3) with constants independent of $\mu$.
Hence, 
\[
\pi^{*,\mu}\in \mathcal{A}'(x_0)
\qquad\text{for every admissible prior }\mu.
\]
Thus, for every admissible prior \(\mu\),
\[
\sup_{\pi\in \mathcal A'(x_0)} \mathbb{E}_{{Q}^\mu}[U(X^\pi(T))]
=
\sup_{\pi\in A(x_0)} \mathbb{E}_{{Q}^\mu}[U(X^\pi(T))],
\qquad {Q}^\mu:=\mu\otimes \mathbb P^W,
\]
because the fixed-prior Bayesian optimizer \(\pi^{*,\mu}\) belongs to \(\mathcal A'(x_0)\).
Combining this identity with the minimax equality on \(\mathcal A'(x_0)\) and weak duality yields the
desired minimax equality on \(\mathcal A(x_0)\).
}
Together with the inequalities at the beginning of the proof, this yields the
desired min–max equality on $\mathcal{A}(x_0)$.

\end{proof}

The value of Theorem \ref{thm:minmax} %
is that for a fixed probability measure ${Q}$, by Theorem \ref{thm:karatzas_zhao_general_compact}, the inner problem $\sup_{\pi \in \mathcal{A}(x_0)}  \mathbb{E}_{{Q}}\left[U(X^\pi(T))\right]$ has a closed-form solution in terms of the distribution of the drift $B$. 
{
In other words, the DRBC problem is equivalent to solving a constrained distributional optimization problem
\begin{align}\label{FWProblem}
&\inf_{{Q} \in \mathcal{U}^{\text{OT}}_{\delta}(\mathbb{P}_0)}\sup_{\pi \in \mathcal{A}(x_0)} 
 \mathbb{E}_{{Q}}\left[U(X^\pi(T))\right]
= \inf_{\mu \in \mathcal{P}(K): D_c(\mu,\mathbb{P}_0) \le \delta}\int_{\mathbb{R}^d} F_{\mu}(T, z) U\left( I\left( \frac{\mathcal{K}_{\mu}(x_0) e^{-rT}}{F_{\mu}(T,z)} \right) \right) \varphi_T(z) dz, 
\end{align}
where 
the uncertainty set on the right-hand side is only concerned with the distribution of $B$ and its objective  no longer depends explicitly on the Brownian motion $W$.} Hence, it suffices to optimize over the nonlinear functional and obtain the extreme probability measure {$\mu^{\text{WC}}$} (denote the functional on the right-hand side of Eq. \eqref{FWProblem} by a specific $J(\mu)$; see Corollary \ref{thm:asymptotic_perturbation_cor} for details):  
\begin{align}\label{FWProblem1new}
\mu^{\text{WC}} {\in} 
\arg\min_{\mu {\in \mathcal{P}(K)}: D_c(\mu,\mathbb{P}_0) \le \delta}{J}(\mu).
\end{align}
{Since Problem \eqref{FWProblem1new} is defined as a lower semicontinuous functional on a weakly compact set, the existence is guaranteed. While, in general, uniqueness is not guaranteed, when the utility is of the power case $U(x) = \frac{1}{\alpha} x^\alpha$, the reduced objective in Equation \eqref{FWProblem} has a form $J(\mu)= \frac{1}{\alpha} (\int_{\mathbb R^d} F_{\mu}(T,z)^{\frac{1}{1-\alpha}}\,\varphi_T(z)\,dz)^{1-\alpha},$ with $\mu \mapsto F_{\mu}(T,\cdot)$ affine
and uniqueness follows from the corresponding strict convexity/concavity of the reduced
power-utility functional, with the optimization direction adjusted according to the sign of
\(1/\alpha\).%
}

Next, we will first give a sensitivity analysis of this distributional optimization problem in greater generality for $J$, and then Problem \eqref{FWProblem1new} can be solved easily in an approximate sense by the following Theorem and the subsequent Corollary \ref{thm:asymptotic_perturbation_cor}. {We emphasize that this corollary provides only a first-order, small-radius characterization of the worst-case prior.}

\begin{theorem}[Nonlinear optimal perturbations]\label{thm:asymptotic_perturbation}
Fix $1<q<\infty$ and let $p$ be its H\"older conjugate.
Let $\|\cdot\|_q$ be a norm on $\R^d$ and $\|\cdot\|_p$ be its dual norm.
Consider the quadratic transport cost $c(x', x) := \tau\,\|x'-x\|_q^2$ with $\tau > 0$
and let $D_c$ denote the corresponding optimal-transport divergence
\[
  D_c(\mu,\nu)
  := \inf_{\pi\in\Pi(\mu,\nu)} \E_\pi\!\left[c(X',X)\right].
\]
For $\delta>0$, set $V_\delta := \sup_{\mu:D_c(\mu,\nu)\le\delta} J(\mu)$, where the functional
$J:\cP(\R^d)\to\R$ satisfies:

\begin{itemize}
  \item[(i)]
  For all $\mu$ in a neighbourhood of $\nu$, there exists a measurable
  $J'_\mu:\R^d\to\R$, differentiable in $x$ with gradient $\nabla J'_\mu(x)$,
  such that for every coupling $\pi$ of $(X',X)\sim(\mu,\nu)$ and
  $\nu_t := (1-t)\nu + t\mu$,
  \[
    J(\mu) - J(\nu)
    = \int_0^1 \E_\pi\big[ J'_{\nu_t}(X') - J'_{\nu_t}(X) \big]\,dt.
  \]

  \item[(ii)]
  (Regularity and growth of the Wasserstein gradient)
  The map $(\mu,x)\mapsto \nabla J'_\mu(x)$ is jointly continuous in the product
  topology (weak topology on $\mu$, Euclidean on $x$) in a neighbourhood
  of $(\nu,\cdot)$.
  Moreover, there exist $1\le r\le 2$ and $C>0$ such that for all such $\mu$ and all $x\in\R^d$,
  \[
    \|\nabla J'_\mu(x)\|_p \le C\big(1+\|x\|^{r-1}\big),
  \]
  and the reference law $\nu$ has a finite $2r$-moment:
  \[
    \E_\nu\|X\|^{2r} < \infty.
  \]
\end{itemize}

Let $g(x) := \nabla J'_\nu(x)$ and assume
$\E_\nu\|g(X)\|_p^2<\infty$.
Then, as $\delta\downarrow 0$,
\[
  V_\delta
  = J(\nu)
    + \sqrt{\frac{\delta}{\tau}}\,
      \Big( \E_\nu\big[\|g(X)\|_p^2\big] \Big)^{1/2}
    + o(\sqrt{\delta}).
\]
{
Moreover, set $
K_1:=\left(E_\nu\|g(X)\|_p^2\right)^{1/2}$. 
If \(K_1=0\), then
$
V_\delta=J(\nu)+o(\sqrt{\delta}),
$ 
and the identity map is asymptotically optimal. 
If \(K_1>0\), 
the Monge-type perturbation $T_\delta$ of the form
\[
    T_\delta(x)
  = x + \sqrt{\frac{\delta}{\tau}}\,
      \bar{\Delta}(x)
\]
is asymptotically optimal, where for each $x$ with $g(x)\neq 0$, we choose a measurable selection 
$u(x) \in \arg\max_{\|u\|_q=1} g(x)\cdot u$ and define $\bar{\Delta}(x) := \frac{\|g(x)\|_p}{K_1}\,u(x)$ with $\bar{\Delta}(x):=0$ when $g(x)=0$. 

}
\end{theorem}

\begin{proof}
Write $g(x) := \nabla J'_\nu(x)$ and 
$
  K_1 := \Big(\E_\nu\|g(X)\|_p^2\Big)^{1/2}.
$ 
If $K_1=0$, then $g(X)=0$ $\nu$-a.s.\ and \eqref{eq:Vdelta-local-opt} below implies
$V_\delta - J(\nu)=o(\sqrt{\delta})$, so the claim holds with the identity map.
We henceforth assume $K_1>0$. 
We will mimic the proof of \cite{BartlDrapeauOblojWiesel2021}.

\medskip
\noindent\textbf{Step 1: Coupling class $\mathcal{C}_\delta$ and $L^2$--scaling.}
For $\delta>0$, define the class of couplings
\[
  \mathcal{C}_\delta
  := \Big\{ \pi\in\cP(\R^d\times\R^d):
       \pi(\R^d,\cdot)=\nu,\;
       \E_\pi\big[c(X', X)\big]\le\delta
     \Big\}.
\]
Write $\Delta:=X'-X$ and note that for $\pi\in\mathcal{C}_\delta$,
\begin{equation}\label{eq:L2-budget}
  \E_\pi\|\Delta\|_q^2 \;\le\; \frac{\delta}{\tau},
  \qquad
  \|\Delta\|_{L^2(\pi)}
  := \big(\E_\pi\|\Delta\|_q^2\big)^{1/2}
  \;\le\; \sqrt{\frac{\delta}{\tau}}.
\end{equation}
In particular, if $\delta_n\downarrow 0$ and $\pi_n\in\mathcal{C}_{\delta_n}$, then
$\|\Delta_n\|_{L^2(\pi_n)}\to 0$. Since $\|\Delta_n\|_q^2\ge 0$ and
$\E\|\Delta_n\|_q^2\to 0$, the family $\{\|\Delta_n\|_q^2\}_n$ is uniformly integrable.

Moreover, for any $\mu$ with $D_c(\mu,\nu)\le\delta$,  there exists at least one coupling
$\pi\in\mathcal{C}_\delta$ with first marginal $\mu$ (an optimal transport plan);
conversely, every $\pi\in\mathcal{C}_\delta$
has some first marginal $\mu$ with $D_c(\mu,\nu)\le\delta$.
Thus, 
\[
  V_\delta - J(\nu)
  = \sup_{\substack{\mu:\,D_c(\mu,\nu)\le\delta}}
      \big(J(\mu)-J(\nu)\big)
  = \sup_{\pi\in\mathcal{C}_\delta}
      \big(J(\mu_\pi)-J(\nu)\big),
\]
where $\mu_\pi$ denotes the first marginal of $\pi$.

\medskip
\noindent\textbf{Step 2: Path identity, Taylor expansion, and linearization.}
Fix $\pi\in\mathcal{C}_\delta$ with first marginal $\mu$ and displacement $\Delta=X'-X$.
For $t\in[0,1]$, set $\nu_t := (1-t)\nu + t\mu$. By Assumption (i),
\begin{equation}\label{eq:path-identity}
  J(\mu)-J(\nu)
  = \int_0^1 \E_\pi\big[ J'_{\nu_t}(X') - J'_{\nu_t}(X) \big]\,dt.
\end{equation}
For each fixed $t$, the map $x\mapsto J'_{\nu_t}(x)$ is differentiable with gradient
$\nabla J'_{\nu_t}(x)$. Using the fundamental theorem of calculus along the segment
$X+s\Delta$, $s\in[0,1]$, we obtain
\[
  J'_{\nu_t}(X+\Delta) - J'_{\nu_t}(X)
  = \int_0^1 \nabla J'_{\nu_t}(X+s\Delta)\cdot\Delta\,ds.
\]
Subtract $\nabla J'_{\nu_t}(X)\cdot\Delta$ and define the remainder
\[
  R_{\delta,t}
  := \int_0^1 \big(\nabla J'_{\nu_t}(X+s\Delta)-\nabla J'_{\nu_t}(X)\big)\cdot\Delta\,ds.
\]
Then
\begin{equation}\label{eq:linearization-plus-R}
  J(\mu)-J(\nu)
  = \int_0^1 \E_\pi\big[\nabla J'_{\nu_t}(X)\cdot\Delta\big]\,dt
    \;+\; \int_0^1 \E_\pi\big[R_{\delta,t}\big]\,dt.
\end{equation}

We now show that the remainder term is $o(\sqrt{\delta})$, \emph{uniformly} over
$\pi\in\mathcal{C}_\delta$ and $t\in[0,1]$. By the Cauchy--Schwarz inequality in $L^2(\pi)$ and the H\"older inequality in $\R^d$,
\begin{align*}
  \E_\pi|R_{\delta,t}|
  &\le \int_0^1
       \E_\pi\big[
         \|\nabla J'_{\nu_t}(X+s\Delta)-\nabla J'_{\nu_t}(X)\|_p\,\|\Delta\|_q
       \big]\,ds\\
  &\le \int_0^1
       \big\|\nabla J'_{\nu_t}(X+s\Delta)-\nabla J'_{\nu_t}(X)\big\|_{L^2(\pi)}
       \,\|\Delta\|_{L^2(\pi)}\,ds.
\end{align*}
Using \eqref{eq:L2-budget},
\[
  \E_\pi|R_{\delta,t}|
  \le \sqrt{\frac{\delta}{\tau}}\,
       \sup_{s\in[0,1]}
       \big\|\nabla J'_{\nu_t}(X+s\Delta)-\nabla J'_{\nu_t}(X)\big\|_{L^2(\pi)}.
\]

\smallskip
\emph{Uniform $L^2$ bounds via polynomial growth.}
By Assumption (ii), there exist $1\le r\le 2$ and $C>0$ such that
\[
  \|\nabla J'_{\nu_t}(x)\|_p \le C\big(1+\|x\|^{r-1}\big),
  \qquad x\in\R^d,\ t\in[0,1],
\]
for all $\mu$ in a neighbourhood of $\nu$; in particular, this holds for all interpolants
$\nu_t=(1-t)\nu+t\mu$ with $D_c(\mu,\nu)$ small. Let $(X',X)\sim\pi\in\mathcal{C}_{\delta_0}$.
Then, by the triangle inequality and the fact that $2(r-1)\le 2$,
\[
  \|X+s\Delta\|^{2(r-1)}
  \le C'\big(\|X\|^{2(r-1)} + \|\Delta\|^{2(r-1)}\big),
  \qquad s\in[0,1],
\]
for some $C'>0$. Since $2(r-1)\le 2$ and $\E_\pi\|\Delta\|^2\le\delta_0/\tau$, Jensen's
inequality yields
\[
  \sup_{\pi\in\mathcal{C}_{\delta_0}}\E_\pi\|\Delta\|^{2(r-1)} < \infty,
\]
and $\E_\nu\|X\|^{2r}<\infty$ implies $\E_\nu\|X\|^{2(r-1)}<\infty$, which implies 
\[
  \sup_{\pi\in\mathcal{C}_{\delta_0}}
  \sup_{s\in[0,1]}
  \E_\pi\|X+s\Delta\|^{2(r-1)}
  <\infty.
\]
Using the growth bound
\[
  \|\nabla J'_{\nu_t}(X+s\Delta)\|_p^2
  \le C''\big(1+\|X+s\Delta\|^{2(r-1)}\big),
\]
we obtain
\[
  \sup_{\pi\in\mathcal{C}_{\delta_0}}
  \sup_{t\in[0,1],\,s\in[0,1]}
  \E_\pi\big[\|\nabla J'_{\nu_t}(X+s\Delta)\|_p^2\big]
  <\infty,
\]
and likewise
\[
  \sup_{t\in[0,1]}
  \E_\nu\big[\|\nabla J'_{\nu_t}(X)\|_p^2\big] <\infty.
\]
Thus the family $\{\nabla J'_{\nu_t}(X+s\Delta): \pi\in\mathcal{C}_{\delta_0}, t,s\in[0,1]\}$
is uniformly square-integrable. The uniformity follows by a standard sequential argument: for any
$\delta_n\downarrow0$ and $\pi_n\in\mathcal{C}_{\delta_n}$,
$\Delta_n\to0$ in $L^2$, and joint continuity together with the
polynomial growth bound and Vitali's theorem yields the required
uniform $L^2$ convergence.

\smallskip
\emph{Spatial term.}
For the spatial part, set
\[
  A_{\delta,t,s}
  := \big\|\nabla J'_{\nu_t}(X+s\Delta)-\nabla J'_{\nu_t}(X)\big\|_{L^2(\pi)}.
\]
By joint continuity of $(\mu,x)\mapsto\nabla J'_\mu(x)$ and the fact that
$\|\Delta\|_{L^2(\pi)}\to 0$ uniformly over $\pi\in\mathcal{C}_\delta$ as $\delta\downarrow 0$,
we have $\nabla J'_{\nu_t}(X+s\Delta)\to\nabla J'_{\nu_t}(X)$ in probability.
Together with the uniform $L^2$-bound just proved and uniform integrability of
$\|\nabla J'_{\nu_t}(X+s\Delta)\|_p^2$, Vitali's theorem yields
\[
  \sup_{\pi\in\mathcal{C}_\delta}
  \sup_{t\in[0,1],\,s\in[0,1]}
  A_{\delta,t,s}
  \;=\; o(1),
  \qquad \delta\downarrow 0.
\]

\smallskip
\emph{Measure term.}
For the measure part, consider
\[
  B_{\delta,t}
  := \big\|\nabla J'_{\nu_t}(X)-g(X)\big\|_{L^2(\pi)}.
\]
Since $X$ has law $\nu$ under any $\pi\in\mathcal{C}_\delta$, we have
\[
  B_{\delta,t}
  = \Big(\E_\nu\big[\|\nabla J'_{\nu_t}(X)-g(X)\|_p^2\big]\Big)^{1/2}.
\]
By joint continuity of $(\mu,x)\mapsto\nabla J'_\mu(x)$ and the uniform $L^2$-bound
from the growth condition, the map
$\mu\mapsto\nabla J'_\mu$ is continuous in $L^2(\nu)$ near $\nu$.
Moreover, for $\mu$ with $D_c(\mu,\nu)\le\delta_0$, all interpolants
$\nu_t=(1-t)\nu+t\mu$ also satisfy $D_c(\nu_t,\nu)\le\delta_0$, so the same
uniform bound applies. Since $\nu_t\to\nu$ weakly and the family is uniformly
square-integrable, Vitali's theorem yields
\[
  \sup_{\substack{\mu:D_c(\mu,\nu)\le\delta_0,\\ t\in[0,1]}}
  B_{\delta,t}
  = o(1),
\]
as $\mu\to\nu$ (hence $\nu_t\to\nu$ in the weak topology). Since here
we are only considering $\mu$ with $D_c(\mu,\nu)\le\delta$ and $\delta\downarrow 0$,
this implies
\[
  \sup_{\pi\in\mathcal{C}_\delta}\sup_{t\in[0,1]} B_{\delta,t}
  \;=\; o(1).
\]

\smallskip
Putting the spatial and measure parts together and using the triangle inequality, we obtain
\[
  \sup_{\pi\in\mathcal{C}_\delta}
  \sup_{t\in[0,1]}
  \sup_{s\in[0,1]}
  \big\|\nabla J'_{\nu_t}(X+s\Delta)-g(X)\big\|_{L^2(\pi)}
  = o(1),
  \qquad \delta\downarrow 0.
\]
In particular,
\[
  \sup_{\pi\in\mathcal{C}_\delta}
  \sup_{t\in[0,1]}
  \sup_{s\in[0,1]}
  \big\|\nabla J'_{\nu_t}(X+s\Delta)-\nabla J'_{\nu_t}(X)\big\|_{L^2(\pi)}
  = o(1),
  \qquad \delta\downarrow 0.
\]
Therefore,
\[
  \sup_{\pi\in\mathcal{C}_\delta}
  \sup_{t\in[0,1]}
  \E_\pi|R_{\delta,t}|
  \;\le\;
  \sqrt{\frac{\delta}{\tau}}\cdot o(1)
  = o(\sqrt{\delta}).
\]
Integrating over $t\in[0,1]$,
\[
  \sup_{\pi\in\mathcal{C}_\delta}
  \left|
    \int_0^1 \E_\pi[R_{\delta,t}]\,dt
  \right|
  = o(\sqrt{\delta}),
  \qquad \delta\downarrow 0.
\]

Next, by the same Vitali-type argument (now without the $s$-shift),
\[
  \int_0^1 \E_\pi\big[\nabla J'_{\nu_t}(X)\cdot\Delta\big]\,dt
  = \E_\pi[g(X)\cdot\Delta] + o(\sqrt{\delta}),
\]
uniformly over $\pi\in\mathcal{C}_\delta$. Combining this with
\eqref{eq:linearization-plus-R}, we conclude that
\begin{equation}\label{eq:almost-linear}
  J(\mu)-J(\nu)
  = \E_\pi[g(X)\cdot\Delta] + o(\sqrt{\delta}),
\end{equation}
with the $o(\sqrt{\delta})$ term uniform over $\pi\in\mathcal{C}_\delta$ (hence over all
feasible $\mu$ with $D_c(\mu,\nu)\le\delta$).

Taking the supremum over $\pi\in\mathcal{C}_\delta$,
\begin{equation}\label{eq:Vdelta-local-opt}
  V_\delta - J(\nu)
  = \sup_{\pi\in\mathcal{C}_\delta}
      \E_\pi[g(X)\cdot\Delta]
    \;+\; o(\sqrt{\delta}),
  \qquad \delta\downarrow 0.
\end{equation}

\medskip
\noindent\textbf{Step 3: Upper bound via H\"older and Cauchy--Schwarz.}
For any $\pi\in\mathcal{C}_\delta$,
\begin{align*}
  \E_\pi[g(X)\cdot\Delta]
  &\le \E_\pi\big[\|g(X)\|_p\,\|\Delta\|_q\big]\Big(\E_\nu\|g(X)\|_p^2\Big)^{1/2}
       \Big(\E_\pi\|\Delta\|_q^2\Big)^{1/2}\le K_1\,\sqrt{\frac{\delta}{\tau}},
\end{align*}
where we used \eqref{eq:L2-budget} and the fact that $X\sim\nu$ under $\pi$.
Thus,
\[
  \sup_{\pi\in\mathcal{C}_\delta}
    \E_\pi[g(X)\cdot\Delta]
  \le K_1\,\sqrt{\frac{\delta}{\tau}},
\]
and \eqref{eq:Vdelta-local-opt} yields
\[
  V_\delta - J(\nu)
  \le \sqrt{\frac{\delta}{\tau}}\,K_1 + o(\sqrt{\delta}).
\]

\medskip
\noindent\textbf{Step 4: Lower bound via a deterministic extremal Monge map.}
For each $x$ with $g(x)\neq 0$, choose
\[
  u(x) \in \arg\max_{\|u\|_q=1} g(x)\cdot u,
\]
and define
\[
  \bar{\Delta}(x) := \frac{\|g(x)\|_p}{K_1}\,u(x),
  \qquad
  \bar{\Delta}(x) := 0 \text{ if } g(x)=0.
\]
{
Note that we choose \(u(\cdot)\) to be a measurable selection. For the \(\ell_q\)-norm, one may take
\[
u_i(x)
=
\frac{\operatorname{sgn}(g_i(x))|g_i(x)|^{p-1}}
{\|g(x)\|_p^{p-1}},
\qquad g(x)\ne0,
\]
and set \(u(x)=0\) when \(g(x)=0\).
}
Then
\[
  \E_\nu\|\bar{\Delta}(X)\|_q^2
  = \frac{1}{K_1^2}\,\E_\nu\|g(X)\|_p^2 = 1,
\]
and
\[
  \E_\nu[g(X)\cdot\bar{\Delta}(X)]
  = \frac{1}{K_1}\,\E_\nu\big[\|g(X)\|_p\,g(X)\cdot u(X)\big]
  = \frac{1}{K_1}\,\E_\nu\|g(X)\|_p^2 = K_1.
\]

Define the Monge map
\[
  T_\delta(x) := x + \sqrt{\frac{\delta}{\tau}}\,\bar{\Delta}(x),
  \qquad
  \mu_\delta := T_{\delta\#}\nu,
\]
and let $\pi_\delta$ be the coupling $(X',X)=(T_\delta(X),X)$ with $X\sim\nu$.
Then
\[
  \E_{\pi_\delta}\|X'-X\|_q^2
  = \frac{\delta}{\tau}\,\E_\nu\|\bar{\Delta}(X)\|_q^2
  = \frac{\delta}{\tau},
\]
so
\[
  D_c(\mu_\delta,\nu)
  \le \E_{\pi_\delta} c(X', X)
  = \tau\,\E_{\pi_\delta}\|X'-X\|_q^2
  = \delta,
\]
and hence $\pi_\delta\in\mathcal{C}_\delta$ and $\mu_\delta$ is feasible.

Applying \eqref{eq:almost-linear} with $\pi_\delta$,
\[
  J(\mu_\delta)-J(\nu)
  = \E_{\pi_\delta}[g(X)\cdot\Delta] + o(\sqrt{\delta})
  = \sqrt{\frac{\delta}{\tau}}\,
      \E_\nu[g(X)\cdot\bar{\Delta}(X)]
      + o(\sqrt{\delta})
  = \sqrt{\frac{\delta}{\tau}}\,K_1 + o(\sqrt{\delta}).
\]
Therefore,
\[
  \liminf_{\delta\downarrow 0}
  \frac{V_\delta - J(\nu)}{\sqrt{\delta}}
  \ge \sqrt{\frac{1}{\tau}}\,K_1.
\]

\medskip
\noindent\textbf{Step 5: Conclusion.}
Combining the upper bound from Step~3 with the lower bound from Step~4,
\[
  \sqrt{\frac{1}{\tau}}\,K_1
  \le \liminf_{\delta\downarrow 0}
     \frac{V_\delta - J(\nu)}{\sqrt{\delta}}
  \le \limsup_{\delta\downarrow 0}
     \frac{V_\delta - J(\nu)}{\sqrt{\delta}}
  \le \sqrt{\frac{1}{\tau}}\,K_1.
\]
Hence
\[
  V_\delta
  = J(\nu) + \sqrt{\frac{\delta}{\tau}}\,K_1 + o(\sqrt{\delta})
  = J(\nu)
    + \sqrt{\frac{\delta}{\tau}}\,
      \Big(\E_\nu\|g(X)\|_p^2\Big)^{1/2}
    + o(\sqrt{\delta}),
\]
and the Monge map $T_\delta$ constructed in Step~4 is asymptotically optimal
and of the stated form.
\end{proof}

When applied to the specific Problem \eqref{FWProblem1new}, the proof reduces to explicitly computing the Wasserstein derivative and checking the regularity conditions.

\begin{corollary}[Asymptotic non-linear optimal perturbation under Wasserstein distance]
\label{thm:asymptotic_perturbation_cor}
Suppose Assumptions \ref{ass:utility} and \ref{ass:compact_B} hold.
Let $K$ be the compact set in Assumption~\ref{ass:compact_B}, and assume that $\operatorname{dist}\!\left(
\operatorname{supp}(\mathbb{P}_0),K^c
\right)>0.$ 
Define the functional
\[
J(\mu)
=
\int_{\mathbb R^d}
F_\mu(T,z)
U\!\left(
I\!\left(
\frac{\mathcal K_\mu(x_0)e^{-rT}}{F_\mu(T,z)}
\right)
\right)
\varphi_T(z)\,dz,  \qquad \mu \in \{{\mu \in \mathcal{P}(K)}: D_c(\mu,\mathbb{P}_0) \le \delta\},
\]
where $F_\mu(T,z)= \mathbb E_\mu[L_T(B,z)]$. 
Define the influence function $H: \R^d \to \R^d$ by
\[
H(b)
:= \nabla J'_{\mathbb{P}_0}(b)
= \int_{\R^{d}} \nabla_b L_T(b,z)\,
U\!\left(
I\!\left(
\frac{\mathcal{K}_{\mathbb{P}_0}(x_0)e^{-rT}}
{F_{\mathbb{P}_0}(T,z)}
\right)
\right)
\varphi_T(z)\,dz,
\]
where $\mathcal{K}_{\mathbb{P}_0}(x_0)$ is the unique solution to the budget constraint under $\mathbb{P}_0$. Let 
\[
\|H\|_{L^2(\mathbb{P}_0)} 
:= \left( \E_{\mathbb{P}_0}\big[\|H(B)\|_2^2\big] \right)^{1/2}.
\]
Assume that $\|H\|_{L^2(\mathbb{P}_0)}>0$.
Then, as $\delta \to 0$, an asymptotically optimal adversarial perturbation is given by the pushforward $\mu^{\text{WC}}_\delta = (\text{Id} + \Delta^{\text{WC}}_\delta)_{\#}\mathbb{P}_0$, where
\[
\Delta_\delta^{\text{WC}}(b)
=
-\sqrt{\delta}\,
\frac{H(b)}{\|H\|_{L^2(\mathbb{P}_0)}}
+ o(\sqrt{\delta})
\quad\text{in }L^2(\mathbb{P}_0).
\]
Furthermore, the corresponding asymptotically optimal value is
\[
\inf_{{\mu \in \mathcal{P}(K)}: D_c(\mu,\mathbb{P}_0) \le \delta} J(\mu)
= J(\mathbb{P}_0) - \sqrt{\delta}\, \|H\|_{L^2(\mathbb{P}_0)} + o(\sqrt{\delta}).
\]
\end{corollary}

\begin{proof}
We show that $J$ satisfies the assumptions of Theorem~\ref{thm:asymptotic_perturbation} 
with $q=p=2$, compute the Wasserstein gradient at $\mathbb{P}_0$, and then apply the theorem
to the functional $-J$ to obtain the stated expansion.

\medskip
\noindent\textbf{Step 1: Computation of the Wasserstein derivative at $\mathbb{P}_0$.}
Consider the perturbation $\mu^\epsilon=(1-\epsilon)\mathbb{P}_0+\epsilon\delta_b$ and write
$k^\epsilon=\mathcal{K}_{\mu^\epsilon}(x_0)$, $F^\epsilon(z)=F_{\mu^\epsilon}(T,z)$, and
\[
J(\mu^\epsilon) = \int_{\R^d} F^\epsilon(z)\,
U\!\left(I\!\left(\frac{k^\epsilon e^{-rT}}{F^\epsilon(z)}\right)\right)\varphi_T(z)\,dz.
\]
Since $F^\epsilon(z)=(1-\epsilon)F_0(z)+\epsilon L_T(b,z)$ and all integrands are dominated by 
$\exp(c\|z\|)\varphi_T(z)$ (by Assumption~\ref{ass:utility} and compactness of $B$),
the dominated convergence theorem applies, allowing differentiation under the integral.

Differentiating $J(\mu^{\epsilon})$ at $\epsilon=0$, and using
$U'(I(y))=y$, $I'(y)=1/U''(I(y))$, together with the differentiated
budget constraint, the terms involving
$\dot{k}=\left.\frac{dk^\epsilon}{d\epsilon}\right|_{\epsilon=0}$
and $I'$ cancel. Hence, up to an additive constant, one obtains the
first variation
\[
\frac{\delta J}{\delta \mathbb{P}_0}(b)
=
\int_{\R^d}
L_T(b,z)\,
U\!\left(
I\!\left(
\frac{k_0e^{-rT}}{F_0(z)}
\right)
\right)
\varphi_T(z)\,dz.
\]
Taking the spatial gradient yields the Wasserstein gradient
\begin{equation}\label{eq:H_final_short}
\nabla J'_{\mathbb{P}_0}(b)
=
\int_{\R^d}
\nabla_bL_T(b,z)\,
U\!\left(
I\!\left(
\frac{\mathcal{K}_{\mathbb{P}_0}(x_0)e^{-rT}}
{F_{\mathbb{P}_0}(T,z)}
\right)
\right)
\varphi_T(z)\,dz
=:H(b),
\end{equation}
matching exactly the influence function stated in the corollary. 

\medskip
\noindent\textbf{Step 2: Verifying regularity and applying Theorem~\ref{thm:asymptotic_perturbation}.}
All remaining assumptions of Theorem~\ref{thm:asymptotic_perturbation} follow immediately
from the compact support of $B$ (Assumption~\ref{ass:compact_B}) and a localization argument, the exponential bounds
on $L_T$ and $\nabla_b L_T$ from Assumption~\ref{ass:utility}, and the continuity of
$\mu \mapsto F_\mu$ and $\mu \mapsto \mathcal{K}_\mu(x_0)$ under weak convergence.  
In particular:

(i) $H\in L^2(\mathbb{P}_0)$ by the uniform exponential bound in~\eqref{eq:H_final_short};  
(ii) $(\mu, b)\mapsto \nabla J'_\mu(b)$ is jointly continuous, again by dominated convergence;  
(iii) $J$ is Gâteaux differentiable along quadratic-cost interpolations, so the linearization formula required by the theorem holds.

Thus, the hypotheses of Theorem~\ref{thm:asymptotic_perturbation} hold with $q=p=2$ and $r=1$.

Applying the theorem to $\widetilde J(\mu)=-J(\mu)$ (so that $\nabla\widetilde J'_{\mathbb{P}_0}=-H$)
gives, as $\delta\to0$,
\[
\sup_{{\mu \in \mathcal{P}(K)}: D_c(\mu,\mathbb{P}_0)\le\delta}\widetilde J(\mu)
= \widetilde J(\mathbb{P}_0) + \sqrt{\delta}\,\|H\|_{L^2(\mathbb{P}_0)} + o(\sqrt{\delta}),
\]
which is equivalent to
\[
\inf_{{\mu \in \mathcal{P}(K)}: D_c(\mu,\mathbb{P}_0)\le\delta}J(\mu)
= J(\mathbb{P}_0) - \sqrt{\delta}\,\|H\|_{L^2(\mathbb{P}_0)} + o(\sqrt{\delta}).
\]
{
The leading-order Monge perturbation is
\[
\Delta_\delta(b)
=
-\sqrt{\delta}\frac{H(b)}{\|H\|_{L^2(\mathbb P_0)}}.
\]
More generally, any perturbation equal to this plus \(o(\sqrt\delta)\) in \(L^2(\mathbb P_0)\)
is asymptotically optimal, provided the resulting pushforward remains supported in \(K\). 
Since $H$ is continuous on the compact set $K$, it is bounded.
Therefore, under
$\operatorname{dist}(\operatorname{supp}(\mathbb P_0),K^c)>0$,
the map $b\mapsto b+\Delta_\delta(b)$ maps
$\operatorname{supp}(P_0)$ into $K$ for all sufficiently small
$\delta$.} This completes the proof. 
\end{proof}

If $\mathbb{P}_0 = \mathbb{P}_n$ for some empirical measure, then we can compute $\Delta^*$ by replacing $\mathbb{P}_0$ by $\mathbb{P}_n$ conditioned on these samples, rather than viewing them as random measures. We are going to construct i.i.d. samples of $B^{(k)}$ in Section \ref{projection}, and then we construct the optimal measure for the original problem via the %
empirical measure 
\[
\mathbb P_n(dx) \;=\; \frac{1}{n}\sum_{k=1}^n \delta_{C^{(k)}}(dx),
\]where for each $k$, $C^{(k)} = B^{(k)} + \Delta^*$ { and $\delta_{C^{(k)}}$ denotes the Dirac measure at $ C^{(k)}$}. Therefore,  this gives a way of solving the Wasserstein constrained distributional optimization problem where {the radius} $\delta > 0$ is small.

\section{Data-Driven Formulation and Choice of Model Parameters}\label{projection}
In this section, we first describe the data-driven version of the DRBC formulation of the stochastic control problem since it is natural to inform the choice of the prior $\mathbb{P}_0$ from the data, and an appropriate empirical measure is the natural candidate. Next, we provide the prescription to choose an uncertainty radius $\delta$ using an asymptotically optimal (as data collected increases) approach. The idea is that this choice should also be based on observed data. We calibrate the distributional ambiguity set to the smallest size that makes the oracle-optimal portfolio (i.e., the one chosen if the true distribution were known) statistically plausible at the desired confidence level. The prescription is then obtained from an asymptotic statistical result for Wasserstein projection onto a nonlinear manifold.

Recall that we assume each asset follows a geometric Brownian motion with random drift:
\[
\frac{dS_i(t)}{S_i(t)} \;=\; B_i\,dt \;+\; \sum_{j=1}^d \sigma_{ij}\,dW_j(t),
\qquad i=1,\dots,d.
\]
Therefore, we have the explicit closed-form formula
\[
\log S_i(t)-\log S_i(0)
= \Big(B_i-\tfrac12\|\sigma_{i\cdot}\|^2\Big)\,t
\;+\;\sum_{j=1}^d \sigma_{ij}\,W_j(t),
\]
where $\|\sigma_{i\cdot}\|^2:=\sum_{j=1}^d \sigma_{ij}^2$.

We adopt a data-driven approximation in which the ground-truth distribution of $B$ is fixed, but at the beginning of the $k$-th non-overlapping window of length $\tilde t$, the realization $B^{(k)}$ of return vectors is redrawn i.i.d. from such ground-truth distribution and constant during the window. Within the $k$-th window, we observe prices of the stock $S$ on a grid $\{0,h,2h,\dots,mh\}$ and we let $\tilde t = mh$. 

Asset $i$ in window $k$ has return $B_i^{(k)}$. An unbiased estimator of $B_i^{(k)}$ (conditional on $B^{(k)}$) is obtained from the endpoint log-return:
\[
\widehat B^{(k)}_i
\;=\;
\frac{\log S_i^{(k)}(\tilde t)-\log S_i^{(k)}(0)}{\tilde t}
\;+\;\tfrac12\|\sigma_{i\cdot}\|^2,
\qquad
\mathbb E\!\left[\widehat B^{(k)}_i \,\middle|\, B^{(k)}\right]=B^{(k)}_i.
\]
Equivalently, averaging one-step log-returns over the window yields
\[
\widehat B^{(k)}_i
\;=\;
\frac{1}{mh}\sum_{\ell=0}^{m-1}
\log\!\left(\frac{S_i^{(k)}((\ell+1)h)}{S_i^{(k)}(\ell h )}\right)
\;+\;\tfrac12\|\sigma_{i\cdot}\|^2.
\]

Collecting  $k \in \{1,...,n\}$ windows yields i.i.d. estimates $ B^{(1)},\dots,B^{(n)}\in\mathbb R^d$. %
We slightly abuse notation: the quantities obtained above are noisy estimates of the realized drifts \(B^{(k)}\), not the realized drifts themselves. 
Strictly speaking, we should apply a deconvolution method, which we did but the results did not change significantly and this approach is much easier to implement. 
{The distributional robustness is intended to absorb the remaining estimation noise in
\(\widehat B^{(k)}\) and other errors arising from imperfect drift estimation.} So, ultimately, we take the nominal prior as the empirical measure
\[
\mathbb P_n(dx) \;=\; \frac{1}{n}\sum_{k=1}^n \delta_{B^{(k)}}(dx),
\]
{where $\delta_{B^{(k)}}$ in the equation denotes the Dirac measure at $ B^{(k)}$.}

The choice of the key parameter $\delta$ is crucial. If $\delta$ is too large, there is too much model ambiguity, and the available data becomes less relevant. If $\delta$ is too small, the effect of robustification is negligible. Therefore, the choice of $\delta$ should not be exogenously defined; rather, it should be endogenously informed by the data. Before presenting the methodology, we first present some technical assumptions. We denote $k^* = \mathcal{K}(x_0)$ as the optimal Lagrangian multiplier in Eq. \eqref{eq:budget}.

In order to choose an appropriate $\delta = \delta_n$, here we follow the idea behind the RWPI approach introduced in \cite{BKW19}. Intuitively, $\delta$ should be chosen such that the set $\mathcal U^{\text{OT}}_{\delta,B}(\mathbb P_n)$ contains all the probability measures that are plausible variations of the data represented by $\mathbb P_n$. 
{According to Theorem \ref{thm:karatzas_zhao_general_compact}, under the ground-truth prior $\mu^*$, the (ex ante) optimal Bayesian policy $\pi^*$ is uniquely characterized by the Lagrange multiplier $k^*=\mathcal{K}_{\mu^*}(x_0)$. In the terminology introduced in Section \ref{intro}, we interpret this optimal policy in a pre-committed (i.e., ex ante) sense. %
Denoting the ground-truth product measure by $\mathbb{Q}^* = \mu^* \otimes \mathbb{P}_W$,} We restate here (real optimal policy $\pi^*$ and $X^*(T)$ has one-to-one correspondence)
$$
 \mathbb{E}_{\mathbb{Q}^*}\left[\hat{\Lambda}(T)X^*(T)\right] =  \int_{\mathbb{R}^d }I\left(\frac{k^*e^{-rT}}{\mathbb{E}_{\mu^*}\left[L_T(B,y)\right]}\right)\varphi_T(y)dy = x_0e^{rT}. %
$$
We can see that the optimal policy $\pi^*$ and the Lagrangian multiplier $k^*$ also have a one-to-one correspondence. Thus, without loss of generality, we assume $k > 0$ is the decision variable and choose the optimal $\delta$ based on whether the optimal $k^*$ is covered. Similar to the notation in \cite{BKW19}, 
{for \(k>0\), we define the set of admissible drift priors
\[
\mathcal M_k
:=
\left\{
\mu \in\mathcal P(K):
\int_{\mathbb R^d}
I\left(
\frac{k e^{-rT}}{\mathbb{E}_\mu[L_T(B,y)]}
\right)\varphi_T(y)\,dy
=
x_0e^{rT}
\right\}.
\]}
{
Now, for a fixed $\mu \in \mathcal{U}^{\text{OT}}_{\delta,B}(\mathbb P_n)$, the set $\left\{k > 0: 
\mu \in \mathcal{M}_k\right\}$ contains all parameter choices that are optimal from the decision maker's point of view. This motivates the definition of the following set
\begin{align*}
\Lambda_{\delta}(\mathbb P_n) &= \left\{k > 0: \mathcal{M}_k \cap \mathcal{U}^{\text{OT}}_{\delta,B}(\mathbb P_n) \neq \varnothing \right\}\\
&=\left\{k > 0: \text{there exists } \mu \in \mathcal{U}^{\text{OT}}_{\delta,B}(\mathbb P_n) \text{ such that }\int_{\mathbb{R}^d }I\left(\frac{ke^{-rT}}{\mathbb{E}_{\mu}\left[L_T(B,y)\right]}\right)\varphi_T(y)dy = x_0e^{rT} \right\},   
\end{align*}
}
which corresponds to all the plausible estimates of $k^*$. Thus, it includes all the parameter choices that are collected by the decision maker as optimal for some distribution in the uncertainty set. As a result, $\Lambda_{\delta}(\mathbb P_n)$ is a natural confidence region for $k^*$. Therefore, $\delta > 0$ should be chosen as the smallest $\delta^*_n$ such that $k^*$ belongs to this region with a given confidence interval. Namely, {let \(\mathbb P_{\mu^*}\) denote the law of the i.i.d. sample
\(B^{(1)},\ldots,B^{(n)}\sim\mu^*\) and then choose
\[
\delta_n^*
=
\min\left\{
\delta>0:
\mathbb P_{\mu^*}(k^*\in\Lambda_\delta(\mathbb P_n))\ge 1-\delta_0
\right\},
\]} 
where %
$\delta_0$ is the user-defined confidence level (typically $\delta_0 = 0.05$).

However, from the definition alone, it is hard to compute $\delta^*_n$. We now provide a simpler representation for $\delta^*_n$ via an auxiliary function called the robust Wasserstein profile (RWP) function:
\begin{align*}
R_n(k)&:=\inf_{\mu \in \mathcal{M}_{k}}D_c(\mathbb{P}_n, \mu) \\
&=\inf\left\{D_c(\mathbb{P}_n, \mu): \int_{\mathbb{R}^d }I\left(\frac{ke^{-rT}}{\mathbb{E}_{\mu}\left[L_T(B,y)\right]}\right)\varphi_T(y)dy = x_0e^{rT}\right\}.    
\end{align*}
Compared with the linear projection in \cite{BKW19}, here the RWP function is defined as a projection of the empirical measure to a nonlinear manifold, so the behavior of the nonlinear RWP function is more complicated. 

In the following theorem, we derive the nonlinear projection asymptotics with great generality, and the specific choice of the uncertainty radius can be obtained by a corollary.

\begin{theorem}[Non-linear projection asymptotics]
\label{thm:projection_law}
Let \(\nu\in \mathcal P_2(\mathbb R^d)\), and let
\(\kappa : \mathcal P_2(\mathbb R^d)\to \mathbb R\) be a functional with
\(\kappa(\nu)=0\).
Assume that \(\kappa\) admits a first variation \(\kappa'_\mu:\mathbb R^d\to\mathbb R\)
such that:

\begin{itemize}
    \item For any probability measure \(\mu\) in a \(W_2\)-neighborhood of \(\nu\),
    and any coupling \((X',X)\sim(\mu,\nu)\), with the linear interpolation $\nu_t := (1-t)\nu + t\mu,$
    we have
    \[
        \kappa(\mu)-\kappa(\nu)
        =
        \int_0^1
        \E\!\left[\kappa'_{\nu_t}(X')-\kappa'_{\nu_t}(X)\right]\,dt.
    \]

    \item For \(\mu\) in a \(W_2\)-neighborhood of \(\nu\), the map
    \(x\mapsto \kappa'_\mu(x)\) is \(C^1\), and its gradient $g_\mu(x):=\nabla_x \kappa'_\mu(x)$
    is jointly continuous in \((\mu,x)\) near \((\nu,\cdot)\).
    {We assume that $\E_\nu[\|g_\nu(X)\|_2^2]>0.$ Moreover, there exist constants \(L_\mu,L_x>0\) and a \(W_2\)-neighborhood
    \(\mathcal U\) of \(\nu\) such that for all \(\mu,\tilde\mu\in\mathcal U\)
    and all \(x,x'\in \mathbb{R}^d\),
    \[
        \|g_\mu(x)-g_{\tilde\mu}(x')\|_2
        \le
L_\mu\,W_2(\mu,\tilde\mu)\,\bigl(1+\|x\|_2+\|x'\|_2\bigr)
        +
        L_x\,\|x-x'\|_2.
    \]}
\end{itemize}
Let \(c:\mathbb R^d\times\mathbb R^d\to[0,\infty]\) be a transport cost satisfying:

\begin{itemize}
\item[(i)] There exist \(\tau>0\), \(r_0>0\), and a function
\(\eta:(0,r_0]\to[0,\infty)\) with \(\eta(r)\to0\) as \(r\to0\) such that
\begin{equation}
\label{eq:uniform_quad_c}
    \bigl|c(x,y)-\tau\|x-y\|_2^2\bigr|
    \le
    \eta(\|x-y\|_2)\,\|x-y\|_2^2
    \qquad\text{whenever }\|x-y\|_2\le r_0.
\end{equation}

{\item[(ii)] The cost admits the global quadratic lower bound
\[
    c(x,y)\ge \tau\|x-y\|_2^2,
    \qquad x,y\in\mathbb R^d.
\]}
\end{itemize}

Let
\[
    D_c(\nu,\mu)
    :=
    \inf_{\pi\in\Pi(\nu,\mu)}
    \int_{\mathbb R^d\times\mathbb R^d} c(x,y)\,\pi(dx,dy)
\]
denote the optimal transport divergence induced by \(c\). For \(z\) in a neighborhood
of \(0\), define the projection cost
\[
    R(z)
    :=
    \inf_{\mu:\kappa(\mu)=z} D_c(\nu,\mu).
\]

Then, as \(z\to0\),
\[
    R(z)
    =
    \frac{\tau z^2}{\E_\nu[\|g_\nu(X)\|_2^2]}
    +o(z^2)
    =
    \frac{\tau z^2}{\E_\nu[\|\nabla\kappa'_\nu(X)\|_2^2]}
    +o(z^2).
\]

Moreover, there exists an asymptotically optimal Monge-type perturbation
of the form
\[
    T_z(x)=x+\Delta_z(x),
    \qquad
    \Delta_z(x)
    =
    \frac{z}{\E_\nu[\|g_\nu(X)\|_2^2]}\,g_\nu(x)+o(z)
    \quad\text{in }L^2(\nu),
\]
which satisfies \(\kappa(T_{z\#}\nu)=z+o(z)\) and attains the above
cost up to \(o(z^2)\).
\end{theorem}

\begin{proof}
{Set $m:=\E_\nu\!\left[\|g_\nu(X)\|_2^2\right],$ where $X\sim \nu.$
By the assumptions given in the theorem, \(m>0\). Moreover, taking
\(\mu=\tilde\mu=\nu\) and \(x'=0\) in the Lipschitz bound gives $\|g_\nu(x)\|_2\le \|g_\nu(0)\|_2 + L_x \|x\|_2$. 
Since \(\nu\in \mathcal P_2(\mathbb R^d)\), this implies
\(\E_\nu[\|g_\nu(X)\|_2^2]<\infty\).
Since \(\mathcal U\) is a \(W_2\)-neighborhood of \(\nu\), there exists \(\rho_0>0\) such that $\{\rho\in\mathcal P_2(\mathbb R^d):W_2(\rho,\nu)<\rho_0\}\subset \mathcal U.$
Fix such a \(\rho_0\).}

\medskip
\noindent
{\textbf{Step 1: A local first-order expansion around \(\nu\).}
Let \(\Delta\in L^2(\nu;\mathbb R^d)\), let $\mu:=(\text{Id}+\Delta)_\#\nu,$
and assume \(\mu\in \mathcal U\). Let \(X\sim \nu\), and set \(X':=X+\Delta(X)\), so
\(X'\sim \mu\). For \(t\in[0,1]\), define \(\nu_t:=(1-t)\nu+t\mu\). By the first-variation
identity and the fundamental theorem of calculus,
\begin{align*}
\kappa(\mu)-\kappa(\nu)
&=
\int_0^1 \E\!\left[\kappa'_{\nu_t}(X')-\kappa'_{\nu_t}(X)\right]\,dt \\
&=
\int_0^1\int_0^1
\E\!\left[\left\langle g_{\nu_t}(X+s\Delta(X)),\Delta(X)\right\rangle\right]\,ds\,dt.
\end{align*}
Hence,
\[
\kappa(\mu)-\kappa(\nu)
=
\E_\nu\!\left[\langle g_\nu(X),\Delta(X)\rangle\right] + \mathcal R(\Delta),
\]
where
\[
\mathcal R(\Delta)
:=
\int_0^1\int_0^1
\E_\nu\!\left[
\left\langle g_{\nu_t}(X+s\Delta(X))-g_\nu(X),\Delta(X)\right\rangle
\right]ds\,dt.
\]
By the Cauchy--Schwarz inequality,
\[
|\mathcal R(\Delta)|
\le
\sup_{s,t\in[0,1]}
\left\|g_{\nu_t}(X+s\Delta(X))-g_\nu(X)\right\|_{L^2(\nu)}
\,
\|\Delta\|_{L^2(\nu)}.
\]
Now $W_2(\mu,\nu)\le \|\Delta\|_{L^2(\nu)}$. 
By coupling \(X\) with \(X+\xi \Delta(X)\), where \(\xi\sim \mathrm{Bernoulli}(t)\)
is independent of \(X\),
\[
W_2(\nu_t,\nu)^2\le t\,\|\Delta\|_{L^2(\nu)}^2
\le \|\Delta\|_{L^2(\nu)}^2.
\]
Therefore, by the Lipschitz assumption, we have 
\begin{align*}
&\left\|g_{\nu_t}(X+s\Delta(X))-g_\nu(X)\right\|_{L^2(\nu)} \\
&\qquad\le
L_\mu\,W_2(\nu_t,\nu)\,
\left\|1+\|X+s\Delta(X)\|_2+\|X\|_2\right\|_{L^2(\nu)}
+
L_x\,\|\Delta\|_{L^2(\nu)}.
\end{align*}
Since
\[
\|X+s\Delta(X)\|_{L^2(\nu)}
\le
\|X\|_{L^2(\nu)}+\|\Delta\|_{L^2(\nu)},
\]
it follows that, for all \(\|\Delta\|_{L^2(\nu)}\) sufficiently small,
\[
\sup_{s,t\in[0,1]}
\left\|g_{\nu_t}(X+s\Delta(X))-g_\nu(X)\right\|_{L^2(\nu)}
\le C\,\|\Delta\|_{L^2(\nu)}
\]
for some constant \(C>0\) depending only on \(\nu,L_\mu,L_x,\rho_0\). Hence
\begin{equation}
\label{eq:local_linearization_kappa}
\kappa(\mu)-\kappa(\nu)
=
\E_\nu\!\left[\langle g_\nu(X),\Delta(X)\rangle\right]
+
O\!\left(\|\Delta\|_{L^2(\nu)}^2\right).
\end{equation}
\medskip
\noindent
\textbf{Step 2: Upper bound via a bounded Monge perturbation.}
For \(M>0\), define the truncation
\[
T_M(v):=
v\,\min\!\left\{1,\frac{M}{\|v\|_2}\right\},
\qquad v\in \mathbb R^d,
\]
with the convention \(T_M(0)=0\). Choose any function \(M(z)\uparrow\infty\) as \(z\to0\)
such that $|z|\,M(z)\to0$ as \(z\to0\).
For instance, \(M(z)=|z|^{-1/2}\) for \(z\neq 0\).
Set
\[
h_z(x):=T_{M(z)}(g_\nu(x)),
\qquad
a_z:=\E_\nu\!\left[\langle g_\nu(X),h_z(X)\rangle\right].
\]
Since \(h_z\to g_\nu\) in \(L^2(\nu)\) as $z \to 0$, we have
\[
a_z\to m,
\qquad
\E_\nu[\|h_z(X)\|_2^2]\to m.
\]
Define
\[
\widetilde\Delta_z(x):=\frac{z}{a_z}\,h_z(x),
\qquad
\widetilde\mu_z:=(\text{Id}+\widetilde\Delta_z)_\#\nu.
\]
Then
\[
\|\widetilde\Delta_z\|_{L^2(\nu)}=O(|z|),
\qquad
\|\widetilde\Delta_z\|_{L^\infty}
\le
\frac{|z|\,M(z)}{|a_z|}
\to0.
\]
In particular, for all sufficiently small \(|z|\), we have $\widetilde\mu_z\in\mathcal U$ and $\|\widetilde\Delta_z\|_{L^\infty}\le r_0.$
Applying \eqref{eq:local_linearization_kappa} with \(\Delta=\widetilde\Delta_z\), we have 
\[
\kappa(\widetilde\mu_z)-\kappa(\nu)
=
\E_\nu\!\left[\langle g_\nu(X),\widetilde\Delta_z(X)\rangle\right]
+
O(z^2)
=
z+O(z^2).
\]
Now consider
\[
\Phi_z(b):=\kappa\!\left((\text{Id}+b\,\widetilde\Delta_z)_\#\nu\right),
\qquad b\in \left[\frac12,\frac32\right].
\]
By the same argument as above, uniformly for \(b\in [1/2,3/2]\), we have 
\[
\Phi_z(b)-\kappa(\nu)
=
b\,z+O(z^2).
\]
Hence, with $\psi_z(b):=\Phi_z(b)-z,$
we have $\psi_z(b)=(b-1)z+O(z^2)$
uniformly for \(b\in [1/2,3/2]\). Choosing \(\Gamma>0\) sufficiently large, for all
small \(|z|\), the two numbers $\psi_z(1-\Gamma |z|)$ and $\psi_z(1+\Gamma |z|)$
have opposite signs. Since \(b\mapsto \Phi_z(b)\) is continuous, the intermediate value theorem
yields some \(b_z\in [1-\Gamma|z|,1+\Gamma|z|]\) such that $\Phi_z(b_z)=z.$
Define
\[
\Delta_z:=b_z\,\widetilde\Delta_z,
\qquad
\mu_z:=(\text{Id}+\Delta_z)_\#\nu.
\]
Then
\[
\kappa(\mu_z)=z,
\qquad
b_z=1+O(|z|),
\qquad
\|\Delta_z\|_{L^\infty}\to0,
\qquad
\|\Delta_z\|_{L^2(\nu)}=O(|z|).
\]
Moreover,
\[
\Delta_z
=
\frac{z}{m}\,g_\nu + o(z)
\qquad\text{in }L^2(\nu),
\]
because \(b_z\to1\), \(a_z\to m\), and \(h_z\to g_\nu\) in \(L^2(\nu)\).
Since \(\|\Delta_z\|_{L^\infty}\le r_0\) for all sufficiently small \(|z|\), the local quadratic
expansion of the cost gives
\[
c(X,X+\Delta_z(X))
=
\tau \|\Delta_z(X)\|_2^2 + r(\Delta_z(X)),
\]
with
\[
|r(h)|\le \eta(\|h\|_2)\,\|h\|_2^2
\qquad\text{for }\|h\|_2\le r_0.
\]
Therefore, we have 
\[
D_c(\nu,\mu_z)
\le
\E_\nu\!\left[c(X,X+\Delta_z(X))\right]
=
\tau \E_\nu[\|\Delta_z(X)\|_2^2]
+
\E_\nu[r(\Delta_z(X))].
\]
Now
\[
\left|\E_\nu[r(\Delta_z(X))]\right|
\le
\eta(\|\Delta_z\|_{L^\infty})\,
\E_\nu[\|\Delta_z(X)\|_2^2]
=
o(z^2),
\]
because \(\|\Delta_z\|_{L^\infty}\to0\) and
\(\E_\nu[\|\Delta_z(X)\|_2^2]=O(z^2)\). Also,
\[
\E_\nu[\|\Delta_z(X)\|_2^2]
=
\frac{z^2}{m}+o(z^2),
\]
since
\[
\E_\nu[\|\widetilde\Delta_z(X)\|_2^2]
=
\frac{z^2}{a_z^2}\,\E_\nu[\|h_z(X)\|_2^2]
=
\frac{z^2}{m}+o(z^2),
\]
and \(b_z=1+O(|z|)\). Hence, we have
\begin{equation}
\label{eq:upper_bound_Rz}
R(z)\le D_c(\nu,\mu_z)
\le
\frac{\tau z^2}{m}+o(z^2).
\end{equation}
\medskip
\noindent
\textbf{Step 3: Lower bound along near minimizers.}
Let \(z_n\to0\). Choose \(\mu_n\in \mathcal P_2(\mathbb R^d)\) such that $\kappa(\mu_n)=z_n$ and $D_c(\nu,\mu_n)\le R(z_n)+o(z_n^2).$
By \eqref{eq:upper_bound_Rz}, $R(z_n)=O(z_n^2)$; 
hence $D_c(\nu,\mu_n)=O(z_n^2).$
Therefore, $D_c(\nu,\mu_n)\ge \tau W_2(\nu,\mu_n)^2,$
so $W_2(\nu,\mu_n)=O(|z_n|).$
In particular, for all large \(n\), we have \(\mu_n\in \mathcal U\).
Let \((X_n',X_n)\) be a coupling of \((\mu_n,\nu)\) such that $\E[\|X_n'-X_n\|_2^2]=W_2(\mu_n,\nu)^2,$
and set $D_n:=X_n'-X_n.$
Then
\[
\|D_n\|_{L^2}=W_2(\mu_n,\nu)=O(|z_n|).
\]
Let \(\nu_{n,t}:=(1-t)\nu+t\mu_n\). By the first-variation identity and the fundamental theorem
of calculus, we have 
\[
z_n
=
\int_0^1\int_0^1
\E\!\left[\left\langle g_{\nu_{n,t}}(X_n+sD_n),D_n\right\rangle\right]\,ds\,dt.
\]
Therefore
\[
z_n
=
\E\!\left[\langle g_\nu(X_n),D_n\rangle\right] + r_n,
\]
where
\[
r_n
:=
\int_0^1\int_0^1
\E\!\left[
\left\langle g_{\nu_{n,t}}(X_n+sD_n)-g_\nu(X_n),D_n\right\rangle
\right]ds\,dt.
\]
By the Cauchy--Schwarz inequality,
\[
|r_n|
\le
\sup_{s,t\in[0,1]}
\left\|g_{\nu_{n,t}}(X_n+sD_n)-g_\nu(X_n)\right\|_{L^2}
\,
\|D_n\|_{L^2}.
\]
Similar to Step 1, we have 
\[
W_2(\nu_{n,t},\nu)\le W_2(\mu_n,\nu)=O(|z_n|),
\]
and
\[
\|X_n+sD_n\|_{L^2}
\le
\|X_n\|_{L^2}+\|D_n\|_{L^2}
=
\|X\|_{L^2}+O(|z_n|).
\]
Hence, the Lipschitz assumption yields
\[
\sup_{s,t\in[0,1]}
\left\|g_{\nu_{n,t}}(X_n+sD_n)-g_\nu(X_n)\right\|_{L^2}
\le
C\,\|D_n\|_{L^2},
\]
and therefore
\[
|r_n|
\le
C\,\|D_n\|_{L^2}^2
=
O(z_n^2)
=
o(|z_n|).
\]
Thus, we have 
\[
z_n
=
\E\!\left[\langle g_\nu(X_n),D_n\rangle\right] + o(|z_n|).
\]
Define the conditional mean displacement
\[
\bar D_n(x):=\E[D_n\mid X_n=x].
\]
Since \(X_n\sim \nu\), we have 
\[
z_n
=
\E_\nu\!\left[\langle g_\nu(X),\bar D_n(X)\rangle\right]
+
o(|z_n|).
\]
Therefore, we have 
\[
\left(\E_\nu\!\left[\langle g_\nu(X),\bar D_n(X)\rangle\right]\right)^2
=
z_n^2+o(z_n^2).
\]
By the Cauchy--Schwarz inequality, we have 
\[
z_n^2+o(z_n^2)
\le
m\,\E_\nu[\|\bar D_n(X)\|_2^2].
\]
On the other hand, by Jensen's inequality, we have 
\[
\E_\nu[\|\bar D_n(X)\|_2^2]
\le
\E[\|D_n\|_2^2]
=
W_2(\mu_n,\nu)^2.
\]
Combining this with assumption (ii), we have 
\[
D_c(\nu,\mu_n)
\ge
\tau W_2(\mu_n,\nu)^2
\ge
\tau \E_\nu[\|\bar D_n(X)\|_2^2]
\ge
\frac{\tau z_n^2}{m}+o(z_n^2).
\]
Since
\[
D_c(\nu,\mu_n)\le R(z_n)+o(z_n^2),
\]
it follows that
\[
R(z_n)\ge \frac{\tau z_n^2}{m}+o(z_n^2).
\]
Because \(z_n\to0\) is arbitrary, this yields
\[
R(z)\ge \frac{\tau z^2}{m}+o(z^2).
\]
Together with \eqref{eq:upper_bound_Rz}, we conclude that
\[
R(z)=\frac{\tau z^2}{m}+o(z^2)
=
\frac{\tau z^2}{\E_\nu[\|g_\nu(X)\|_2^2]}+o(z^2).
\]
Finally, the perturbation \(T_z(x):=x+\Delta_z(x)\) constructed in Step 2 satisfies
\[
\Delta_z
=
\frac{z}{m}\,g_\nu+o(z)
\qquad\text{in }L^2(\nu),
\]
and
\[
D_c(\nu,T_{z\#}\nu)
=
\frac{\tau z^2}{m}+o(z^2).
\]
Since \(\kappa(T_{z\#}\nu)=z\), it in particular satisfies the weaker statement
\[
\kappa(T_{z\#}\nu)=z+o(z).
\]
Thus \(T_z\) is asymptotically optimal, and the proof is complete.}
\end{proof}

{
In the following corollary, we prove, under the sampling law \(\mathbb P_{\mu^*}\),
\[
nR_n(k^*)\Rightarrow \Upsilon,\qquad n\to\infty.
\]
Consequently, if \(\eta_{1-\delta_0}\) denotes the \((1-\delta_0)\)-quantile of \(\Upsilon\)
and is a continuity point of its distribution, then
\[
\lim_{n\to\infty}
\mathbb P_{\mu^*}
\left(
R_n(k^*)\le \frac{\eta_{1-\delta_0}}{n}
\right)
=
\mathbb P_{\mu^*}(\Upsilon\le \eta_{1-\delta_0})
=
1-\delta_0.
\]
}
{
{The next corollary relies on a local uniform version of Theorem \ref{thm:projection_law} to justify the random
plug-in. Under the compact-support assumption, this uniformity follows by the same argument
as in Theorem \ref{thm:projection_law}, with constants chosen uniformly in a sufficiently small \(W_2\)-neighborhood
of \(\mu^*\).} 
Moreover, the compactness can be relaxed to sub-Gaussian, and the computation is shown in the supplementary materials for such a generalization.}

\begin{corollary}\label{Coro2}
Suppose Assumptions \ref{ass:utility} and \ref{ass:compact_B} hold. {We assume that $\operatorname{dist}(\operatorname{supp}\mu^*,K^c)>0$, $\E_{\mu^*}\!\left[\|\nabla \kappa'_{\mu^*}(B)\|_2^2\right] > 0,$
and the scalar map $u \mapsto -I'\!\left(\frac{k^*e^{-rT}}{u}\right)\frac{k^*e^{-rT}}{u^2}$ is locally Lipschitz on $(0,\infty)$.} Define the RWP function for the Merton budget constraint as
$$
R_n(k) = \inf \left\{ D_c(\mathbb{P}_n, \mu) : \int_{\R^{d}} I\left( \frac{k e^{-rT}}{\E_{\mu}[L_T(B, y)]} \right) \varphi_T(y)  dy = x_0 e^{rT} \right\}.
$$ Then with the cost function $c(x,y) = ||x-y||_2^2$, {under the sampling law \(\mathbb P_{\mu^*}\),}
$$
n R_n(k^*) \Rightarrow \Upsilon, \quad \text{as } n \to \infty,
$$
where $\Upsilon$ is a non-negative random variable given by
\begin{equation}\label{preL}
\Upsilon = \frac{Z^2}{\E_{\mu^*}[\|\nabla \kappa'_{\mu^*}(B)\|_2^2]},    
\end{equation}
with $Z \sim \mathcal{N}(0, \sigma^2)$ and
\begin{equation}\label{sigma}
\sigma^2 = \var_{\mu^*}\left( \int_{\R^d} g'(F(y)) L_T(B, y) \varphi_T(y) dy \right),    
\end{equation}
where $F(y) = \E_{\mu^*}[L_T(B, y)]$,
\[
g'(F) = -I'\!\left( \frac{k^* e^{-rT}}{F} \right) \frac{k^* e^{-rT}}{F^2},
\]and $$\nabla \kappa'_{ \mu^*}(b) = -\int_{\R^d} I'\left( \frac{k^* e^{-rT}}{F(y)} \right) \cdot \frac{k^* e^{-rT}}{F(y)^2} \cdot \nabla_b L_T(b, y) \varphi_T(y) dy.$$

\end{corollary}

{\begin{remark}
For a power utility $U(x)=x^\alpha/\alpha$ with $\alpha<1$ and $\alpha\neq 0$, the map $u \mapsto -I'\!\left(\frac{k^*e^{-rT}}{u}\right)\frac{k^*e^{-rT}}{u^2}$
is locally Lipschitz on $(0,\infty)$. In addition, if $\mu^*$ is compactly supported with a continuous density and $\nabla \kappa'_{\mu^*}\not\equiv 0$, then $\E_{\mu^*}\!\left[\|\nabla \kappa'_{\mu^*}(B)\|_2^2\right] > 0.$
\end{remark}}

\begin{proof}
We adapt Theorem~\ref{thm:projection_law} to the DRBC Merton setting.  
Recall
\[
\kappa_k(\mu)
= \int_{\R^d} I\!\left(\frac{k e^{-rT}}{\E_{\mu}[L_T(B,y)]}\right)\varphi_T(y)\,dy
- x_0 e^{rT},
\qquad 
\mathcal{M}_k = \{\mu:\kappa_k(\mu)=0\},
\]
and
\[
R_n(k)=\inf\{D_c(\mathbb{P}_n,\mu):\mu \in\mathcal{M}_k\}.
\]
Let $\mu^*$ be the true model and $k^*$ be the associated multiplier with 
$\kappa_{k^*}(\mu^*)=0$.

\medskip
\noindent\textbf{Step 1: Wasserstein derivative of $\kappa_{k^*}$.}
As in Corollary~\ref{thm:asymptotic_perturbation_cor}, consider {$\mu^\epsilon=(1-\epsilon)\mu^*+\epsilon\delta_b$} and write 
$F(y)=\E_{\mu^*}[L_T(B,y)]$.  
Differentiating under the integral (justified by Assumptions~\ref{ass:utility}--\ref{ass:compact_B}) yields
\[
\frac{\delta\kappa_{k^*}}{\delta\mu^*}(b)
=
-\int_{\R^d} 
I'\!\left(\frac{k^* e^{-rT}}{F(y)}\right)
\frac{k^* e^{-rT}}{F(y)^2}
L_T(b,y)\,\varphi_T(y)\,dy.
\]
Therefore, the Wasserstein gradient is
\begin{equation}\label{eq:kappa_grad_final_camera}
\nabla\kappa'_{\mu^*}(b)
=
-\int_{\R^d} 
I'\!\left(\frac{k^* e^{-rT}}{F(y)}\right)
\frac{k^* e^{-rT}}{F(y)^2}\,
\nabla_b L_T(b,y)\,\varphi_T(y)\,dy.
\end{equation}
\medskip
\noindent\textbf{Step 2: Linearization and CLT for $\kappa_{k^*}(\mathbb P_n)$.}
{
Let $F_n(y):=\mathbb E_{\mathbb P_n}[L_T(B,y)]$ and $F(y):=\mathbb E_{\mu^*}[L_T(B,y)]$. Then 
define $\Phi(f)
:=
\int_{\mathbb R^d}
I\!\left(\frac{k^*e^{-rT}}{f(y)}\right)\varphi_T(y)\,dy
-
x_0e^{rT},$
so that $\Psi(\mathbb P_n)=\kappa_{k^*}(\mathbb P_n)=\Phi(F_n)$ and $\Psi(\mu^*)=\kappa_{k^*}(\mu^*)=\Phi(F).$
Let $\mathcal{H}:= L^2(\varphi_T)$ be the separable Hilbert space with the norm $h \in \mathcal{H}$, $\left \Vert h \right\Vert^2 = \int_{\mathbb{R}^d}h^2(y)\varphi_T(y)dy$. By the compact support of $\mu^*$ and the exponential bounds on $L_T(\cdot,y)$, both $F_n$ and $F$ belong to $\mathcal H$ almost surely, and $\Phi$ is Hadamard differentiable at $F$, tangentially to $\mathcal H$, with derivative
\[
\Phi_F'(v)
=
\int_{\mathbb R^d}
g'(F(y))\,v(y)\,\varphi_T(y)\,dy,
\qquad 
g'(F)
=
-\,I'\!\left(\frac{k^* e^{-rT}}{F}\right)\frac{k^* e^{-rT}}{F^2}.
\]
Moreover,
\[
\sqrt n\big(F_n-F\big)
=
\frac{1}{\sqrt n}\sum_{i=1}^n
\Big(L_T(B^{(i)},\cdot)-F(\cdot)\Big)
\Rightarrow \mathbb G
\qquad\text{in }\mathcal H,
\]
for some centered Gaussian element $\mathbb G$. Therefore, by the functional delta method
\cite{shapiro2009stochastic},
\begin{align*}
\sqrt n\big(\Psi(\mathbb P_n)-\Psi(\mu^*)\big)
&=
\Phi_F'\!\big(\sqrt n(F_n-F)\big)+o_{\mu^*}(1) \\
&=
\sqrt n\int_{\mathbb R^d}
g'(F(y))\big(F_n(y)-F(y)\big)\varphi_T(y)\,dy
+o_{\mu^*}(1) \\
&=
\sqrt n\int_{\mathbb R^d}
g'(F(y))
\left(
\int_{\mathbb R^d} L_T(b,y)\,(\mathbb P_n-\mu^*)(db)
\right)
\varphi_T(y)\,dy
+o_{\mu^*}(1) \\
&=
\sqrt n\int_{\mathbb R^d}
\left(
\int_{\mathbb R^d}
g'(F(y))\,L_T(b,y)\,\varphi_T(y)\,dy
\right)
(\mathbb P_n-\mu^*)(db)
+o_{\mu^*}(1) \\
&=
\sqrt n\int_{\mathbb R^d}
\frac{\delta \kappa_{k^*}}{\delta \mu^*}(b)\,
(\mathbb P_n-\mu^*)(db)
+o_{\mu^*}(1).
\end{align*}
Since $\mathbb P_n-\mu^*$ is a signed measure with total mass zero, we may replace
$\frac{\delta \kappa_{k^*}}{\delta \mu^*}$ by its centered version
\[
\psi(b)
:=
\frac{\delta \kappa_{k^*}}{\delta \mu^*}(b)
-
\mathbb E_{\mu^*}\!\left[
\frac{\delta \kappa_{k^*}}{\delta \mu^*}(B)
\right].
\]
Therefore, 
\[
\sqrt n\big(\kappa_{k^*}(\mathbb P_n)-\kappa_{k^*}(\mu^*)\big)
=
\frac{1}{\sqrt n}\sum_{i=1}^n \psi(B^{(i)})
+o_{\mu^*}(1).
\]
By the compact support of $\mu^*$ and the exponential bound on $L_T(\cdot,y)$, we have
$\psi(B)\in L^2(\mu^*)$. Hence, the classical CLT yields
\[
\frac{1}{\sqrt n}\sum_{i=1}^n \psi(B^{(i)}) \Rightarrow Z,
\qquad Z\sim \mathcal N(0,\sigma^2),
\]
with
\[
\sigma^2
=
\var_{\mu^*}\!\left(
\int_{\mathbb R^d}
g'(F(y))\,L_T(B,y)\,\varphi_T(y)\,dy
\right).
\]
Combining the last two displays gives
\[
\sqrt n\big(\kappa_{k^*}(\mathbb P_n)-\kappa_{k^*}(\mu^*)\big)
\Rightarrow Z.
\]
Since $\kappa_{k^*}(\mu^*)=0$, if we write $z_n:=\kappa_{k^*}(\mathbb P_n)$,
then $\sqrt n\, z_n \Rightarrow Z,$
and in particular $z_n=O_{\mu^*}(n^{-1/2})$.
For the present compact-support setting, the same conclusion may also be obtained directly from the first-order von Mises expansion of $\Psi$ at $\mu^*$ together with the classical CLT.}

\medskip
{\noindent\textbf{Step 3: Random plug-in via the local projection law.}
Let $\kappa := \kappa_{k^*}$ and $m(\mu) := \E_{\mu}\!\left[\|\nabla \kappa'_{\mu}(B)\|_2^2\right].$
By the nondegeneracy assumption, $m(\mu^*)>0$. Moreover, by the continuity established in Step 1, there exist $\rho_0>0$ and $c_0>0$ such that $m(\mu) \ge c_0$ whenever $W_2(\mu,\mu^*)<\rho_0.$}

{Because all admissible priors are supported on the same compact set, the proof of Theorem~\ref{thm:projection_law} applies with constants chosen uniformly for deterministic base laws $\mu$ in a sufficiently small $W_2$-neighborhood of $\mu^*$. By $\operatorname{dist}(\operatorname{supp}\mu^*,K^c)>0$, the bounded local perturbations remain in $K$ for all sufficiently small perturbation levels. Hence, the same local projection expansion holds when the optimization is restricted to $\mathcal{P}(K)$. In particular, for such $\mu$ and for $\delta$ near $0$, if we define
\[
R_{\mu}(\delta)
:=
\inf\left\{
D_c(\mu,\tilde{\mu})
:
\kappa(\tilde{\mu})-\kappa(\mu)=\delta
\right\},
\]
then, since here $c(x,y)=\|x-y\|_2^2$, the same proof as Theorem~\ref{thm:projection_law} yields
\[
R_{\mu}(\delta)
=
\frac{\delta^2}{m(\mu)}
+
o(\delta^2),
\qquad \delta\to 0,
\]
where the remainder is uniform for $\mu$ in that neighborhood.
Now recall that $z_n = \kappa(\mathbb{P}_n)=\kappa_{k^*}(\mathbb{P}_n)$. Then by Step 2, we have $\sqrt{n}\, z_n \Rightarrow Z,$
and hence $z_n = O_{\mu^*}(n^{-1/2})$. Since $\mu^*$ is compactly supported, $\mathbb{P}_n$ is supported on the same compact set almost surely, and $\mathbb{P}_n \to \mu^*$ in $W_2$ almost surely. Therefore, with probability tending to one, $\mathbb{P}_n$ belongs to the above neighborhood of $\mu^*$ and $|z_n|$ is sufficiently small.
Observe that
\[
R_n(k^*)
=
\inf\left\{
D_c(\mathbb{P}_n,\tilde{\mu})
:
\kappa(\tilde{\mu})=0
\right\}.
\]
Since $\kappa(\mathbb{P}_n)=z_n$, this can be rewritten as
\[
R_n(k^*)
=
\inf\left\{
D_c(\mathbb{P}_n,\tilde{\mu})
:
\kappa(\tilde{\mu})-\kappa(\mathbb{P}_n)=-z_n
\right\}
=
R_{\mathbb{P}_n}(-z_n).
\]
Applying the above deterministic local projection law with base measure $\mathbb{P}_n$ and level $-z_n$, we obtain
\[
R_n(k^*)
=
\frac{z_n^2}{m(\mathbb{P}_n)}
+
o_{\mu^*}(z_n^2).
\]
By the continuity of $m(\cdot)$ at $\mu^*$, $m(\mathbb{P}_n)\to m(\mu^*)
=
\E_{\mu^*}\!\left[\|\nabla \kappa'_{\mu^*}(B)\|_2^2\right]$ in probability.
Therefore,
\[
R_n(k^*)
=
\frac{z_n^2}{\E_{\mu^*}[\|\nabla \kappa'_{\mu^*}(B)\|_2^2]}
+
o_{\mu^*}(z_n^2).
\]
Multiplying by $n$ and using Step 2 together with Slutsky's theorem gives
\[
nR_n(k^*)
=
\frac{(\sqrt{n}\, z_n)^2}{\E_{\mu^*}[\|\nabla \kappa'_{\mu^*}(B)\|_2^2]}
+
o_{\mu^*}(1)
\Rightarrow
\frac{Z^2}{\E_{\mu^*}[\|\nabla \kappa'_{\mu^*}(B)\|_2^2]}.
\]
This proves the desired weak convergence.}

\end{proof}

Based on the above discussion, we can give the following recipe for computing the optimal DRBC policies with the Wasserstein uncertainty set with choosing the optimal $\delta$ based on the data (for a certain time window). This can be viewed as a nonlinear extension of the RWPI method (e.g., see \cite{BlanchetChenZhou2021, Blanchet2021WassersteinDRO}).
\begin{itemize}
\label{bullpoint}
    \item (1) Collect return data $\{B^{(i)}\}_{i = 1, \ldots, n}$.
    \item (2) Use the collected data $\{B^{(i)}\}_{i = 1, \ldots, n}$ to solve the equation $\int_{\mathbb{R}^d }I\left(\frac{ke^{-rT}}{\mathbb{E}_{\mathbb{P}_n}\left[L_T(B,y)\right]}\right)\varphi_T(y)dy = x_0e^{rT}$ (this corresponds to estimating $\mathcal{K}(x_0)$). Denote the solution by $\hat{k}$.
    \item (3) Obtain independent samples $Y_1, \ldots, Y_N$ from $\mathcal{N}(0,TI_d)$. Compute $\hat{F}(Y_i)  = \E_{\mathbb{P}_n}[L_T(B, Y_i)]$ for each sample using the collected data $\{B^{(i)}\}_{i = 1, \ldots, n}$, and then compute {$\E_{\mathbb{P}_n}[\|\nabla \kappa'_{\mathbb{P}_n}(B)\|_2^2]$} as an estimate of the denominator of $\Upsilon,$ where {$\nabla \kappa'_{\mathbb{P}_n}(b)$} is computed by the Monte Carlo method with $\hat{F}(Y_i)$, $L_T(b, Y_i)$, and $\hat{k}$ derived in Step (2).

    \item (4) Estimate (\ref{sigma}) by the same Monte Carlo method as in Step (3) and denote the estimated variance by $\hat{\sigma}$. Define 
    {$\hat{\Upsilon} = \frac{\hat{Z}^2}{\E_{\mathbb{P}_n}[\|\nabla \kappa'_{\mathbb{P}_n}(B)\|_2^2]}$}, where $\hat{Z} \sim \mathcal{N}(0, \hat{\sigma}^2)$.

    \item (5) Obtain independent samples $\Upsilon_1, \dots, \Upsilon_K$ from $\hat{\Upsilon}$. Let $\hat{\eta}_{.95}$ be the $95\%$ quantile of the sample collections $\Upsilon_1, \dots, \Upsilon_K$ ($K$ is a sequence such that $K \to \infty$ as $n \to \infty$, for example, $K = \log n$; see Algorithm 1 in \cite{Blanchet2021WassersteinDRO}). 
    \item (6) Set $\delta = \frac{\hat{\eta}_{.95}}{n}$ and approximate the solution to problem \eqref{FWProblem} using Corollary \ref{thm:asymptotic_perturbation_cor} and obtain a worst-case drift prior {\(\widehat\mu^{\text{wc}}\)} for the drift $B$.
    \item (7) Plug in this {\(\widehat\mu^{\text{wc}}\)} into the closed-form formula of the optimal fraction with the observed stock prices and interest rates in Theorem \ref{thm:karatzas_zhao_general_compact} replacing $\mu$ by {\(\widehat\mu^{\text{wc}}\)}.
\end{itemize}

\begin{remark}
{In Steps (3)-(4), the unknown population quantities in Corollary \ref{Coro2} are replaced by plug-in estimators based on $\mathbb{P}_n$ and $\hat{k}$. More precisely, define $\widehat{F}_n(y):=\E_{\mathbb{P}_n}[L_T(B,y)]$ and $\widehat m_n:=\E_{\mathbb{P}_n}\!\left[\left\|\nabla \kappa'_{\mathbb{P}_n,\hat{k}}(B)\right\|_2^2\right],$
where $\nabla \kappa'_{\mathbb{P}_n,\hat{k}}(b)
=
-\int_{\mathbb{R}^d}
I'\!\left(\frac{\hat{k}e^{-rT}}{\widehat{F}_n(y)}\right)
\frac{\hat{k}e^{-rT}}{\widehat{F}_n(y)^2}
\,\nabla_b L_T(b,y)\,\varphi_T(y)\,dy,$
and the integral is evaluated numerically by Monte Carlo using the Gaussian samples
$Y_1,\ldots,Y_N\sim\mathcal{N}(0,TI_d)$. Likewise, $\hat{\sigma}^2$ is computed from the corresponding
plug-in version of \eqref{sigma}, and Step (4) uses $\hat{\Upsilon}=\frac{\hat{Z}^2}{\widehat m_n}$ and $\hat{Z}\sim\mathcal{N}(0,\hat{\sigma}^2).$}

{This plug-in implementation does not affect the first-order asymptotic calibration. Indeed, under
Assumptions \ref{ass:utility} and \ref{ass:compact_B}, we have
$\mathbb{P}_n\to\mu^*$ in $W_2$ almost surely and $\hat{k}\to k^*$. Moreover, by the common
compact support and the continuity properties used in Corollaries \ref{thm:asymptotic_perturbation_cor}
and \ref{Coro2}, the map $(\mu,k,b)\longmapsto \nabla \kappa'_{\mu,k}(b)$
is jointly continuous with a common integrable envelope in a neighborhood of
$(\mu^*,k^*)$. Consequently,
\[
\widehat m_n
\to
\E_{\mu^*}\!\left[\left\|\nabla \kappa'_{\mu^*,k^*}(B)\right\|_2^2\right]
\qquad\text{and}\qquad
\hat{\sigma}^2\to \sigma^2
\]
in probability. Therefore, replacing $(\mu^*,k^*)$ by $(\mathbb{P}_n,\hat{k})$ contributes only a
lower-order plug-in error and does not change the first-order radius calibration. If, in addition,
the inner Monte Carlo sizes $N$ and $K$ are taken to infinity with $n$, then the additional simulation
error is asymptotically negligible.}
\end{remark}

\section{Synthetic Experiment}\label{syn}

\subsection{Understanding Model Parameters}
\label{param}

In this section, we generate synthetic data in a high-dimensional setting and understand how parameters affect the performance of different models. We let the ground-truth drift be 
\begin{equation}
    B_{it} = \frac{B_0}{2}\left(1 + 2\cos(2\pi\kappa_i t)\right) \label{eq:highdimdrift}
\end{equation} where $i$ indexes the stock. Each stock's $\kappa_i$ is sampled from the same Gaussian distribution. We set the total number of stocks to be 20. The synthetic stock data is generated using (\ref{stock}) with drift (\ref{eq:highdimdrift}), and we set the volatility matrix to be diagonal for simplicity. Throughout the paper, we always assume time 1 to be 1 year, and the unit time is $dt$. Since our model is continuous-time, $dt$ should be very granular. The trading rule is simple: for every 22 $dt$, calculate data-driven parameters and portfolio weights based on data of the previous 2520 $dt$ on a rolling basis and trade. Bayesian Merton and DRBC update portfolio weights at every \(dt\), whereas the other methods use static weights over each rebalancing period. We evaluate the performance using the Sharpe ratio with interest rate $r=1\%$ on the last 252 $dt$'s wealth. We vary $B_0,dt$ and the distribution to sample $\kappa$ to understand the models better, and we would like to use the intuition here to guide our experiments on real data. 

All experiment procedures follow the algorithm in Step \eqref{bullpoint} above. To be more precise, we estimate {the nominal prior} from the data using the “Consecutive Drift” approach. We then use the projection method to estimate $\delta$. After that, we compute $\Delta(B)$, which is not a projection, but rather a perturbation according to Theorem \ref{thm:asymptotic_perturbation}. Then we get the modified prior $B^*=B+\Delta(B)$ and then apply formulae in Theorem \ref{thm:karatzas_zhao_general_compact} with the modified prior distributions to get portfolio weights. Experiment results are shown in Table~\ref{sharpe_highdim}, which is the average of 100 simulations (We use Sharpe Ratio here as a performance metric. Tables with terminal utility are reported in Appendix \ref{utility_highdim}).  In every simulation, we use a different seed to sample both $\kappa_i$ and all the random numbers used in Theorem \ref{thm:karatzas_zhao_general_compact}. The value of \(1/(dt\times252)\) indicates how many subperiods each trading day is divided into. Notice that the ground-truth drift is a periodic function with periodicity $\frac{1}{\kappa}$. The drift distribution is estimated as follows: update of drift estimation happens every 30 $dt$, each time we use previous 2520 $dt$ to estimate the drift distribution with batched, disjoint time windows. Specifically, we split the 2520 $dt$ into 10 non-overlapping periods, and estimate annualized return with each stock's period cumulative return as in Section \ref{projection}. We call this “Consecutive Drift”. We call the market condition smooth when $\kappa \sim \mathcal{N}(0,1)$ and volatile when $\kappa \sim \mathcal{N}(12,10)$. During smooth market conditions, models with frequent rebalancing, like Bayesian Merton, DRBC, and DRC, benefit from finer time trading resolution, while DRMV shows worse results. During volatile times, such benefits are not significant. In other words, if one believes the economy will grow steadily with no crisis for a long time, yet one also wants to avoid huge downside risk caused by noise, then DRBC could be a reasonable choice.

\begin{table}[H]
\centering
\small
\renewcommand{\arraystretch}{1.25}
\begin{tabular}{cccccccc}
\toprule
\multicolumn{3}{c}{Parameters} & \multicolumn{5}{c}{Average Sharpe Ratio} \\
\midrule
$B_0$ & $1/(dt\times252)$ & $\kappa$ & Bayesian Merton & DRBC & DRMV\_no\_rf & DRMV\_rf & DRC \\
\midrule

0.2 & 6  & $\mathcal{N}(0,1)$ 
& 0.855 & 0.868 & 1.040 & 1.017 & 0.528 \\
&     &        
& (2.340) & (2.339) & (2.719) & (2.721) & (2.609) \\

0.2 & 11 & $\mathcal{N}(0,1)$ 
& 0.880 & 0.893 & 0.887 & 0.861 & 0.641 \\
&     &        
& (3.128) & (3.126) & (3.474) & (3.585) & (3.479) \\

0.4 & 6  & $\mathcal{N}(0,1)$ 
& 2.050 & 2.058 & 2.282 & 2.363 & 1.657 \\
&     &        
& (2.481) & (2.483) & (2.788) & (2.837) & (2.567) \\

\textbf{0.4} & \textbf{11} & $\bm{\mathcal{N}(0,1)}$ 
& \textbf{2.132} & \textbf{2.137} & \textbf{1.989} & \textbf{2.104} & \textbf{1.978} \\
&     &        
& \textbf{(3.120)} & \textbf{(3.126)} & \textbf{(3.547)} & \textbf{(3.607)} & \textbf{(3.402)} \\

0.2 & 6  & $\mathcal{N}(12,10)$ 
& 0.720 & 0.728 & 0.922 & 0.902 & 0.422 \\
&     &        
& (2.312) & (2.310) & (2.645) & (2.683) & (2.608) \\

0.2 & 11 & $\mathcal{N}(12,10)$ 
& 0.720 & 0.739 & 0.870 & 0.840 & 0.366 \\
&     &        
& (3.267) & (3.267) & (3.500) & (3.637) & (3.576) \\

0.4 & 6  & $\mathcal{N}(12,10)$ 
& 1.796 & 1.801 & 2.034 & 2.106 & 1.399 \\
&     &        
& (2.398) & (2.398) & (2.654) & (2.702) & (2.622) \\

0.4 & 11 & $\mathcal{N}(12,10)$ 
& 1.809 & 1.829 & 1.953 & 2.012 & 1.334 \\
&     &        
& (3.403) & (3.412) & (3.631) & (3.640) & (3.640) \\

\bottomrule
\end{tabular}
\caption{Sharpe Ratio comparison across parameter settings over 100 simulations. Means reported with standard deviations in parentheses.}
\label{sharpe_highdim}
\end{table}

We also try different ways to estimate the drift distribution and how projection works in DRBC and benchmark methods. Details of DRMV and DRC are discussed in Section \ref{drmv_detail} and \ref{drc_detail}. For drift estimation, we try day-of-week aggregation within larger windows (e.g., averages of “Mondays”, “Tuesdays”, etc.); we call it “Type Drift”. In the experiments, all $dt$ are divided into 10 types to make a more fair comparison with the Consecutive Drift approach.  For the perturbation $\Delta$, we have two choices: change $\Delta$ every $dt$ by changing the plan time, or keep $\Delta$ static and make the plan time the same as drift update frequency. The average Sharpe Ratio results are shown in Table~\ref{tab:sharpe} (terminal utility results can be found in appendix \ref{utility_highdim}). We use $B_0=0.4, dt=\frac{1}{252 \times 11},\kappa\sim \mathcal{N}(0,1)$ suggested by Section \ref{param}. The consecutive drift and the static projection achieve the best performance. A reasonable explanation for the “Consecutive” approach better than the “Type” approach is that the implicit number of “types” appears difficult to estimate (e.g., if we think that there is weak seasonality, the number of types should be around 5-7, but likely there are other time-patterns that do not align for all stocks). Only if the type number perfectly echoes with $B_{it}$'s period can we get better drift estimation, yet it is a rare case in practice. The “Consecutive” approach smooths the drift within the time batch and has a better empirical practice overall. Regarding the use of the stating projection, since here we use a very small $dt$, the time-varying projection is likely influenced by extreme values at the end of the drift update period. Simulation results follow our intuition, though the difference is small. 

\begin{table}[htbp]
\centering
\small
\caption{Sharpe Ratios on different drift estimation and projection methods over 100 simulations (means on first line; standard deviations in parentheses on the next line)}
\label{tab:sharpe}
\begin{tabular}{llccccc}
\toprule
Drift & Projection  & Bayesian Merton & DRBC & DRMV\_no\_rf & DRMV\_rf & DRC \\
\midrule
\multirow{2}{*}{Consecutive} & \multirow{2}{*}{Static} 
  & 2.1319 & 2.1374 & 1.9890 & 2.1037 & 1.9778 \\
 &
  & (3.1203) & (3.1261) & (3.5466) & (3.6072) & (3.4019) \\
\multirow{2}{*}{Consecutive} & \multirow{2}{*}{Time-varying} 
  & 2.1319 & 2.1331 & 1.9890 & 2.1037 & 1.9778 \\
 &
  & (3.1203) & (3.1268) & (3.5466) & (3.6072) & (3.4019) \\

\multirow{2}{*}{Type} & \multirow{2}{*}{Static}
  & 1.9815 & 1.9835 & 1.9890 & 2.1037 & 1.9995 \\
&
  & (3.4129) & (3.4132) & (3.5466) & (3.6072) & (3.4382) \\

\multirow{2}{*}{Type} & \multirow{2}{*}{Time-varying} 
  & 1.9815 & 1.9814 & 1.9890 & 2.1037 & 1.9995 \\
 & 
  & (3.4129) & (3.4142) & (3.5466) & (3.6072) & (3.4382) \\

\bottomrule
\end{tabular}
\end{table}

\subsection{Role of Radius}
In this synthetic experiment, we show how different radii change the performance for DRC and DRBC compared to optimal strategies, which implicitly proves the importance of data-driven radius determination in Section \ref{projection}. The experiment is a multi-dimensional setting with the Wasserstein ball as the uncertainty set, and we use the projection approach in Corollary \ref{Coro2}.

\begin{figure}[htbp]
    \centering
    \includegraphics[width=0.8\textwidth]{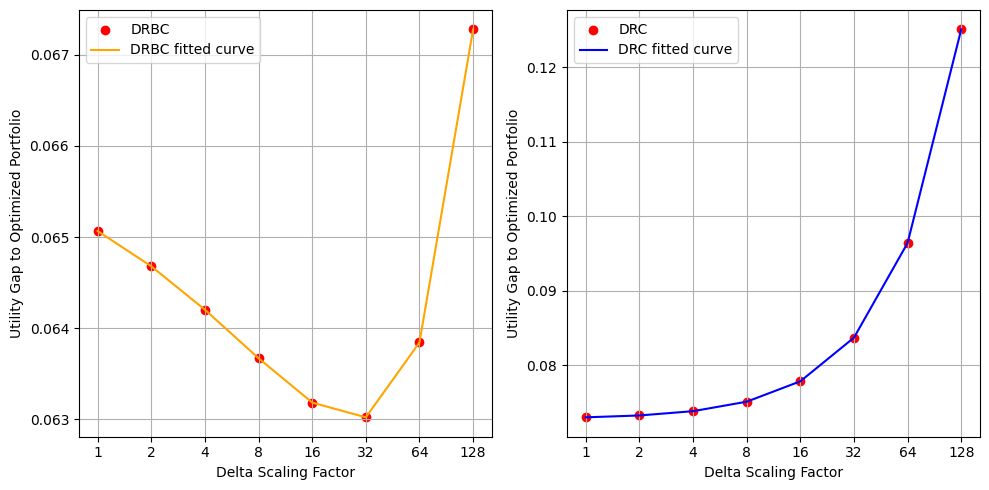}
    \caption{Expectation of terminal utility gap versus $\delta$ scaling factor for DRC and DRBC in a high-dimensional case with Wasserstein uncertainty.}
    \label{fig:gaphighd}
\end{figure}

The data generation part is similar to Section \ref{param}. The drift terms are sampled from (\ref{eq:highdimdrift}), and the synthetic stock data is sampled using (\ref{stock}) with 3024 $dt$. We use Merton's formula 
\begin{equation}
\label{mertonformula}
    \frac{\pi(t)}{X(t)}=\frac{1}{1-\alpha}(\sigma\sigma^\top)^{-1}(B_t-r\mathbf{1})
\end{equation} to calculate the optimal high-dimensional portfolio strategy (policy) and calculate the oracle average terminal utility over one hundred simulated paths. The drift distribution is estimated with batched, disjoint time windows. In this way, we have 10 support vectors of the drift distribution, and we assume the drift is uniform on these 10 supports. Here we use Wasserstein uncertainty and the algorithm in Step (\ref{bullpoint}) to get the portfolio weights. We do not manually set $\delta$ levels, but use a scale factor to scale up the data-driven base $\delta$ to explore the best radius since the data-driven $\delta$ is rough. Since every point in Figure \ref{fig:gaphighd} is an average of one hundred simulations, and our data-driven $\delta$ calculation gives different base $\delta$ every simulation, we do not show exact $\delta$ in use but show the scaling factor instead, with the average base $\delta$ at the level of $10^{-3}$. According to Figure \ref{fig:gaphighd}, we observe that for DRBC, the radius $\delta$ needs to be carefully calibrated to achieve the best terminal utility. For DRC, a larger $\delta$ always leads to a worse performance.

{This experiment also helps address the policy-level interpretation of the RWPI-based radius choice in Section \ref{projection}. The theory there calibrates $\delta_n$ through coverage of the oracle parameter $k^*$, rather than through a direct utility bound for the resulting robust policy. Figure \ref{fig:gaphighd} provides complementary empirical evidence at the decision level: the left panel reports the expected terminal-utility gap between the DRBC policy and the oracle policy as the radius is varied around the data-driven baseline. We observe that this gap remains relatively stable for moderate radius choices and becomes worse when the radius is taken too large, suggesting that the RWPI-calibrated ambiguity size places the robust policy in a regime that is reasonably close to oracle performance in our experiments. This experiment also shows that around the RWPI-calibrated baseline, the DRBC policy remains relatively stable for moderate radius choices and deteriorates only when the ambiguity radius is taken too large.}

\section{Real-data Experiment}\label{real}
We use the same dataset as  \cite{BlanchetChenZhou2021} to compare overall performance.

\subsection{ Experiment Design and Data Preparation}
Inspired by the synthetic data experiment, which shows that more granular time period helps continuous-time models, we choose to set the trading frequency daily, which is more granular than monthly in \cite{BlanchetChenZhou2021}. We get real data from Wharton Research Database Service (WRDS). The dataset contains all Standard \& Pool's 500 constituents data from 2017-01-01 to 2024-12-31. We choose this time period for two reasons: timeliness and variety of market events. {During this period, the market experienced a stable uptrend, COVID-19, inflation concerns, market recovery, and the boom of artificial intelligence, making this period highly unpredictable and well suited for testing the ability of different strategies to adapt to changing market conditions. }

The experiment is done with a rolling time window of one month. We use 5 years of the previous data as the training set for DRBC to get the uncertainty radius $\hat{\delta}$, the empirical distribution $\hat{B}$, and values needed for benchmark strategies like DRC and Bayesian Merton. Here for simplicity, we let the $\delta$ scaling factor to be 1. For stock selection, we randomly sample 20 stocks from S\&P 500 constituents in the past 5 years. The real trading period starts from 2022-01-01. DRBC and Bayesian Merton strategies trade daily and follow trading rules in Theorem \ref{thm:karatzas_zhao_general_compact}, with $T$ to be two months. Other strategies get static portfolio weights for one month. To avoid future information, specifically in $Y_t$, we do not trade on the first day of every month and then start using $Y_t$ of the previous day. {Since the interest rate changed much during our trading period, increasing rapidly from about 0 percent to 5 percent and staying at that level,} we simplify by assuming all trades in the trading period with interest rate 5\%. For evaluation, we report the annualized return, the standard deviation and the Sharpe ratio for the whole time series, and assume the interest rate to be 4\%.

\subsubsection{DRMV} \label{drmv_detail}
We include the Wasserstein-robust mean-variance model of \cite{BlanchetChenZhou2021}.
Given empirical distribution 
\(\mathbb{P}_n = \frac1n \sum_{i=1}^n \delta_{R_i}\),
they consider all return distributions inside the Wasserstein ball
\[
\mathcal{U}_\delta(\mathbb{P}_n)
= \{\mu : D_c(\mathbb{P}_n, {\mu})\le \delta \}, 
\qquad 
c(u,v)=\|u-v\|_q^2.
\]
The robust Markowitz problem is
\begin{equation}
\min_{\phi}\;\; \max_{{\mu}\in \mathcal{U}_\delta(\mathbb{P}_n)} 
\phi^\top \operatorname{Var}_{{\mu}}(R)\phi
\quad 
\text{s.t. }\;
\mathbf{1}^\top\phi=1,\;
\min_{{\mu} \in\mathcal{U}_\delta(\mathbb{P}_n)} \mathbb{E}_{\mu}[R]^\top\phi \ge \bar\alpha.
\label{eq:drmv-primal-short}
\end{equation}
A key result is that \eqref{eq:drmv-primal-short} is \emph{exactly equivalent} to a regularized empirical problem:
\begin{equation}
\begin{aligned}
\min_{\phi}\quad &
\sqrt{\phi^\top \hat\Sigma\,\phi} + \sqrt{\delta}\,\|\phi\|_p,\\
\text{s.t.}\quad &
\mathbf{1}^\top\phi = 1,\qquad
\hat\mu^\top\phi \ge \bar\alpha + \sqrt{\delta}\,\|\phi\|_p,
\end{aligned}
\label{eq:drmv-regularized-short}
\end{equation}
where \(\hat\mu=\mathbb{E}_{\mathbb{P}_n}[R]\), \(\hat\Sigma=\operatorname{Var}_{\mathbb{P}_n}(R)\), and \(1/p+1/q=1\). 
Thus, Wasserstein robustness leads to a theoretically justified norm penalty \(\|\phi\|_p\), and the ambiguity radius \(\delta\) and target \(\bar\alpha\) are chosen via data-driven Wasserstein profile inference.

We use two sets of DRMV algorithm, the original one in \cite{BlanchetChenZhou2021} and the one with a risk-free asset. For the second one, we view the interest rate as the last entry of the return vector, choose the $\delta$ in the same way as \cite{BlanchetChenZhou2021}, and change the annual target return $\rho$ from $10 \%$ to $10.5 \%$. To avoid trivial results, we manually add a small noise to the interest rate.

\subsubsection{DRC}\label{drc_detail}

Classical DRC formulations \cite{HansenSargent2001,hansen2008robustness} introduce an adversary who, at every time~$t$, perturbs the model by selecting a worst-case probability measure within a $\phi$-divergence ball.  
Given a baseline model ${\mathbb{P}_0}$, at each time step the adversary selects ${\mu \ll \mathbb{P}_0}$ satisfying
\[
D_\phi(\mu\|\mathbb{P}_0)\;\le\;\delta ,
\]
where $D_\phi$ is the $\phi$--divergence generated by a convex function $\phi$ with $\phi(1)=0$.  
The controller then solves the dynamic game
\begin{equation}
\sup_{\pi\in\mathcal A(x_0)}
\inf_{{Q}^\mu\in\mathcal U_\delta^{\text{DRC}}}
\;\mathbb{E}_{{ Q^\mu}}\!\left[U(X_T)\right],
\label{eq:drc-primal}
\end{equation}
where $\mathcal U_\delta^{\text{DRC}}$ denotes the time replenished uncertainty set: the adversary is allowed to choose a new worst-case distribution at every instant.

Here, DRC is a high-dimensional implementation of \cite{Blanchet2025Duality}, where an optimization problem induced by the Hamilton-Jacobi-Bellman-Isaacs (HJBI) equation is solved to get the drift estimation, and then the Merton portfolio weight formula (\ref{mertonformula}) is used to get the portfolio weights. For the same practical concern in shorting stocks as in DRBC, we also assume that the short position of each stock cannot exceed half of the wealth. For the radius in DRC, we directly use the radius from DRBC for simplicity. Note that this is not the optimal radius for DRC. 

\subsubsection{DRBC}
To calculate empirical $\hat{B}$ and $\hat{\delta}$, we use the previous 5 years of data. We use the Ledoit-Wolf \cite{ledoit2004well} algorithm to estimate the covariance matrix, and we calibrate the items in the inverse of the volatility matrix to be within a range for numerical stability. As directly scaling up the daily return to the annual return might lead to unreasonable values, we also calibrate on the empirical center to keep outliers in a limit. Since $B$ can be roughly regarded as log return plus volatility, the limit is chosen based on the knowledge of financial markets. Besides, given the difficulty to naked short-selling in practice, we assume that the short position of each stock cannot exceed half of wealth.

\subsection{Comparison and Discussions}

\begin{figure}[htbp]
    \centering
    \includegraphics[width=0.8\textwidth]{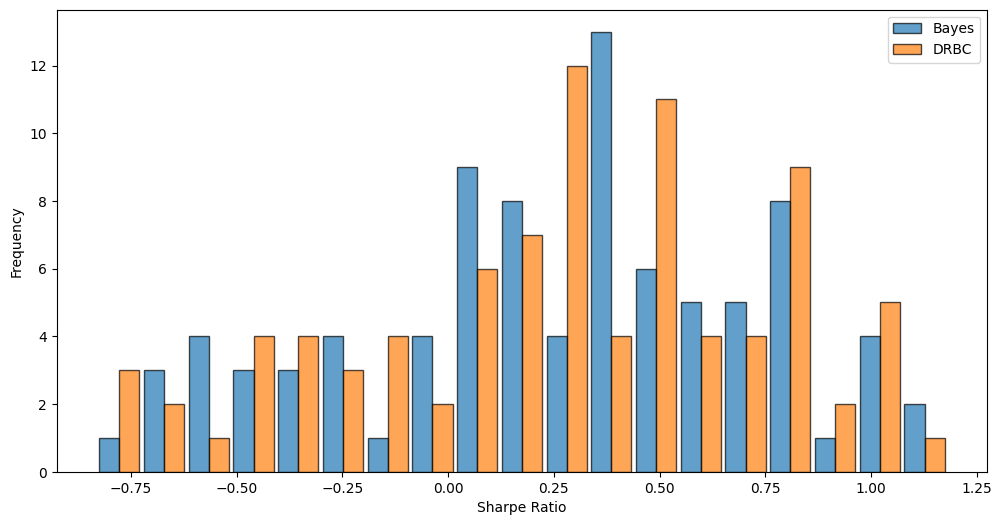}%
    \caption{Histogram of Sharpe Ratios for Bayesian Merton and DRBC}
    \label{fig:kara_drbc}
\end{figure}

We compare the histogram of the Sharpe Ratio for the aforementioned methods. Figure \ref{fig:kara_drbc} shows the distribution of the Sharpe Ratio for Bayesian Merton and DRBC. It is clear that the DRBC distribution shifts rightward and is more concentrated on a positive Sharpe Ratio. %
This suggests that priors in financial markets are highly unstable and that distributionally robust methods can be useful.

\begin{figure}[htbp]
    \centering
    \includegraphics[width=0.8\textwidth]{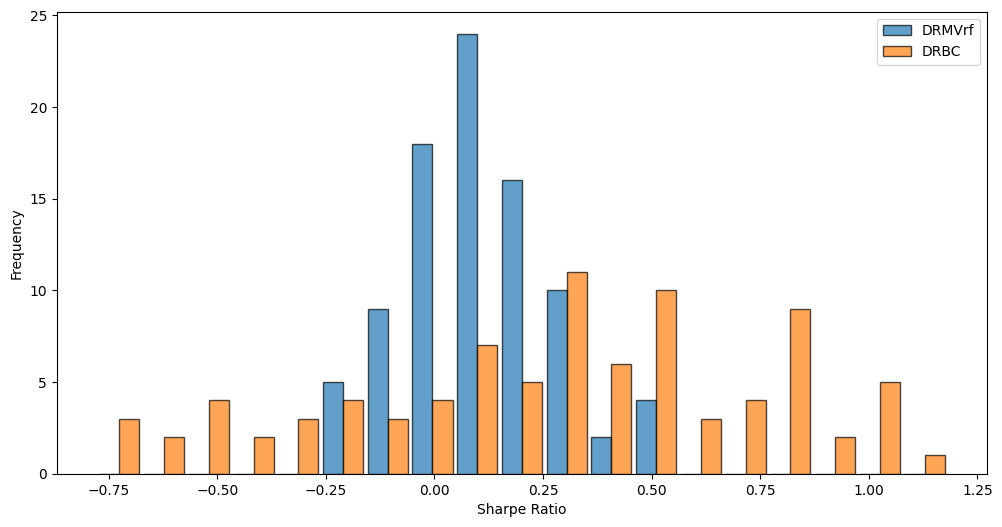}%
    \caption{Histogram of Sharpe Ratios for DRMV with a riskless asset and DRBC}
    \label{fig:hist_drmvrf_drbc}
\end{figure}

\begin{figure}[htbp]
    \centering
    \includegraphics[width=0.8\textwidth]{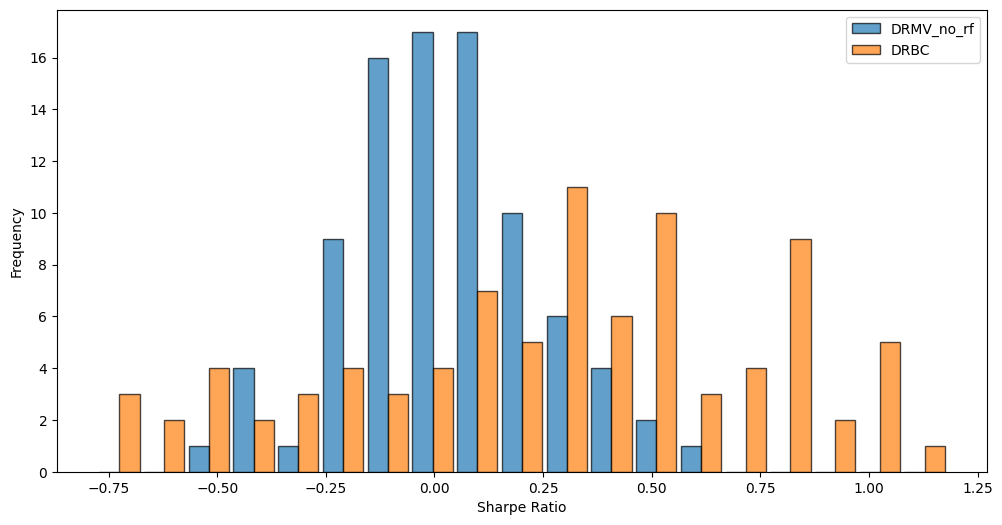}%
    \caption{Histogram of Sharpe Ratios for DRMV without a riskless asset and DRBC}
    \label{fig:hist_drmv_drbc}
\end{figure}

A comparison between DRBC and two types of DRMV appears in Figures \ref{fig:hist_drmvrf_drbc} and  \ref{fig:hist_drmv_drbc}. Both DRMV methods have more concentrated distributions than DRBC, which is reasonable since DRMV implementations do not allow short selling and borrowing money, and thus are more stable than DRBC. Despite not being concentrated, DRBC achieves a larger average Sharpe Ratio than DRMV during our volatile test period. It is also worth noting that DRBC is a continuous-time model while DRMV is not. According to \cite{BlanchetChenZhou2021}, continuous-time models tend to perform worse than discrete-time models, which further supports the effectiveness of DRBC.

\begin{figure}[htbp]
    \centering
    \includegraphics[width=0.8\textwidth]{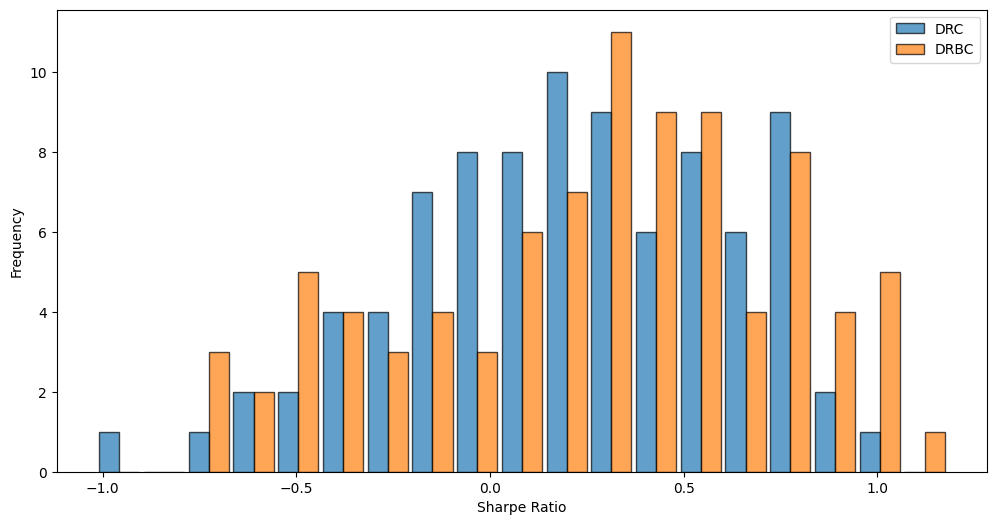}%
    \caption{Histogram of Sharpe Ratios for DRC and DRBC}
    \label{fig:hist_drc_drbc}
\end{figure}

Likewise, as in Figure \ref{fig:hist_drc_drbc}, DRBC outperforms DRC with a more right-shifted distribution, with both larger maximum and minimum Sharpe Ratios. Figure \ref{fig:hist_drc_drbc} also supports our claim that DRBC can reduce over-pessimism in DRC.

\section{Conclusions and Future Work}
In this paper, we revisited Merton's continuous-time portfolio selection problem through a distributionally robust Bayesian control lens. 
By placing a single ambiguity set on the drift prior, rather than adopting time-rectangular uncertainty on the data-generating process, we preserve the Bayesian learning structure while mitigating the over-pessimism often induced by dynamic robust control.
A Sion-type minimax swap reduces the DRBC game to a non-linear distributional optimization problem over the drift prior, allowing us to retain the Karatzas--Zhao \cite{KZ98}'s closed-form characterization of the value function and optimal policy for each candidate prior. {If the conditions required for the Sion minimax swap fail and a duality gap remains, weak
duality still gives the swapped problem as an upper bound on the robust value.} %
Therefore, the approach that we propose in this paper, we believe, could be broadly applicable to Bayesian control formulations in which the prior is imposed on certain model parameters. 

There are natural directions for future work. One of them is precisely to investigate our proposed approach, as we hinted in the previous paragraph, in distributionally robust Bayesian control settings. This paper provides a blueprint in the setting of Merton's model because it is both elegant and of significant interest (both academically and practically). A broader extension would need to be more algorithmic, dealing with efficient methods for solving the Bayesian control problem jointly with the evaluation of the required sensitivities. {Beyond the closed-form setting, we see two natural extension routes. One is a special-structure route, in which other Bayesian control models may still admit semi-closed-form or reduced-form evaluation under a fixed prior, allowing the same outer prior-robustification. The other, and more general, route is algorithmic: replace the closed-form inner solver by a numerical solver for the fixed-prior Bayesian control problem, and then optimize the outer prior-robust layer numerically.} We believe that such extensions are worth pursuing.

\section{Appendix: Additional Experiment Results}
In this section, we show additional experiment results in Sections \ref{syn} and \ref{real}.

\begin{table}[htbp]
\label{utility_highdim}
\centering
\small
\renewcommand{\arraystretch}{1.25}
\begin{tabular}{cccccccc}
\toprule
\multicolumn{3}{c}{Parameters} & \multicolumn{5}{c}{Average Terminal Utility} \\
\midrule
$B_0$ & $1/(dt\times252)$ & $\kappa$ & Bayesian Merton & DRBC & DRMV\_no\_rf & DRMV\_rf & DRC \\
\midrule

0.2 & 6  & $\mathcal{N}(0,1)$ 
& -0.389 & -0.383 & -0.319 & -0.320 & -0.787 \\
&     &        
& (0.386) & (0.373) & (0.033) & (0.030) & (1.396) \\

0.2 & 11 & $\mathcal{N}(0,1)$ 
& -0.356 & -0.353 & -0.327 & -0.327 & -0.400 \\
&     &        
& (0.225) & (0.218) & (0.024) & (0.022) & (0.298) \\

0.4 & 6  & $\mathcal{N}(0,1)$ 
& -0.277 & -0.275 & -0.303 & -0.305 & -0.769 \\
&     &        
& (0.345) & (0.338) & (0.033) & (0.030) & (1.975) \\

\textbf{0.4} & \textbf{11} & $\bm{\mathcal{N}(0,1)}$
& \textbf{-0.289} & \textbf{-0.288} & \textbf{-0.319} & \textbf{-0.319} & \textbf{-0.322} \\
&     &        
& \textbf{(0.216)} & \textbf{(0.209)} & \textbf{(0.023)} & \textbf{(0.022)} & \textbf{(0.308)} \\

0.2 & 6  & $\mathcal{N}(12,10)$ 
& -0.403 & -0.397 & -0.321 & -0.322 & -0.773 \\
&     &        
& (0.390) & (0.376) & (0.033) & (0.030) & (1.267) \\

0.2 & 11 & $\mathcal{N}(12,10)$ 
& -0.365 & -0.361 & -0.327 & -0.327 & -0.417 \\
&     &        
& (0.229) & (0.221) & (0.024) & (0.022) & (0.306) \\

0.4 & 6  & $\mathcal{N}(12,10)$ 
& -0.304 & -0.300 & -0.306 & -0.308 & -0.708 \\
&     &        
& (0.358) & (0.347) & (0.031) & (0.029) & (1.660) \\

0.4 & 11 & $\mathcal{N}(12,10)$ 
& -0.313 & -0.311 & -0.319 & -0.319 & -0.372 \\
&     &        
& (0.236) & (0.230) & (0.024) & (0.022) & (0.338) \\

\bottomrule
\end{tabular}
\caption{Terminal utility comparison across parameter settings over 100 simulations. Means reported with standard deviations in parentheses.}
\end{table}

\begin{table}[htbp]
\centering
\small
\caption{Terminal utility on different drift estimation and projection methods over 100 simulations (means on first line; standard deviations in parentheses on the next line)}
\label{tab:utility}
\begin{tabular}{llccccc}
\toprule
Drift & Projection  & Bayesian Merton & DRBC & DRMV\_no\_rf & DRMV\_rf & DRC \\
\midrule
\multirow{2}{*}{Consecutive} & \multirow{2}{*}{Static} 
  & -0.2893 & -0.2883 & -0.3186 & -0.3189 & -0.3224 \\
 & 
  & (0.2156) & (0.2092) & (0.0234) & (0.0216) & (0.3077) \\
\multirow{2}{*}{Consecutive} & \multirow{2}{*}{Time-varying} 
  & -0.2893 & -0.2889 & -0.3186 & -0.3189 & -0.3224 \\
 & 
  & (0.2156) & (0.2113) & (0.0234) & (0.0216) & (0.3077) \\

\multirow{2}{*}{Type} & \multirow{2}{*}{Static} 
  & -0.2986 & -0.2982 & -0.3186 & -0.3189 & -0.3287 \\
 & 
  & (0.2204) & (0.2192) & (0.0234) & (0.0216) & (0.3351) \\

\multirow{2}{*}{Type} & \multirow{2}{*}{Time-varying} 
  & -0.2986 & -0.2984 & -0.3186 & -0.3189 & -0.3287 \\
 & 
  & (0.2204) & (0.2196) & (0.0234) & (0.0216) & (0.3351) \\

\bottomrule
\end{tabular}
\end{table}

\subsection{Discussions for Real Data Experiment}

In Section \ref{real}, we mention the short constraint we add to make sure every single stock's short position cannot exceed half of the wealth. This is reasonable in practice, yet we want to investigate how many times this constraint is triggered in DRC, DRBC, and Bayesian methods, as well as the leverage condition among these methods. 

Across all simulations with different data-generating parameters, the short constraint happens in more than 90\% of trades for the DRC method, which is a common disadvantage for Merton-like problems since it always gives extreme values due to the estimation error of average return and variance matrix and the unconstrained optimization. Short constraint appears in about 70\% of trades in the Bayesian setting and about 40\% in the DRBC setting, illustrating that Bayesian methods by \cite{KZ98} can mitigate extreme weights in Merton-like problems and DRBC can further decrease extreme weights. 

We also check the overall leverage ratio for three algorithms. DRC on average uses about 13x leverage, with outliers like 31x. Bayesian uses about 7x and DRBC 5x on average. The leverage ratio highly depends on how drift B is estimated, which depends on the training period. If finer constraints on both short and long positions are added, it will give results with less leverage, yet it cannot fully reveal the strength of the DRBC method. 

{
\section*{Acknowledgments}
The authors would like to thank the editors and two anonymous referees for their valuable comments and insightful suggestions that help us greatly improve the paper. The authors acknowledge the use of GPT 5.5 for editorial revisions and writing. We accept responsibility for final product. J. Blanchet gratefully acknowledges support from DoD through ONR N00014-24-1-2655, and support from NSF via grants 2312204 and 2403007. Y. Liu acknowledges financial support from the National Natural Science Foundation of China (Grant No. 12401624), Guangdong Science and Technology Program (Grant No. 2024QN11X076), Shenzhen Science and Technology Program (Grant Nos. RCBS20231211090814028, JCYJ20250604141203005, 2025TC0010) and The Chinese University of Hong Kong (Shenzhen) University Development Fund (Grant No. UDF01003336) and is partly supported by the Guangdong Provincial Key Laboratory of Mathematical Foundations for Artificial Intelligence (Grant No. 2023B1212010001). 
}

\bibliography{example_paper}
\bibliographystyle{plain}

\newpage
\appendix
\onecolumn

\begin{center}
{\Large\bfseries Supplementary Material for \textit{``Bayesian Distributionally Robust Merton Problem with Nonlinear Wasserstein Projections''}\par}
\vspace{1em}
Jose Blanchet \quad Jiayi Cheng \quad Hao Liu \quad Yang Liu
\end{center}

This supplement establishes the well-posedness of the DRBC construction by
showing that its uncertainty set is independent of the prior. It also extends
the main-text results for power utility to priors with noncompact support.

{
\section{Prior-Independent Admissible Controls}\label{sec:priorind}
Recall that we place the model on the canonical product space
$\Omega=\mathbb{R}^d\times C([0,T];\mathbb{R}^d)$ equipped with the Borel
$\sigma$-algebra
$\mathcal{B}(\Omega)
=
\mathcal B(\mathbb R^d)\otimes
\mathcal B\!\bigl(C([0,T];\mathbb R^d)\bigr)$,
with canonical coordinates $B(\omega)=b$, $W(\omega)=w$, and
$\omega=(b,w)\in\Omega.$
For any prior $\mu$ on $\mathbb R^d$, define ${Q}^{\mu}:=\mu\otimes\mathbb{P}_W.$
For $\omega=(b,w)\in\Omega$, define the observation process $Y$ and the
stock-price process $S$ pathwise by
$Y(t,\omega):=\sigma^{-1}(b-r\mathbf 1)t+w(t)$, $0\le t\le T$, and $S_i(t,\omega)
=
s_{0,i}\exp\!\Big(
\bigl(r-\tfrac12(\sigma\sigma^\top)_{ii}\bigr)t
+(\sigma Y(t,\omega))_i
\Big), 0\le t\le T,$
respectively. For each $t\in[0,T]$, let both the truncated spaces $C_t^Y:=C([0,t];\mathbb R^d)$
and $C_t^S
:=
\bigl\{
s\in C([0,t];(0,\infty)^d):s(0)=s_0
\bigr\}$
be equipped with their Borel $\sigma$-fields induced by the uniform norm.
Define the truncated paths $Y^{[t]}(\omega)
:=
Y(\cdot,\omega)\big|_{[0,t]}\in C_t^Y$
and $S^{[t]}(\omega)
:=
S(\cdot,\omega)\big|_{[0,t]}\in C_t^S.$
Define the maps
\[
\Psi_t:C_t^S\to C_t^Y,
\qquad
(\Psi_t s)(u)
:=
\sigma^{-1}\!\left(
\log\Bigl(\frac{s(u)}{s(0)}\Bigr)
-
\Bigl(
r\mathbf 1-\tfrac12\operatorname{diag}(\sigma\sigma^\top)
\Bigr)u
\right),
\qquad 0\le u\le t,
\]
and
\[
\Phi_t:C_t^Y\to C_t^S,
\qquad
(\Phi_t y)_i(u)
:=
s_{0,i}\exp\!\left(
\bigl(r-\tfrac12(\sigma\sigma^\top)_{ii}\bigr)u
+(\sigma y(u))_i
\right),
\qquad 0\le u\le t.
\]
Here, all vector logarithms and divisions are understood componentwise.
Therefore, the maps $\Psi_t$ and $\Phi_t$ are Borel measurable inverses
of each other. Moreover,
$Y^{[t]}=\Psi_t(S^{[t]})$ and
$S^{[t]}=\Phi_t(Y^{[t]})$.
Therefore, the raw natural $\sigma$-fields generated by $S$ and $Y$
coincide.
Next, for any prior $\mu$ on $\mathbb R^d$, we define the law of the
observation path on $C_t^Y$ by $\nu_t^\mu
:=
{Q}^{\mu}\circ(Y^{[t]})^{-1},$
and denote the Wiener measure on $C_t^Y$ by $\gamma_t
:=
\mathbb{P}_W
\circ
\bigl(w\mapsto w|_{[0,t]}\bigr)^{-1}$.
For fixed $b\in\mathbb R^d$, define ${Q}^{\delta_b}
:=
\delta_b\otimes\mathbb{P}_W$ and $\nu_t^b
:=
{Q}^{\delta_b}\circ(Y^{[t]})^{-1},$
where $\delta_b$ is the Dirac mass at point $b$.
Then, for every Borel set $A\subseteq C_t^Y$, we have $\nu_t^\mu(A)
=
\int_{\mathbb R^d}\nu_t^b(A)\,\mu(db).$
Moreover, under ${Q}^{\delta_b}$, the process $Y$ is a Brownian
motion with deterministic drift
$u\mapsto\sigma^{-1}(b-r\mathbf 1)u$ on $[0,t]$.
Hence, $\nu_t^b\sim\gamma_t$
for every $b\in\mathbb R^d$ by the Cameron--Martin theorem.
Consequently, $\nu_t^\mu\sim\gamma_t$
for every prior $\mu$, and therefore all measures $\nu_t^\mu$ have the
same null sets on $C_t^Y$.
Let
\[
\mathfrak N_t
:=
\bigl\{
N\in\mathcal B(C_t^Y):
\nu_t^\mu(N)=0
\text{ for one, equivalently every, prior }\mu
\bigr\}
\]
be the common family of observation-null sets on the path space. Define
the corresponding null sets on $\Omega$ by
\[
\mathcal N_t^{\mathrm{obs}}
:=
\bigl\{
M\subseteq\Omega:
\exists N\in\mathfrak N_t
\text{ such that }
M\subseteq(Y^{[t]})^{-1}(N)
\bigr\}.
\]
We then define the completed observation filtration by
\[
\mathcal{F}^S(t)
:=
\mathcal F_t^{\mathrm{obs}}
=
\sigma(S(u):0\le u\le t)\vee\mathcal N_t^{\mathrm{obs}}
=
\sigma(Y(u):0\le u\le t)\vee\mathcal N_t^{\mathrm{obs}}.
\]
By construction, the filtration
$\{\mathcal F_t^{\mathrm{obs}}\}_{t\in[0,T]}$
is independent of the prior choice. If desired, one may further replace
$\mathcal F_t^{\mathrm{obs}}$ by its right-continuous modification $\mathcal F_{t+}^{\mathrm{obs}}
:=
\bigcap_{s>t}\mathcal F_s^{\mathrm{obs}}.$
}

\section{Sub-Gaussian Extensions of the Main Results}
\label{sec:subgaussian-overview}

In Sections~\ref{Th2}, \ref{Th3}, and~\ref{PCLT}, we study a commonly used
utility function in the case where the random drift $B$ is no longer
compactly supported. This shows that the compactness assumption in the
main text is not essential for the validity of our results; it is
imposed there mainly to avoid lengthy technical arguments involving
exponential moment bounds and Gaussian integrals.

Throughout these sections, we fix a CRRA utility
\begin{equation}
  \label{eq:CRRA-appendix}
  u(x) = \frac{x^{\alpha}}{\alpha}, 
  \qquad x>0,\quad \alpha<1,\ \alpha\neq 0,
\end{equation}
and work under a sub-Gaussian assumption on the prior distribution of
$B$. More precisely, we assume:

\begin{assumption}[Sub-Gaussian prior on $B$]
\label{ass:subgaussian-B}
There exists $\gamma_0>0$ such that
\[
  \E_{\mathbb{P}_0}\!\left[\exp\bigl(\gamma^2 \|B\|_2^2\bigr)\right]
  < \infty
  \qquad \text{for all } \gamma < \gamma_0 .
\]
\end{assumption}

{This condition controls the Gaussian integrals involving the likelihood
ratio $L_T(B,y)$ and its derivatives.}

\section{Generalization of Minimax Theorem with Non-compact Support}
\label{Th2}

{
We explain how the minimax theorem extends to the CRRA utility
\eqref{eq:CRRA-appendix} under Assumption~\ref{ass:subgaussian-B},
whenever the sub-Gaussian parameter gives the Gaussian moment bounds used below. We focus on the role of
compactness. The required estimates are of the same
completion-of-squares and exponential-moment type as those carried out
in Sections~\ref{Th3} and~\ref{PCLT}, and are not repeated here.
}

\medskip\noindent
{
\textbf{Step 1 (regular class).}
In the non-compact setting, we use the exponentially weighted analogue
$\mathcal A'_{\mathrm{sg}}(x_0)\subset\mathcal A(x_0)$ of the regular
class $\mathcal A'(x_0)$ used in the compact-support proof. The
slice-continuity condition (R1) is imposed on $\mathbb R^d$, while the
uniform bounds over the compact set $K$ in (R2)--(R3) are replaced by
exponentially weighted bounds. More precisely, for some sufficiently
small $\kappa>0$ and suitable $r_+>1$ and $r_->p_U$, we require
$b\mapsto\pi^b(\cdot,w)$ to be continuous from $\mathbb R^d$ into
$L^2([0,T];\mathbb R^d)$ for $\mathbb P_W$-almost every $w$, together
with bounds of the form $\mathbb E_{\mathbb P_W}\!\left[
\int_0^T
\sup_{b\in\mathbb R^d}
e^{-\kappa\|b\|_2^2}
\|\pi^b(t,W)\|_2^2\,dt
\right]<\infty$ and $\sup_{b\in\mathbb R^d}
e^{-\kappa\|b\|_2^2}
\mathbb E_{\mathbb P_W}\!\left[
(X_T^{\pi,b})^{r_+}+(X_T^{\pi,b})^{-r_-}
\right]<\infty.$
The constants may be fixed uniformly so that
$\mathcal A'_{\mathrm{sg}}(x_0)$ is convex. For
$\pi\in\mathcal A'_{\mathrm{sg}}(x_0)$, define $g_\pi(b)
:=
\mathbb E_{\mathbb P_W}
\!\left[U(X_T^{\pi,b})\right].$
}

\medskip\noindent
\textbf{Step 2 and Step 3 (moment bounds and continuity).}
{
The variation-of-constants argument from the compact-support proof,
applied on bounded neighborhoods of $b$, gives the
$L^2(\mathbb P_W)$-continuity of
$b\mapsto X_T^{\pi,b}$. The exponentially weighted positive and inverse
moment bounds give uniform integrability and an estimate of the form $|g_\pi(b)|
\leq
C_\pi e^{c_\pi\|b\|_2^2},$
where $c_\pi>0$ may be chosen within the available sub-Gaussian moment
range. Hence $g_\pi$ is continuous. The uniform exponential-moment
bound on the ambiguity set established in Step~4 below then implies
continuity of $\mu
\longmapsto
\int_{\mathbb R^d}g_\pi(b)\,\mu(db)$
under weak convergence on the ambiguity set. The same weighted moment
bounds yield the continuity-in-$\pi$ argument of Step~3 in the main
text. Thus Steps~2 and~3 have the same structure as in the
compact-support proof; only the compact uniform bounds are replaced by
exponentially weighted ones.
}

\medskip\noindent
\textbf{Step 4 (use of Sion's Theorem).}
In the sub-Gaussian extension in Sections~\ref{Th3} and~\ref{PCLT},
the optimal transport cost is no longer the quadratic cost
$\|\Delta\|_2^2$ used in the main text. Instead, following the
Gaussian-integrability estimates derived in this appendix, the natural
cost becomes 
\begin{equation}
\label{eq:exp-cost}
c_\tau(\Delta)
:=
e^{\tau\|\Delta\|_2^2}-1,
\end{equation}
where the parameter $\tau>0$ depends explicitly on
$\alpha$, $T$, and $\|\sigma^{-1}\|_F^2$ through the sub-Gaussian tail
bounds of the prior distribution of $B$.

We emphasize that this change of transport cost does \emph{not} affect
the minimax argument. The prior-level ambiguity set is $\mathcal U_{\delta,B}^{\mathrm{OT}}(\mathbb P_0)
:=
\left\{
\mu\in\mathcal P(\mathbb R^d):
D_{c_\tau}(\mu,\mathbb P_0)\leq\delta
\right\},$
and the corresponding full-model ambiguity set is $\mathcal U_{\delta}^{\mathrm{OT}}(\mathbb P_0)
:=
\left\{
Q^\mu=\mu\otimes\mathbb P_W:
\mu\in\mathcal U_{\delta,B}^{\mathrm{OT}}(\mathbb P_0)
\right\}.$
The function $\Delta\mapsto c_\tau(\Delta)$ is convex, so the
prior-level ambiguity set, and hence the induced full-model ambiguity
set, is convex.

{
The exponential transport budget also propagates the sub-Gaussian
moment of the nominal prior uniformly over the ambiguity set. 
}
Moreover, because $c_\tau(\Delta)$ grows superlinearly in
$\|\Delta\|_2$, the corresponding OT balls are tight and relatively
compact under the topology induced by optimal transport.
{
Together with the lower semicontinuity of the transport cost, this
gives weak compactness of the ambiguity set. The preceding uniform
exponential-moment bound and the growth estimate for $g_\pi$ give the
required continuity in the model variable. Since the payoff functional
is linear in $Q^\mu$, concave in the control $\pi$, and
continuous in each variable on $\mathcal A'_{\mathrm{sg}}(x_0)
\times
\mathcal U_{\delta}^{\mathrm{OT}}(\mathbb P_0),$
all assumptions of Sion's minimax theorem are satisfied.
}
Hence the change of cost from $\|\Delta\|_2^2$ to $c_\tau(\Delta)$ does
not alter the validity of the main-text min--max swap.

\medskip\noindent
\textbf{Step 5 (feedback form and admissibility).}
{
For each fixed prior
$\mu\in\mathcal U_{\delta,B}^{\mathrm{OT}}(\mathbb P_0)$, the
Karatzas--Zhao feedback construction is well defined under the same
exponential-moment conditions and solves the corresponding singleton
Bayesian problem. The same Gaussian-integral estimates used above give
the weighted analogues of (R1)--(R3): slice continuity in $b$, an
exponentially weighted integrable envelope for the sliced controls, and
weighted positive and inverse moments of the optimal terminal wealth.
Consequently, the singleton optimizer belongs to
$\mathcal A'_{\mathrm{sg}}(x_0)$. The closing argument from Step~5 of
the main-text proof therefore extends the min--max equality from
$\mathcal A'_{\mathrm{sg}}(x_0)$ to the full admissible class
$\mathcal A(x_0)$.
}
\section{Generalization of Nonlinear Perturbation Theorem with Non-compact Support}\label{Th3}
{This section extends the nonlinear perturbation theorem beyond compactly
supported priors. The estimates also justify the minimax extension in
Section~\ref{Th2}.}

After the minimax swap, for the power utility considered here, the
remaining prior-level optimization problem becomes
\begin{align}\label{FWProblem1}
\mu^{\mathrm{WC}} \in
\begin{cases}
\displaystyle
\arg\min_{\mu\in
\mathcal U_{\delta,B}^{\mathrm{OT}}(\mathbb P_0)}
\mathcal J(\mu),
& \text{if } \alpha\in(0,1),\\[2mm]
\displaystyle
\arg\max_{\mu\in
\mathcal U_{\delta,B}^{\mathrm{OT}}(\mathbb P_0)}
\mathcal J(\mu),
& \text{if } \alpha<0.
\end{cases}
\end{align}
Here
\[
\mathcal J(\mu)
:=
\int_{\mathbb R^d}
\left(
\mathbb E_{\mu}[L_T(B,z)]
\right)^{\frac{1}{1-\alpha}}
\varphi_T(z)\,dz.
\]

Throughout this section, write $\beta:=\frac{\alpha}{1-\alpha},$ $a(b):=\sigma^{-1}(b-r\mathbf 1),$ $m:=\sigma^{-1}r\mathbf 1,$ so that $a(b)=\sigma^{-1}b-m.$
{For a displacement $\Delta\in\mathbb R^d$ and scale $\tau>0$, define the
optimal transport cost by}
\[
\ c_\tau(\Delta)\;:=\;e^{\,\tau\,\|\Delta\|_2^2}\;-\;1.
\]
{This choice supplies the integrability needed to interchange integration,
differentiation, and limits without assuming compact support.}

Assume that
\[
\tau>4T\|\sigma^{-1}\|_F^2.
\]

\begin{assumption}\label{subGaussianFinal}
Suppose there exists $\gamma_0>0$ such that
\[
\E_{\mathbb P_0}
\left[
\exp\!\left(\gamma^2\|B\|_2^2\right)
\right]
<\infty
\qquad
\text{for every }\gamma<\gamma_0,
\]
and
\[
\gamma_0^2
>
T\|\sigma^{-1}\|_F^2
\max\left\{
\frac12
\left(
\beta^2-\beta+2
+
\sqrt{
(\beta^2-\beta+2)^2+8(\beta^2+\beta)
}
\right),
\,8
\right\}.
\]

When $\beta>0$, assume additionally that there exist $r_\beta>\max\{1,2\beta\}$
and conjugate exponents $p_\beta,q_\beta>1$ satisfying $\frac1{p_\beta}+\frac1{q_\beta}=1.$
Assume further that $2T(r_\beta^2-r_\beta)
\|\sigma^{-1}\|_F^2\,p_\beta
<
\gamma_0^2$ and $2T(r_\beta^2-r_\beta)
\|\sigma^{-1}\|_F^2\,q_\beta
<
\tau.$
\end{assumption}

We begin with the nonlinear perturbation theorem in this setting and stating several technical lemmas.
\begin{theorem}\label{Taylor}
Assume that for the fixed $\alpha$, we have the corresponding cost
function and $\tau$ as in the definition, and suppose
Assumption~\ref{subGaussianFinal} holds. We define
\begin{equation}\label{eq:Hb}
H(b)
:=
\frac{1}{1-\alpha}
\int_{\mathbb R^d}
\nabla_bL_T(b,y)
\left(
\E_{\mathbb P_0}[L_T(B,y)]
\right)^{\frac{\alpha}{1-\alpha}}
\varphi_T(y)\,dy
\end{equation}
and
\[
\|H\|_{L^2(\mathbb P_0)}
:=
\left(
\E_{\mathbb P_0}
\left[
\|H(B)\|_2^2
\right]
\right)^{1/2}.
\]

Suppose $\|H\|_{L^2(\mathbb P_0)}>0$.
When $\alpha\in(0,1)$, as $\delta\to0$, there exists an
asymptotically optimal feasible deterministic pushforward perturbation
such that
\[
\Delta_\delta^*(b)
=
-\sqrt{\frac{\delta}{\tau}}\,
\frac{H(b)}{\|H\|_{L^2(\mathbb P_0)}}
+
o(\sqrt{\delta})
\quad\text{in }L^2(\mathbb P_0),
\]
and
\[
\inf_{\mu\in
\mathcal U_{\delta,B}^{\mathrm{OT}}(\mathbb P_0)}
\mathcal J(\mu)
=
\mathcal J(\mathbb P_0)
-
\sqrt{\frac{\delta}{\tau}}\,
\|H\|_{L^2(\mathbb P_0)}
+
o(\sqrt{\delta}).
\]

When $\alpha<0$, as $\delta\to0$, there exists an asymptotically
optimal feasible deterministic pushforward perturbation such that
\[
\Delta_\delta^*(b)
=
\sqrt{\frac{\delta}{\tau}}\,
\frac{H(b)}{\|H\|_{L^2(\mathbb P_0)}}
+
o(\sqrt{\delta})
\quad\text{in }L^2(\mathbb P_0),
\]
and
\[
\sup_{\mu\in
\mathcal U_{\delta,B}^{\mathrm{OT}}(\mathbb P_0)}
\mathcal J(\mu)
=
\mathcal J(\mathbb P_0)
+
\sqrt{\frac{\delta}{\tau}}\,
\|H\|_{L^2(\mathbb P_0)}
+
o(\sqrt{\delta}).
\]
\end{theorem}
The proof essentially contains two parts. The first part is Lemma \ref{lem:L22-fro}, giving the required integrability conditions, and the second part is to linearize the problem and then solve the constrained linear problem, similar to the case in the main body.

\begin{lemma}\label{lem:L22-fro}
Under Assumption~\ref{subGaussianFinal},
\begin{itemize}
    \item[(1)]
    For any $\alpha<1$ and $\alpha\neq0$, the vector field in
    \eqref{eq:Hb} is well-defined $\mathbb P_0$-almost surely and
    satisfies $\E_{\mathbb P_0}
    \left[
    \|H(B)\|_2^2
    \right]
    <\infty.$

    \item[(2)]
    For any $s\in[0,1]$, let $(B,B+\Delta)$ be distributed according
    to a coupling $\pi$ whose first marginal is $\mathbb P_0$, and satisfying $\E_\pi
\left[
\exp\!\left(\tau\|\Delta\|_2^2\right)
\right]
\leq 2.$ We define
    \[
    F_s(B,\Delta)
    :=
    \int_{\mathbb R^d}
    \nabla_bL_T(B+s\Delta,y)
    \left(
    \E_{\pi}[L_T(B+s\Delta,y)]
    \right)^{\frac{\alpha}{1-\alpha}}
    \varphi_T(y)\,dy.
    \]

    Then there exist $\eta>0$, a constant $C_1>0$, and a measurable
    random variable $C(B)>0$ such that
    \[
    2(1+\eta)C_1<\tau,
    \qquad
    \E_{\mathbb P_0}
    \left[
    C(B)^{2(1+\eta)}
    \right]
    <\infty,
    \]
    and, $\pi$-almost surely,
    \[
    \|F_s(B,\Delta)\|_2^2
    \leq
    C(B)
    \exp\!\left(C_1\|\Delta\|_2^2\right),
    \qquad s\in[0,1].
    \]
\end{itemize}
\end{lemma}

\begin{proof}
We first prove part~(1). Recall that $\beta:=\frac{\alpha}{1-\alpha}.$
Define $A(y):=\E_{\mathbb P_0}[L_T(B,y)].$
We begin with the case $\alpha<0$.
Direct computations give
\[
\nabla_b L_T(b,y)
=
L_T(b,y)\,\sigma^{-T}\!\bigl(y-Ta(b)\bigr),
\]
and
\[
L_T(b,y)^2\varphi_T(y)
=
\exp\!\bigl(T\|a(b)\|_2^2\bigr)
\varphi_T\!\bigl(y-2Ta(b)\bigr).
\]
By the Jensen inequality,
\[
A(y) \;=\; \E_{\mathbb{P}_0}\big[e^{\langle a(B),y\rangle-\frac{T}{2}\|a(B)\|_2^2}\big]
\;\ge\;
\exp\!\Big( \langle \E_{\mathbb{P}_0}[a(B)],y\rangle - \tfrac{T}{2}\,\E_{\mathbb{P}_0}\left[\|a(B)\|_2^2\right] \Big).
\]
Raising to the negative power $\beta<0$ reverses the inequality and gives
\begin{equation}\label{eq:est_A_beta}
A(y)^{\beta}\;\le\;C_a\,\exp\!\big(\beta\,\langle v,y\rangle\big),
\qquad
C_a:=\exp\!\Big(-\tfrac{\beta T}{2}\,\E_{\mathbb{P}_0}\left[\|a(B)\|_2^2\right]\Big),\ \ 
v:=\E_{\mathbb{P}_0}[a(B)].
\end{equation}
Since
\[
\|\nabla_b L_T(B,y)\|_2^2
= L_T(B,y)^2\,\|\sigma^{-T}(y-T a(B))\|_{2}^2
\le \,\|\sigma^{-1}\|_F^2\,L_T(B,y)^2\,\|y-T a(B)\|_2^2,
\] 
we use $\|u+v\|_2^2\le 2(\|u\|_2^2+\|v\|_2^2)$ and have %
\begin{equation*}
\|\nabla_b L_T(B,y)\|_2^2
\;\le\; C\,\|\sigma^{-1}\|_F^2\,L_T(B,y)^2\,\big(1+\|y\|_2^2+T^2\|a(B)\|_2^2\big),
\end{equation*}where $C > 0$ is a constant (we use $C$ to absorb constants if there is no confusion). Therefore,
\[
\begin{aligned}
&\int_{\R^d} A(y)^{\beta}\,\E_{\mathbb{P}_0}\big[\|\nabla_b L_T(B,y)\|_2^2\big]\;\varphi_T(y)\,dy \\
&\ \ \le\ C\,\|\sigma^{-1}\|_F^2\,\E_{\mathbb{P}_0}\!\Big[\int_{\R^d} L_T(B,y)^2\,(1+\|y\|_2^2+T^2\|a(B)\|_2^2)\,e^{\beta\langle v,y\rangle}\,\varphi_T(y)\,dy\Big]\\
&\ \ =\ C\,\|\sigma^{-1}\|_F^2\,\E_{\mathbb{P}_0}\!\Big[e^{T\|a(B)\|_2^2}
\int_{\R^d} (1+\|y\|_2^2+T^2\|a(B)\|_2^2)\,e^{\beta\langle v,y\rangle}\,\varphi_T\!\big(y-2T a(B)\big)\,dy\Big].
\end{aligned}
\]
The inner integral is 
\begin{align*}
    &\int_{\mathbb{R}^d} \bigl(1+\|y\|_2^{2}+T^{2}\|a(B)\|_2^{2}\bigr)\,
e^{\beta\langle v,y\rangle}\,
\varphi_{T}\bigl(y-2T\,a(B)\bigr)\,dy\\
&=
\exp\!\Big(2T\beta\,\langle v,a(B)\rangle + \tfrac{T}{2}\|\beta v\|_2^{2}\Big)\,
\Big[\,1 + T d + \|2T a(B) + T\beta v\|_2^{2} + T^{2}\|a(B)\|_2^{2}\,\Big].
\end{align*}
For a fixed $\varepsilon > 0$ (which will be chosen later), by the Young inequality, 
\[
2T\beta\,\langle v,a(B)\rangle
\leq
2T|\beta|\,|\langle v,a(B)\rangle|
\leq
\varepsilon\|a(B)\|_2^2
+
\frac{T^2\beta^2\|v\|_2^2}{\varepsilon}.
\]
Hence, 
\[
e^{\,2T\beta\langle v,a(B)\rangle+\frac{T}{2}\|\beta v\|_2^2}
\le C_{\mathrm{tilt}}(\varepsilon)\,e^{\,\varepsilon\|a(B)\|_2^2},
\]
with
\[
C_{\mathrm{tilt}}(\varepsilon):=\exp\!\Big(\tfrac{T}{2}\|\beta v\|_2^2+\tfrac{T^2\beta^2\|v\|_2^2}{\varepsilon}\Big).
\]
Using the Young inequality in the quadratic bracket, we have
\[
\|2Ta+T\beta v\|_2^2\le (1+\varepsilon)\,(2T)^2\|a\|_2^2+\Big(1+\tfrac1\varepsilon\Big)T^2\beta^2\|v\|_2^2,
\]
and further get
\[
1+Td+\|2Ta(B)+T\beta v\|_2^2+T^2\|a(B)\|_2^2
\le C_0(\varepsilon)+C_1(\varepsilon)\,\|a(B)\|_2^2,
\]where $C_0(\varepsilon):=1+Td+\Big(1+\tfrac{1}{\varepsilon}\Big)T^2\beta^2\|v\|_2^2,$ and $C_1(\varepsilon):=T^2\big(5+4\varepsilon\big).$
Therefore, the inner integral is upper bounded by
$C_{\mathrm{tilt}}(\varepsilon)\,\bigl(C_0(\varepsilon)+C_1(\varepsilon)\,\|a(B)\|_2^2\bigr)\,e^{\,\varepsilon\|a(B)\|_2^2}.$

Therefore, there exists a constant $C > 0$ such that
\begin{align*}
    &\int_{\R^d} A(y)^{\beta}\,\E_{\mathbb{P}_0}\big[\|\nabla_b L_T(B,y)\|_2^2\big]\;\varphi_T(y)\,dy \\
    &\leq C\,\|\sigma^{-1}\|_F^2\,C_{\mathrm{tilt}}(\varepsilon)\;
\E_{\mathbb{P}_0}\!\Big[\bigl(C_0(\varepsilon)+C_1(\varepsilon)\|a(B)\|_2^2\bigr)\,e^{\,(T+\varepsilon)\|a(B)\|_2^2}\Big].
\end{align*}

Let
\[
\varepsilon_*:=\tfrac12\!\left(\frac{\gamma_0^2}{2\|\sigma^{-1}\|_F^2}-T\right)
\quad\text{(well-defined and $>0$ if }\gamma_0^2>2T\|\sigma^{-1}\|_F^2).
\]
Since $a(B)=\sigma^{-1}B-m$, we have $\|a(B)\|_2^2\le 2\|\sigma^{-1}\|_F^2\|B\|_2^2+2\|m\|_2^2$. Hence, 
\begin{equation}\label{goal}
\bigl(C_0(\varepsilon_*)+C_1(\varepsilon_*)\|a(B)\|_2^2\bigr)\,e^{\,(T+\varepsilon_*)\|a(B)\|_2^2}
\;\le\;K_m\,(1+c\|B\|_2^2)\,e^{\,C'\|B\|_2^2},    
\end{equation}
for constants $K_m,c>0$ and
\[
C':=2\,(T+\varepsilon_*)\,\|\sigma^{-1}\|_F^2
= T\|\sigma^{-1}\|_F^2+\frac{\gamma_0^2}{2}\;<\gamma_0^2.
\]

To see (\ref{goal}), recall that \(a(B)=\sigma^{-1}B - m\), and denote
\[
\kappa_1 := 2\,\|\sigma^{-1}\|_F^2,\qquad \kappa_0 := 2\,\|m\|_2^2.
\]
Then
$$\|a(B)\|_2^2 \;=\; \|\sigma^{-1}B - m\|_2^2
\;\le\; 2\|\sigma^{-1}B\|_2^2 + 2\|m\|_2^2
\;\le\; \kappa_1 \|B\|_2^2 + \kappa_0,$$
and
\begin{align*}
(T+\varepsilon_*)\|a(B)\|_2^2
&\le (T+\varepsilon_*)\big(\kappa_1 \|B\|_2^2 + \kappa_0\big)
= \underbrace{2(T+\varepsilon_*)\|\sigma^{-1}\|_F^2}_{=:~C'} \|B\|_2^2 \;+\; 2(T+\varepsilon_*)\|m\|_2^2,
\end{align*}
where
$C'= 2(T+\varepsilon_*)\|\sigma^{-1}\|_F^2
= T\|\sigma^{-1}\|_F^2 + \frac{\gamma_0^2}{2} \;<\; \gamma_0^2.$ 
Hence, 
$$e^{(T+\varepsilon_*)\|a(B)\|_2^2}
\;\le\; e^{\,2(T+\varepsilon_*)\|m\|_2^2}\; e^{\,C'\|B\|_2^2}.$$
For the polynomial prefactor, 
\begin{align*}
C_0(\varepsilon_*) + C_1(\varepsilon_*)\|a(B)\|_2^2
&\le C_0(\varepsilon_*) + C_1(\varepsilon_*)\big(\kappa_1 \|B\|_2^2 + \kappa_0\big) \\
&= \underbrace{\big(C_0(\varepsilon_*) + C_1(\varepsilon_*)\kappa_0\big)}_{=:~K_1}
\;+\; \underbrace{C_1(\varepsilon_*)\kappa_1}_{=:~K_2}\,\|B\|_2^2 \\
&= K_1\Big(1 + \frac{K_2}{K_1}\,\|B\|_2^2\Big)
= K_1\big(1 + c\,\|B\|_2^2\big), 
\end{align*}
where
\[
c := \frac{C_1(\varepsilon_*)\kappa_1}{C_0(\varepsilon_*) + C_1(\varepsilon_*)\kappa_0} \;>\; 0 .
\]
Therefore, we get \eqref{goal} with
\[
K_m \;:=\; K_1\, e^{\,2(T+\varepsilon_*)\|m\|_2^2}
\;=\; \big(C_0(\varepsilon_*) + C_1(\varepsilon_*)\kappa_0\big)\, e^{\,2(T+\varepsilon_*)\|m\|_2^2}.
\]
Since \(C' < \gamma_0^2\), we pick any \(\gamma\) such that \(C' < \gamma^2 < \gamma_0^2\).
For \(t\ge 0\) and \(\delta := \gamma^2 - C' > 0\), the elementary bound
\[
1 + c t \;\le\; \Big(1 + \frac{c}{e\,\delta}\Big)\, e^{\delta t}
\]
implies
\[
(1 + c\,\|B\|_2^2)\, e^{\,C'\|B\|_2^2}
\;\le\; K\, e^{\,\gamma^2 \|B\|_2^2},
\qquad K := 1 + \frac{c}{e(\gamma^2 - C')}.
\]
Therefore, under the sub-Gaussian assumption
\(\E_{\mathbb{P}_0}\big[e^{\gamma^2\|B\|_2^2}\big] < \infty\) for all \(\gamma < \gamma_0\),
\[
\E_{\mathbb{P}_0}\!\Big[\bigl(C_0(\varepsilon_*) + C_1(\varepsilon_*)\|a(B)\|_2^2\bigr)\, e^{(T+\varepsilon_*)\|a(B)\|_2^2}\Big]
\;\le\; K_m\,K\,\E_{\mathbb{P}_0}\big[e^{\gamma^2\|B\|_2^2}\big] \;<\; \infty.
\]
In order to show the well-definedness of $H$, we need to show
\begin{equation}\label{CS}
    \int_{\mathbb{R}^d}A(y)^{\beta}\varphi_T(y)dy < \infty
\end{equation}since from the Cauchy-Schwarz inequality, for $\mathbb{P}_0$ almost every $b$,
\begin{align*}
    &\left\Vert H(b) \right\Vert_2 \leq \frac{1}{1-\alpha}\int_{\mathbb R^d}A(y)^{\beta}\left\Vert\nabla_b L_T(b,y)\right\Vert_2\,\,\varphi_T(y)\,dy\\
    &\leq \frac{1}{1-\alpha} \left(\int_{\mathbb{R}^d}\left\Vert\nabla_b L_T(b,y)\right\Vert_2^2A(y)^{\beta}\,\varphi_T(y)dy\right)^{\frac{1}{2}}\left(\int_{\mathbb{R}^d}A(y)^{\beta}\,\varphi_T(y)dy\right)^{\frac{1}{2}},
\end{align*}
{The first term is finite $\mathbb P_0$-almost surely. Squaring, taking
$\mathbb P_0$-expectations, and applying Tonelli's theorem complete the proof.}

When $\alpha < 0$, we recall Eq. \eqref{eq:est_A_beta} and have %
\[
\int_{\R^d} A(y)^{\beta}\,\varphi_T(y)\,dy
\;\le\;
C_a \int_{\R^d} e^{\beta\langle v,y\rangle}\,\varphi_T(y)\,dy.
\]
Completion of squares gives
\[
\int_{\R^d} e^{\beta\langle v,y\rangle}\,\varphi_T(y)\,dy
=(2\pi T)^{-d/2}\!\!\int_{\R^d}
\exp\!\left(\beta\langle v,y\rangle-\frac{\|y\|_2^2}{2T}\right)dy
=\exp\!\left(\frac{T}{2}\|\beta v\|_2^2\right).
\]
Indeed,
\[
\beta\langle v,y\rangle-\frac{\|y\|_2^2}{2T}
=-\frac{1}{2T}\Big\|\,y-T\beta v\,\Big\|_2^2+\frac{T}{2}\|\beta v\|_2^2.
\]
Hence, 
\[
\int_{\R^d} A(y)^{\beta}\,\varphi_T(y)\,dy
\;\le\;
C_a\,\exp\!\left(\frac{T}{2}\|\beta v\|_2^2\right)
<\infty.
\]
This finishes the proof of the case when $\alpha<0$.
Next we divide the positive $\alpha$ into two cases.

\noindent\emph{Case 1: $0<\beta\le 1$ (equivalently $0<\alpha\le \tfrac12$).}
Let $\lambda_T(dy):=\varphi_T(y)\,dy$; then $\lambda_T$ is a probability measure since
\[
\int_{\R^d}\varphi_T(y)\,dy=1.
\]
Set $f(y):=A(y)=\E_{\mathbb P_0}[L_T(B,y)]\geq0$ and
\[
\phi(x):=x^{\beta},\qquad x\ge 0.
\]
For $0<\beta\le 1$, the map $\phi$ is concave and increasing, so the Jensen 
inequality for concave functions gives
\[
\phi\!\left(\int_{\R^d} f(y)\,\lambda_T(dy)\right)\;\ge\;\int_{\R^d}\phi(f(y))\,\lambda_T(dy),
\]
i.e.
\begin{equation}\label{eq:jensen-easy}
\int_{\R^d} A(y)^{\beta}\,\varphi_T(y)\,dy
\;\le\;
\Big(\int_{\R^d} A(y)\,\varphi_T(y)\,dy\Big)^{\beta}.
\end{equation}

It remains to compute $\int_{\R^d} A(y)\,\varphi_T(y)\,dy$. By the Fubini theorem,
\[
\int_{\R^d} A(y)\,\varphi_T(y)\,dy
=\E_{\mathbb{P}_0}\!\left[\int_{\R^d} L_T(B,y)\,\varphi_T(y)\,dy\right] = 1
\]
from a completion of squares argument. Thus, plugging into \eqref{eq:jensen-easy} yields the bound
\[\ \int_{\R^d} A(y)^{\beta}\,\varphi_T(y)\,dy \;\le\; 1,\qquad 0<\beta\le 1. \ 
\]

\medskip
\noindent\emph{Case 2: $\beta>1$ (equivalently $\tfrac12<\alpha<1$).}
Since the map $x\mapsto x^{\beta}$ is convex on $(0,\infty)$, the Jensen inequality yields
\[
A(y)^{\beta}=\big(\E_{\mathbb{P}_0}[L_T(B,y)]\big)^{\beta}\le \E_{\mathbb{P}_0}\big[L_T(B,y)^{\beta}\big].
\]
Integrating and exchanging expectation and integral,
\[
\int_{\R^d} A(y)^{\beta}\,\varphi_T(y)\,dy
\;\le\;
\E_{\mathbb{P}_0}\!\left[\int_{\R^d} L_T(B,y)^{\beta}\,\varphi_T(y)\,dy\right].
\]
For fixed $b$, completing the square gives
\[
\int_{\R^d} L_T(b,y)^{\beta}\,\varphi_T(y)\,dy
=\exp\!\left(\frac{T}{2}(\beta^2-\beta)\,\|a(b)\|_2^2\right).
\]
Therefore
\[
\int_{\R^d} A(y)^{\beta}\,\varphi_T(y)\,dy
\;\le\;
\E_{\mathbb{P}_0}\!\left[\exp\!\left(\frac{T}{2}(\beta^2-\beta)\,\|a(B)\|_2^2\right)\right].
\]
Using $\|a(B)\|_2^2\le 2\|\sigma^{-1}\|_F^2\,\|B\|_2^2+2\|m\|_2^2$, we obtain
\begin{equation}\label{eq:est_int_A_beta}
\int_{\R^d} A(y)^{\beta}\,\varphi_T(y)\,dy
\;\le\;
\exp\!\Big(T(\beta^2-\beta)\,\|m\|_2^2\Big)\;
\E_{\mathbb{P}_0}\!\left[\exp\!\Big(T(\beta^2-\beta)\,\|\sigma^{-1}\|_F^2\,\|B\|_2^2\Big)\right].
\end{equation}
Hence, under the sub-Gaussian moment assumption
\[
\E_{\mathbb{P}_0}\!\left[e^{\gamma^2\|B\|_2^2}\right]<\infty\quad\text{for all }\gamma<\gamma_0,
\]
the right-hand side of Eq. \eqref{eq:est_int_A_beta} is finite whenever
\begin{equation}\label{eq:easy-threshold}
\ T(\beta^2-\beta)\,\|\sigma^{-1}\|_F^2 \;<\; \gamma_0^2.
\end{equation}

Next, we focus on the bound when $\alpha\in(0,1)$. Define
\[
\mathcal I_\beta
:=
\int_{\R^d}
A(y)^\beta
\E_{\mathbb P_0}
\left[
\|\nabla_bL_T(B,y)\|_2^2
\right]
\varphi_T(y)\,dy.
\]
Write
\[
G(y)
:=
\E_{\mathbb P_0}
\left[
\|\nabla_bL_T(B,y)\|_2^2
\right],
\]
and set
\[
I_0
:=
\int_{\R^d}G(y)\varphi_T(y)\,dy,
\qquad
I_1
:=
\int_{\R^d}A(y)G(y)\varphi_T(y)\,dy.
\]

When $0<\beta<1$, H\"older's inequality with exponents
$p=1/\beta$ and $q=1/(1-\beta)$ gives
\begin{equation}\label{eq:Ibeta-interpolate}
\begin{aligned}
\mathcal I_\beta
&=
\int_{\R^d}
\bigl(A(y)G(y)\bigr)^\beta
G(y)^{1-\beta}
\varphi_T(y)\,dy \\
&\leq
\left(
\int_{\R^d}A(y)G(y)\varphi_T(y)\,dy
\right)^\beta
\left(
\int_{\R^d}G(y)\varphi_T(y)\,dy
\right)^{1-\beta} \\
&=
I_1^\beta I_0^{1-\beta}.
\end{aligned}
\end{equation}
When $\beta=1$, we simply have
\[
\mathcal I_1=I_1.
\]
Thus, it suffices to show that
\[
I_0<\infty
\qquad\text{and}\qquad
I_1<\infty.
\]

Recall that there exists a constant $C > 0$ such that
\[
\|\nabla_b L_T(B,y)\|_2^2
\;\le\;C\,\|\sigma^{-1}\|_F^2\,L_T(B,y)^2\,(1+\|y\|_2^2+T^2\|a(B)\|_2^2).
\]
Using the Tonelli theorem,
\[
\begin{aligned}
I_0
&\le C\,\|\sigma^{-1}\|_F^2\,\E_{\mathbb{P}_0}\!\Big[e^{T\|a(B)\|_2^2}
\!\!\int_{\R^d}\!(1+\|y\|_2^2+T^2\|a(B)\|_2^2)\,\varphi_T(y-2Ta(B))\,dy\Big].
\end{aligned}
\]For the inner integral, 
fix \(a\in\R^d\) and set the notation
\[
J(a):=\int_{\R^d}\bigl(1+\|y\|_2^2+T^2\|a\|_2^2\bigr)\,\varphi_T(y-2Ta)\,dy.
\]
Let \(Z\sim\mathcal N(0,T I_d)\). Since \(\varphi_T(y-2Ta)\,dy\) is the law of \(Y:=Z+2Ta\), we have 
\[
J(a)=\E\!\left[1+\|Y\|_2^2+T^2\|a\|_2^2\right]
=1+\E\left[\|Z+2Ta\|_2^2\right]+T^2\|a\|_2^2 = 1+dT+5T^2\|a\|_2^2.
\]
With \(a=a(B)\), the inner integral in \(I_0\) equals \(1+dT+5T^2\|a(B)\|_2^2\), and therefore
\[
I_0 \le C\,\|\sigma^{-1}\|_F^2\,
\E_{\mathbb{P}_0}\!\Big[e^{T\|a(B)\|_2^2}\big(1+dT+5T^2\|a(B)\|_2^2\big)\Big].
\]
The estimate for \(I_0\) requires
\[
\mathbb{E}_{\mathbb{P}_0}\!\left[\exp\!\big(T\|a(B)\|_2^{2}\big)\right] < \infty,
\]
which is actually a sub-Gaussian assumption on \(B\) with condition
\[
2T\,\|\sigma^{-1}\|_F^{2} \;<\; \gamma_0^{2}.
\]
Next, we upper bound $I_1$. By the Young inequality \(\langle u,v\rangle\le \tfrac{\varepsilon}{2}\|u\|_2^2+\tfrac{1}{2\varepsilon}\|v\|_2^2\) with a fixed $\varepsilon > 0$ (which will be chosen later),
\[
A(y)\;\le\;\E_{\mathbb{P}_0}\!\Big[e^{-\frac{T-\varepsilon}{2}\|a(B)\|_2^2}\Big]\;
\exp\!\Big(\frac{\|y\|_2^2}{2\varepsilon}\Big)
\;=\;C_A(\varepsilon)\,e^{c\|y\|_2^2},
\qquad c:=\frac{1}{2\varepsilon}\in\Big(0,\frac{1}{2T}\Big).
\]
Here \(C_A(\varepsilon)=\E_{\mathbb{P}_0}\!\big[e^{\frac{\varepsilon-T}{2}\|a(B)\|_2^2}\big]\).
Hence, by the Tonelli theorem,
\[
\begin{aligned}
I_1
&\le C\,\|\sigma^{-1}\|_F^2\,C_A(\varepsilon)\,
\E_{\mathbb{P}_0}\!\Big[e^{T\|a(B)\|_2^2}\!\int_{\R^d}(1+\|y\|_2^2+T^2\|a(B)\|_2^2)e^{c\|y\|_2^2}\varphi_T(y-2Ta(B))\,dy\Big].
\end{aligned}
\]
Fix \(a\in\R^d\) and \(c\in(0,\frac{1}{2T})\). Let
\[
K_c(a):=\int_{\R^d}\bigl(1+\|y\|_2^2+T^2\|a\|_2^2\bigr)\,e^{c\|y\|_2^2}\,\varphi_T(y-2Ta)\,dy.
\]
With \(Y\sim\mathcal N(2Ta,TI_d)\) we have \(K_c(a)=\E[(1+\|Y\|_2^2+T^2\|a\|_2^2)e^{c\|Y\|_2^2}]\).
Define
\[
F(c,a):=\int_{\R^d} e^{c\|y\|_2^2}\,\varphi_T(y-2Ta)\,dy
=(1-2cT)^{-d/2}\exp\!\Big(\frac{4T^2 c}{1-2cT}\,\|a\|_2^2\Big).
\]
Then
\[
\frac{\partial}{\partial c}F(c,a)
=\int_{\R^d}\|y\|_2^2 e^{c\|y\|_2^2}\varphi_T(y-2Ta)\,dy
=F(c,a)\!\left(\frac{dT}{1-2cT}+\frac{4T^2\|a\|_2^2}{(1-2cT)^2}\right),
\]
so
\[
\begin{aligned}
K_c(a)
&=(1+T^2\|a\|_2^2)\,F(c,a)+\frac{\partial}{\partial c}F(c,a)\\
&=F(c,a)\!\left[1+T^2\|a\|_2^2+\frac{dT}{1-2cT}+\frac{4T^2\|a\|_2^2}{(1-2cT)^2}\right]
\;\le\;C_1(c,T,d)\,\bigl(1+\|a\|_2^2\bigr)\,F(c,a).
\end{aligned}
\]
Therefore, with \(c=\frac{1}{2\varepsilon}\in(0,\frac{1}{2T})\),
\[
\begin{aligned}
I_1
&\le C\,\|\sigma^{-1}\|_F^2\,C_A(\varepsilon)\,
\E_{\mathbb{P}_0}\!\Big[e^{T\|a(B)\|_2^2}\,K_c\big(a(B)\big)\Big]\\
&\le C(\varepsilon,T,d,\sigma)\,
\E_{\mathbb{P}_0}\!\Big[(1+\|a(B)\|_2^2)\,
\exp\!\Big(\Big(T+\frac{2T^2}{\varepsilon-T}\Big)\|a(B)\|_2^2\Big)\Big],
\end{aligned}
\]where \[
C(\varepsilon,T,d,\sigma)
= C\;\|\sigma^{-1}\|_{F}^{2}\;C_{A}(\varepsilon)\;
\Big(1-\tfrac{T}{\varepsilon}\Big)^{-d/2}
\left(1+\frac{dT}{1-\tfrac{T}{\varepsilon}}+T^{2}
+\frac{4T^{2}}{\big(1-\tfrac{T}{\varepsilon}\big)^{2}}\right).
\]
To make the bound for \(I_1\) finite, we only need an \(\varepsilon>T\) such that 
\[
\frac{\varepsilon-T}{2}\;<\;\frac{\gamma_0^2}{2\|\sigma^{-1}\|_F^2}
\qquad\text{and}\qquad
T+\frac{2T^2}{\varepsilon-T}\;<\;\frac{\gamma_0^2}{2\|\sigma^{-1}\|_F^2}.
\]
These two inequalities simultaneously hold whenever
\[
4T^{2}\,\|\sigma^{-1}\|_{F}^{4}\;+\;2T\,\|\sigma^{-1}\|_{F}^{2}\,\gamma_{0}^{2}\;<\;\gamma_{0}^{4},
\]which is equivalent to 
\[
 \gamma_0^2 \;>\; (1+\sqrt{5})\,T\,\|\sigma^{-1}\|_F^2.
\]
Under this %
condition%
, an explicit admissible choice is
\[
\varepsilon^{*}\;:=\;T+\frac{1}{2}\left(
\frac{2T^{2}}{\frac{\gamma_{0}^{2}}{2\|\sigma^{-1}\|_{F}^{2}}-T}
+\frac{\gamma_{0}^{2}}{\|\sigma^{-1}\|_{F}^{2}}
\right)\;>\;T.
\]
For this \(\varepsilon^{*}\) we have
\[
C_{A}(\varepsilon^{*})
=\mathbb{E}_{\mathbb{P}_{0}}\!\left[\exp\!\left(\frac{\varepsilon^{*}-T}{2}\,\|a(B)\|_{2}^{2}\right)\right]<\infty,
\qquad
\mathbb{E}_{\mathbb{P}_0}\!\left[(1+\|a(B)\|_2^2)\,
\exp\!\left(\Big(T+\frac{2T^2}{\varepsilon^{*}-T}\Big)\|a(B)\|_2^2\right)\right]<\infty,
\]
by the sub-Gaussian integrability above. Hence,  \(I_1<\infty\).

Finally, we upper bound $\mathcal I_\beta$ when $\beta > 1$.
By convexity of the map $x \mapsto x^{\beta}$ when $\beta > 1$ and the Jensen inequality, 
\[
A(y)^{\beta}\;=\;\big(\mathbb E_{\mathbb P_0}[L_T(B,y)]\big)^{\beta}\;\le\;\mathbb E_{\mathbb P_0}\!\big[L_T(B,y)^{\beta}\big].
\]
Fix \(\varepsilon>\beta^2 T\) and use the Young inequality
\(\langle u,v\rangle\le \tfrac{\varepsilon}{2}\|u\|_2^2+\tfrac{1}{2\varepsilon}\|v\|_2^2\) with \(u=\beta a(B)\), \(v=y\):
\[
L_T(B,y)^{\beta}
=\exp\!\Big(\beta\langle a(B),y\rangle-\tfrac{\beta T}{2}\|a(B)\|_2^2\Big)
\le \exp\!\Big(\tfrac{\varepsilon-\beta T}{2}\|a(B)\|_2^2\Big)\,e^{c\|y\|_2^2},
\quad c:=\frac{\beta^2}{2\varepsilon}\in\Big(0,\frac{1}{2T}\Big).
\]
Taking expectation in \(B\) yields
\[
A(y)^{\beta}\;\le\;C_{\beta}(\varepsilon)\,e^{c\|y\|_2^2},
\quad
C_{\beta}(\varepsilon):=
\mathbb E_{\mathbb P_0}\!\Big[\exp\!\Big(\tfrac{\varepsilon-\beta T}{2}\|a(B)\|_2^2\Big)\Big].
\]
Thus, exactly as in the \(I_1\) bound for \(\beta\in(0,1)\), we obtain 
\[
\mathcal I_\beta
\;\le\;C(\varepsilon,T,d,\sigma)\;
\mathbb E_{\mathbb P_0}\!\Big[(1+\|a(B)\|_2^2)\,
\exp\!\Big(\Big(T+\frac{2T^2\beta^2}{\varepsilon-\beta^2T}\Big)\|a(B)\|_2^2\Big)\Big],
\]
for some finite constant \(C(\varepsilon,T,d,\sigma)\) whenever \(c\in(0,1/(2T))\).
Therefore, to have \(\mathcal I_\beta<\infty\), it suffices to pick \(\varepsilon>\beta^2T\) such that
\[
\frac{\varepsilon-\beta T}{2}\;<\;\frac{\gamma_0^2}{2\|\sigma^{-1}\|_F^2}
\qquad\text{and}\qquad
T+\frac{2T^2\beta^2}{\varepsilon-\beta^2T}\;<\;\frac{\gamma_0^2}{2\|\sigma^{-1}\|_F^2}.
\tag{$\star$}
\]
These two inequalities simultaneously hold precisely when
\[
\frac{\gamma_0^2}{\|\sigma^{-1}\|_F^2}\;>\;\frac{T}{2}\Big(\beta^2-\beta+2
\;+\;\sqrt{\,(\beta^2-\beta+2)^{2}+8(\beta^2+\beta)\,}\Big)\;.
\]
Under the displayed condition, choose any
\[
\varepsilon\in\Big(\,\beta^2T+\frac{2T^2\beta^2}{\frac{\gamma_0^2}{2\|\sigma^{-1}\|_F^2}-T}\,,\;\beta T+\frac{\gamma_0^2}{\|\sigma^{-1}\|_F^2}\,\Big),
\]
which is a nonempty interval; then both parts of \((\star)\) hold and the expectation above is finite. Hence \(\mathcal I_\beta<\infty\) for all \(\beta>1\) under the stated sub-Gaussian slack.

For the proof of part~(2), let $\pi$ be a coupling whose first marginal
is $\mathbb P_0$, write $(B,B+\Delta)\sim\pi,$
and suppose that $\E_\pi
\left[
\exp\!\left(\tau\|\Delta\|_2^2\right)
\right]
\leq 2.$
Throughout this part, expectations involving both $B$ and $\Delta$ are
taken with respect to $\pi$.

For $s\in[0,1]$, define
\[
I_\beta^\pi(s)
:=
\int_{\R^d}
\left(
\E_\pi[L_T(B+s\Delta,y)]
\right)^{2\beta}
\varphi_T(y)\,dy.
\]
By the Cauchy--Schwarz inequality in $y$,
\begin{align}
\|F_s(B,\Delta)\|_2^2
&\leq
\left(
\int_{\R^d}
\|\nabla_bL_T(B+s\Delta,y)\|_2^2
\varphi_T(y)\,dy
\right)
I_\beta^\pi(s).
\label{eq:Fs-two-factors}
\end{align}

We first bound the first factor in \eqref{eq:Fs-two-factors}. Recall that
\[
\|\nabla_bL_T(b,y)\|_2^2
\leq
\|\sigma^{-1}\|_F^2
L_T(b,y)^2
\left(
1+\|y\|_2^2+T^2\|a(b)\|_2^2
\right).
\]
Using
\[
L_T(b,y)^2\varphi_T(y)
=
\exp\!\left(T\|a(b)\|_2^2\right)
\varphi_T(y-2Ta(b))
\]
and the corresponding Gaussian moment formula, we obtain
\[
\int_{\R^d}
\|\nabla_bL_T(B+s\Delta,y)\|_2^2
\varphi_T(y)\,dy
\leq
C
\left(
1+\|a(B+s\Delta)\|_2^2
\right)
\exp\!\left(
T\|a(B+s\Delta)\|_2^2
\right).
\]
Since
\[
a(B+s\Delta)=a(B)+s\sigma^{-1}\Delta,
\]
we have
\[
\|a(B+s\Delta)\|_2^2
\leq
2\|a(B)\|_2^2
+
2\|\sigma^{-1}\|_F^2\|\Delta\|_2^2.
\]
Consequently,
\begin{align}
&\int_{\R^d}
\|\nabla_bL_T(B+s\Delta,y)\|_2^2
\varphi_T(y)\,dy \nonumber\\
&\quad\leq
C
\left(
1+2\|a(B)\|_2^2
+2\|\sigma^{-1}\|_F^2\|\Delta\|_2^2
\right)
\exp\!\left(
2T\|a(B)\|_2^2
+2T\|\sigma^{-1}\|_F^2\|\Delta\|_2^2
\right).
\label{eq:first-factor-pre-envelope}
\end{align}

By Assumption~\ref{subGaussianFinal}, $\tau>4T\|\sigma^{-1}\|_F^2$ and $\gamma_0^2>8T\|\sigma^{-1}\|_F^2.$
We may therefore choose $\eta>0$ sufficiently small such that
\[
4(1+\eta)T\|\sigma^{-1}\|_F^2<\tau
\]
and
\[
8(1+\eta)T\|\sigma^{-1}\|_F^2<\gamma_0^2.
\]
Next, choose $\varepsilon_\Delta>0$ sufficiently small that
\[
2(1+\eta)
\left(
2T\|\sigma^{-1}\|_F^2+\varepsilon_\Delta
\right)
<\tau,
\]
and define
\[
C_1
:=
2T\|\sigma^{-1}\|_F^2+\varepsilon_\Delta.
\]
For $x\geq0$,
\[
1+2\|\sigma^{-1}\|_F^2x
\leq
K_\Delta e^{\varepsilon_\Delta x}
\]
for a finite constant $K_\Delta$. Also,
\[
1+2\|a(B)\|_2^2
+2\|\sigma^{-1}\|_F^2\|\Delta\|_2^2
\leq
\left(
1+2\|a(B)\|_2^2
\right)
\left(
1+2\|\sigma^{-1}\|_F^2\|\Delta\|_2^2
\right).
\]
It follows from \eqref{eq:first-factor-pre-envelope} that
\begin{equation}\label{eq:first-factor-envelope}
\int_{\R^d}
\|\nabla_bL_T(B+s\Delta,y)\|_2^2
\varphi_T(y)\,dy
\leq
C_0(B)
\exp\!\left(C_1\|\Delta\|_2^2\right),
\end{equation}
where
\[
C_0(B)
:=
CK_\Delta
\left(
1+2\|a(B)\|_2^2
\right)
\exp\!\left(
2T\|a(B)\|_2^2
\right).
\]

Using
\[
\|a(B)\|_2^2
\leq
2\|\sigma^{-1}\|_F^2\|B\|_2^2
+2\|m\|_2^2,
\]
we obtain, for a finite constant $K$,
\[
C_0(B)^{2(1+\eta)}
\leq
K
\left(
1+\|B\|_2^2
\right)^{2(1+\eta)}
\exp\!\left(
8(1+\eta)T
\|\sigma^{-1}\|_F^2
\|B\|_2^2
\right).
\]
The strict inequality
\[
8(1+\eta)T\|\sigma^{-1}\|_F^2<\gamma_0^2
\]
and the sub-Gaussian assumption therefore imply
\begin{equation}\label{eq:C0-moment}
\E_{\mathbb P_0}
\left[
C_0(B)^{2(1+\eta)}
\right]
<\infty.
\end{equation}

It remains to bound $I_\beta^\pi(s)$ uniformly over the couplings
considered above.

Suppose first that $\beta<0$. Define
\[
a_s(B,\Delta)
:=
a(B)+s\sigma^{-1}\Delta.
\]
By Jensen's inequality,
\begin{align*}
\E_\pi[L_T(B+s\Delta,y)]
&=
\E_\pi
\left[
\exp\!\left(
\langle a_s(B,\Delta),y\rangle
-
\frac{T}{2}\|a_s(B,\Delta)\|_2^2
\right)
\right] \\
&\geq
\exp\!\left(
\left\langle
\E_\pi[a_s(B,\Delta)],y
\right\rangle
-
\frac{T}{2}
\E_\pi
\left[
\|a_s(B,\Delta)\|_2^2
\right]
\right).
\end{align*}
Since $2\beta<0$, raising both sides to the power $2\beta$ reverses the
inequality. Integrating the resulting bound against $\varphi_T$ gives
\[
I_\beta^\pi(s)
\leq
\exp\!\left(
(2\beta^2-\beta)T
\E_\pi
\left[
\|a_s(B,\Delta)\|_2^2
\right]
\right).
\]
Moreover,
\[
\E_\pi
\left[
\|a_s(B,\Delta)\|_2^2
\right]
\leq
2\E_{\mathbb P_0}
\left[
\|a(B)\|_2^2
\right]
+
2\|\sigma^{-1}\|_F^2
\E_\pi
\left[
\|\Delta\|_2^2
\right].
\]
By Jensen's inequality and the exponential-cost bound,
\[
\exp\!\left(
\tau\E_\pi\|\Delta\|_2^2
\right)
\leq
\E_\pi
\left[
\exp\!\left(
\tau\|\Delta\|_2^2
\right)
\right]
\leq 2,
\]
and therefore
\[
\E_\pi\|\Delta\|_2^2
\leq
\frac{\log 2}{\tau}.
\]
It follows that there exists a finite constant $K_{\beta,-}$,
independent of $\pi$ and $s$, such that
\begin{equation}\label{eq:Ibeta-negative-uniform}
I_\beta^\pi(s)
\leq
K_{\beta,-},
\qquad
\beta<0.
\end{equation}

Now suppose that $\beta>0$. Choose
$r_\beta,p_\beta,q_\beta$ as in
Assumption~\ref{subGaussianFinal}. Thus,
\[
r_\beta>\max\{1,2\beta\},
\qquad
\frac{1}{p_\beta}+\frac{1}{q_\beta}=1,
\]
and
\[
2T(r_\beta^2-r_\beta)
\|\sigma^{-1}\|_F^2p_\beta
<
\gamma_0^2,
\]
while
\[
2T(r_\beta^2-r_\beta)
\|\sigma^{-1}\|_F^2q_\beta
<
\tau.
\]
By the Lyapunov and Jensen inequalities,
\[
I_\beta^\pi(s)
\leq
\left(
\E_\pi
\left[
\int_{\R^d}
L_T(B+s\Delta,y)^{r_\beta}
\varphi_T(y)\,dy
\right]
\right)^{2\beta/r_\beta}.
\]
Completing the square gives
\[
\int_{\R^d}
L_T(b,y)^{r_\beta}
\varphi_T(y)\,dy
=
\exp\!\left(
\frac{T}{2}
(r_\beta^2-r_\beta)
\|\sigma^{-1}b-m\|_2^2
\right).
\]
Since
\[
\|\sigma^{-1}(B+s\Delta)-m\|_2^2
\leq
4\|\sigma^{-1}\|_F^2
\left(
\|B\|_2^2+\|\Delta\|_2^2
\right)
+2\|m\|_2^2,
\]
{H\"older's inequality} yields
\begin{align*}
I_\beta^\pi(s)
&\leq
K_{\beta,+}
\left(
\E_{\mathbb P_0}
\left[
\exp\!\left(
2T(r_\beta^2-r_\beta)
\|\sigma^{-1}\|_F^2p_\beta
\|B\|_2^2
\right)
\right]
\right)^{\frac{2\beta}{r_\beta p_\beta}} \\
&\quad\times
\left(
\E_\pi
\left[
\exp\!\left(
2T(r_\beta^2-r_\beta)
\|\sigma^{-1}\|_F^2q_\beta
\|\Delta\|_2^2
\right)
\right]
\right)^{\frac{2\beta}{r_\beta q_\beta}},
\end{align*}
where $K_{\beta,+}<\infty$ is independent of $\pi$ and $s$.
The first expectation is finite by the sub-Gaussian assumption. The
coefficient in the second expectation is strictly smaller than $\tau$,
and hence
\[
\E_\pi
\left[
\exp\!\left(
2T(r_\beta^2-r_\beta)
\|\sigma^{-1}\|_F^2q_\beta
\|\Delta\|_2^2
\right)
\right]
\leq
\E_\pi
\left[
\exp\!\left(
\tau\|\Delta\|_2^2
\right)
\right]
\leq 2.
\]
Therefore, there exists a finite constant $\widetilde K_{\beta,+}$,
independent of $\pi$ and $s$, such that
\begin{equation}\label{eq:Ibeta-positive-uniform}
I_\beta^\pi(s)
\leq
\widetilde K_{\beta,+},
\qquad
\beta>0.
\end{equation}

Combining \eqref{eq:Ibeta-negative-uniform} and
\eqref{eq:Ibeta-positive-uniform}, there exists a finite constant
$K_\beta$, independent of $\pi$ and $s$, such that $I_\beta^\pi(s)\leq K_\beta.$
Together with \eqref{eq:Fs-two-factors} and
\eqref{eq:first-factor-envelope}, this gives
\[
\|F_s(B,\Delta)\|_2^2
\leq
K_\beta C_0(B)
\exp\!\left(
C_1\|\Delta\|_2^2
\right),
\qquad
s\in[0,1],
\]
$\pi$-almost surely.

Finally, set $C(B):=K_\beta C_0(B).$
Then \eqref{eq:C0-moment} gives $\E_{\mathbb P_0}
\left[
C(B)^{2(1+\eta)}
\right]
<\infty,$
and by construction, $2(1+\eta)C_1<\tau.$
This proves part~(2).
\end{proof}

\begin{proof}[Proof of Theorem~\ref{Taylor}]
We first prove the result for $\alpha<0$, in which case the
prior-level problem is $\sup_{\mu\in
\mathcal U_{\delta,B}^{\mathrm{OT}}(\mathbb P_0)}
\mathcal J(\mu).$

For $\delta>0$, define the collection of couplings
\[
\mathcal C_\delta
:=
\left\{
\pi\in\mathcal P(\mathbb R^d\times\mathbb R^d):
\pi(\cdot,\mathbb R^d)=\mathbb P_0,\ 
\int_{\mathbb R^d\times\mathbb R^d}
\left(e^{\tau\|x-y\|_2^2}-1\right)\pi(dx,dy)\le\delta
\right\}.
\]
For $\pi\in\mathcal C_\delta$, write the coupled pair as $(B,B+\Delta)\sim\pi.$
Thus $B\sim\mathbb P_0$, while the law of $B+\Delta$ is the
corresponding perturbed prior. Throughout the proof, expectations
involving both $B$ and $\Delta$ are taken with respect to $\pi$.

By Jensen's inequality,
\[
\exp\!\left(
\tau\E_\pi\|\Delta\|_2^2
\right)
\le
\E_\pi
\left[
\exp\!\left(
\tau\|\Delta\|_2^2
\right)
\right]
\le
1+\delta.
\]
Therefore,
\begin{equation}\label{eq:Delta-L2-small}
\E_\pi\|\Delta\|_2^2
\le
\frac{\log(1+\delta)}{\tau}
=
\frac{\delta}{\tau}
+
o(\delta),
\qquad
\delta\downarrow0.
\end{equation}

For $s\in[0,1]$, define
\[
G_s^\pi(B,\Delta)
:=
\frac{1}{1-\alpha}
\int_{\mathbb R^d}
\nabla_bL_T(B+s\Delta,y)
\left(
\E_\pi[L_T(B+s\Delta,y)]
\right)^{\frac{\alpha}{1-\alpha}}
\varphi_T(y)\,dy.
\]
By the integral form of Taylor's theorem and Fubini's theorem, whose
use is justified by Lemma~\ref{lem:L22-fro},
\begin{align}\label{eq:J-coupling-expansion}
\mathcal J\bigl(\mathcal L_\pi(B+\Delta)\bigr)
-
\mathcal J(\mathbb P_0)
=
\int_0^1
\E_\pi
\left[
\left\langle
G_s^\pi(B,\Delta),\Delta
\right\rangle
\right]
ds.    
\end{align}
Consequently, by Cauchy--Schwarz,
\begin{align}\label{eq:J-upper-cs}
\left|
\mathcal J\bigl(\mathcal L_\pi(B+\Delta)\bigr)
-
\mathcal J(\mathbb P_0)
\right|
\le
\left(
\E_\pi\|\Delta\|_2^2
\right)^{1/2}
\int_0^1
\left(
\E_\pi
\left[
\|G_s^\pi(B,\Delta)\|_2^2
\right]
\right)^{1/2}
ds.
\end{align}

We next show that the integral on the right-hand side converges to
$\|H\|_{L^2(\mathbb P_0)}$, uniformly over small-cost couplings. Let
$\delta_n\downarrow0$, let $\pi_n\in\mathcal C_{\delta_n}$ be arbitrary,
and write $(B,B+\Delta_n)\sim\pi_n.$
By \eqref{eq:Delta-L2-small}, $\E_{\pi_n}\|\Delta_n\|_2^2\longrightarrow0.$
For $s\in[0,1]$, set $X_n(s)
:=
\left\|
G_s^{\pi_n}(B,\Delta_n)
\right\|_2^2.$
For all sufficiently large $n$, we have $\delta_n\le1$, and therefore $\E_{\pi_n}
\left[
e^{\tau\|\Delta_n\|_2^2}
\right]
\le
1+\delta_n
\le2.$
Lemma~\ref{lem:L22-fro}, with the fixed factor
$(1-\alpha)^{-2}$ absorbed into $C(B)$, gives constants
$\eta>0$ and $C_1>0$ and a measurable random variable $C(B)>0$
such that
\[
2(1+\eta)C_1<\tau,
\qquad
\E_{\mathbb P_0}
\left[
C(B)^{2(1+\eta)}
\right]
<\infty,
\]
and
\[
X_n(s)
\le
C(B)
\exp\!\left(
C_1\|\Delta_n\|_2^2
\right),
\qquad
s\in[0,1].
\]
Hence, by Cauchy--Schwarz,
\begin{align*}
\E_{\pi_n}
\left[
X_n(s)^{1+\eta}
\right]
&\le
\E_{\pi_n}
\left[
C(B)^{1+\eta}
e^{(1+\eta)C_1\|\Delta_n\|_2^2}
\right]
\\
&\le
\left(
\E_{\mathbb P_0}
\left[
C(B)^{2(1+\eta)}
\right]
\right)^{1/2}
\\
&\qquad\times
\left(
\E_{\pi_n}
\left[
e^{2(1+\eta)C_1\|\Delta_n\|_2^2}
\right]
\right)^{1/2}
\\
&\le
\left(
\E_{\mathbb P_0}
\left[
C(B)^{2(1+\eta)}
\right]
\right)^{1/2}
(1+\delta_n)^{1/2}.
\end{align*}
Therefore,
\begin{equation}\label{eq:X-uniform-moment}
\sup_{n\ge1}\sup_{s\in[0,1]}
\E_{\pi_n}
\left[
X_n(s)^{1+\eta}
\right]
<\infty.
\end{equation}

We now identify the pointwise limit. For fixed $s\in[0,1]$ and
$y\in\mathbb R^d$, the functions $b\longmapsto L_T(b,y)$ and $b\longmapsto\nabla_bL_T(b,y)$
are continuous. Moreover, the negative quadratic term in the
definition of $L_T$ implies the required boundedness for fixed $y$.
Since $\Delta_n\to0$ in $\pi_n$-probability, $L_T(B+s\Delta_n,y)-L_T(B,y)
\longrightarrow0$
and $\nabla_bL_T(B+s\Delta_n,y)
-
\nabla_bL_T(B,y)
\longrightarrow0$
in $\pi_n$-probability. It follows that
\[
\E_{\pi_n}
\left[
L_T(B+s\Delta_n,y)
\right]
\longrightarrow
\E_{\mathbb P_0}
\left[
L_T(B,y)
\right].
\]
Using the same Gaussian-integrability estimates as in
Lemma~\ref{lem:L22-fro} to pass through the $y$-integral, we obtain $G_s^{\pi_n}(B,\Delta_n)
\longrightarrow
H(B)$
in $\pi_n$-probability. Hence, $X_n(s)
\longrightarrow
\|H(B)\|_2^2$
in $\pi_n$-probability.

By \eqref{eq:X-uniform-moment} and Vitali's theorem,
\[
\E_{\pi_n}[X_n(s)]
\longrightarrow
\E_{\mathbb P_0}
\left[
\|H(B)\|_2^2
\right]
=
\|H\|_{L^2(\mathbb P_0)}^2
\]
for every $s\in[0,1]$. Moreover,
\eqref{eq:X-uniform-moment} provides an integrable bound independent of
$s$. Therefore, by dominated convergence,
\begin{equation}\label{eq:Gs-limit}
\int_0^1
\left(
\E_{\pi_n}[X_n(s)]
\right)^{1/2}
ds
\longrightarrow
\|H\|_{L^2(\mathbb P_0)}.
\end{equation}

Because the sequence $\delta_n\downarrow0$ and the couplings
$\pi_n\in\mathcal C_{\delta_n}$ were arbitrary, combining
\eqref{eq:Delta-L2-small}, \eqref{eq:J-upper-cs}, and
\eqref{eq:Gs-limit} yields
\begin{equation}\label{eq:uniform-two-sided-bound}
\sup_{\mu\in
\mathcal U_{\delta,B}^{\mathrm{OT}}(\mathbb P_0)}
\left|
\mathcal J(\mu)-\mathcal J(\mathbb P_0)
\right|
\le
\sqrt{\frac{\delta}{\tau}}\,
\|H\|_{L^2(\mathbb P_0)}
+
o(\sqrt{\delta}).
\end{equation}
In particular, when $\alpha<0$,
\begin{equation}\label{eq:sup-upper-bound}
\sup_{\mu\in
\mathcal U_{\delta,B}^{\mathrm{OT}}(\mathbb P_0)}
\mathcal J(\mu)
\le
\mathcal J(\mathbb P_0)
+
\sqrt{\frac{\delta}{\tau}}\,
\|H\|_{L^2(\mathbb P_0)}
+
o(\sqrt{\delta}).
\end{equation}

We next construct feasible deterministic pushforward perturbations
attaining the matching first-order bound. Write $h_H:=\|H\|_{L^2(\mathbb P_0)}>0.$
For $M>0$, define $H_M(b)
:=
H(b)\mathbf 1_{\{\|H(b)\|_2\le M\}},$ $h_M
:=
\|H_M\|_{L^2(\mathbb P_0)}.$
For all sufficiently large $M$, $h_M>0$. Set $u_M(b)
:=
\frac{H_M(b)}{h_M}.$
Then
\[
\E_{\mathbb P_0}
\left[
\|u_M(B)\|_2^2
\right]
=1,
\qquad
\|u_M\|_\infty<\infty.
\]

For fixed $M$, define
\[
\Gamma_M(t)
:=
\E_{\mathbb P_0}
\left[
\exp\!\left(
\tau t^2\|u_M(B)\|_2^2
\right)-1
\right],
\qquad
t\ge0.
\]
Since $u_M$ is bounded, dominated convergence gives
\[
\Gamma_M(t)
=
\tau t^2+o_M(t^2),
\qquad
t\downarrow0.
\]
Let $t_{\delta,M}>0$ satisfy $\Gamma_M(t_{\delta,M})=\delta.$
Then
\begin{equation}\label{eq:t-delta-M}
t_{\delta,M}
=
\sqrt{\frac{\delta}{\tau}}
\bigl(1+o_M(1)\bigr),
\qquad
\delta\downarrow0.
\end{equation}

For the maximizing problem $\alpha<0$, define $\mu_{\delta,M}^{+}
:=
(\operatorname{Id}+t_{\delta,M}u_M)_\#
\mathbb P_0.$
The deterministic coupling $\bigl(B,B+t_{\delta,M}u_M(B)\bigr)$
has exponential transport cost exactly $\delta$, so
\[
\mu_{\delta,M}^{+}
\in
\mathcal U_{\delta,B}^{\mathrm{OT}}(\mathbb P_0).
\]

For fixed $M$, the integral Taylor expansion and
Lemma~\ref{lem:L22-fro} yield
\begin{align*}
\mathcal J(\mu_{\delta,M}^{+})
-
\mathcal J(\mathbb P_0)
&=
t_{\delta,M}
\E_{\mathbb P_0}
\left[
H(B)^\top u_M(B)
\right]
+
o_M(t_{\delta,M}).
\end{align*}
Moreover,
\[
\E_{\mathbb P_0}
\left[
H(B)^\top u_M(B)
\right]
=
\frac{
\E_{\mathbb P_0}
\left[
\|H(B)\|_2^2
\mathbf 1_{\{\|H(B)\|_2\le M\}}
\right]
}{
h_M
}
=
h_M.
\]
Consequently, by \eqref{eq:t-delta-M},
\begin{equation}\label{eq:fixed-M-value}
\mathcal J(\mu_{\delta,M}^{+})
-
\mathcal J(\mathbb P_0)
=
\sqrt{\frac{\delta}{\tau}}\,
h_M
+
o_M(\sqrt{\delta}).
\end{equation}
Therefore, for every fixed sufficiently large $M$,
\[
\liminf_{\delta\downarrow0}
\frac{
\displaystyle
\sup_{\mu\in
\mathcal U_{\delta,B}^{\mathrm{OT}}(\mathbb P_0)}
\mathcal J(\mu)
-
\mathcal J(\mathbb P_0)
}{
\sqrt{\delta}
}
\ge
\frac{h_M}{\sqrt{\tau}}.
\]
Since
\[
h_M\uparrow h_H
\qquad
\text{as }M\to\infty,
\]
we conclude that
\[
\liminf_{\delta\downarrow0}
\frac{
\displaystyle
\sup_{\mu\in
\mathcal U_{\delta,B}^{\mathrm{OT}}(\mathbb P_0)}
\mathcal J(\mu)
-
\mathcal J(\mathbb P_0)
}{
\sqrt{\delta}
}
\ge
\frac{h_H}{\sqrt{\tau}}.
\]
Combining this with \eqref{eq:sup-upper-bound} gives
\[
\sup_{\mu\in
\mathcal U_{\delta,B}^{\mathrm{OT}}(\mathbb P_0)}
\mathcal J(\mu)
=
\mathcal J(\mathbb P_0)
+
\sqrt{\frac{\delta}{\tau}}\,
\|H\|_{L^2(\mathbb P_0)}
+
o(\sqrt{\delta}).
\]

It remains to obtain one feasible family with the stated
$L^2(\mathbb P_0)$ asymptotic direction. Since
\[
H_M\longrightarrow H
\quad\text{in }L^2(\mathbb P_0),
\qquad
h_M\longrightarrow h_H,
\]
we have
\[
u_M
\longrightarrow
\frac{H}{h_H}
\quad\text{in }L^2(\mathbb P_0).
\]
Choose $M_k\uparrow\infty$ and then choose
$\delta_k\downarrow0$ sufficiently fast so that, whenever
$0<\delta\le\delta_k$,
\[
\left|
\frac{t_{\delta,M_k}}
{\sqrt{\delta/\tau}}
-1
\right|
\le
\frac1k
\]
and the $o_{M_k}(\sqrt{\delta})$ remainder in
\eqref{eq:fixed-M-value} has absolute value at most
$\sqrt{\delta}/k$. Define
\[
M(\delta):=M_k
\qquad
\text{when }
\delta_{k+1}<\delta\le\delta_k,
\]
and set $\Delta_\delta^*(b)
:=
t_{\delta,M(\delta)}
u_{M(\delta)}(b).$
Then every $(\operatorname{Id}+\Delta_\delta^*)_\#
\mathbb P_0$
is feasible, asymptotically optimal, and
\[
\frac{\Delta_\delta^*}
{\sqrt{\delta/\tau}}
\longrightarrow
\frac{H}{\|H\|_{L^2(\mathbb P_0)}}
\quad
\text{in }L^2(\mathbb P_0).
\]
Equivalently,
\[
\Delta_\delta^*(b)
=
\sqrt{\frac{\delta}{\tau}}\,
\frac{H(b)}
{\|H\|_{L^2(\mathbb P_0)}}
+
o(\sqrt{\delta})
\quad
\text{in }L^2(\mathbb P_0).
\]

Finally, suppose $\alpha\in(0,1)$. The uniform two-sided estimate
\eqref{eq:uniform-two-sided-bound} gives
\[
\inf_{\mu\in
\mathcal U_{\delta,B}^{\mathrm{OT}}(\mathbb P_0)}
\mathcal J(\mu)
\ge
\mathcal J(\mathbb P_0)
-
\sqrt{\frac{\delta}{\tau}}\,
\|H\|_{L^2(\mathbb P_0)}
+
o(\sqrt{\delta}).
\]
For fixed $M$, use instead the feasible perturbation $\mu_{\delta,M}^{-}
:=
(\operatorname{Id}-t_{\delta,M}u_M)_\#
\mathbb P_0.$
The same Taylor expansion gives
\[
\mathcal J(\mu_{\delta,M}^{-})
-
\mathcal J(\mathbb P_0)
=
-
\sqrt{\frac{\delta}{\tau}}\,
h_M
+
o_M(\sqrt{\delta}).
\]
Letting $M\to\infty$ yields the matching upper bound for the infimum,
and therefore
\[
\inf_{\mu\in
\mathcal U_{\delta,B}^{\mathrm{OT}}(\mathbb P_0)}
\mathcal J(\mu)
=
\mathcal J(\mathbb P_0)
-
\sqrt{\frac{\delta}{\tau}}\,
\|H\|_{L^2(\mathbb P_0)}
+
o(\sqrt{\delta}).
\]
Applying the same diagonal construction with the negative sign gives a
feasible asymptotically optimal family satisfying
\[
\Delta_\delta^*(b)
=
-
\sqrt{\frac{\delta}{\tau}}\,
\frac{H(b)}
{\|H\|_{L^2(\mathbb P_0)}}
+
o(\sqrt{\delta})
\quad
\text{in }L^2(\mathbb P_0).
\]
This completes the proof.
\end{proof}

\section{Generalization of Nonlinear Projection Theorem with Non-compact Support}\label{PCLT}
{This section derives a stochastic asymptotic upper bound for the RWPI under
nonlinear projection when $B$ need not have compact support. The argument is
independent of Wasserstein geometry; technical lemmas appear at the end. Let
$\alpha<1$, $\alpha\ne0$, and
$c_\tau(\Delta):=e^{\tau\|\Delta\|_2^2}-1$. Let
$B^{(1)},\ldots,B^{(n)}$ be i.i.d.\ calibrated drift samples}
with common distribution $\mu^\star$, and let
$\mathbb P_{\mu^\star}$ denote the corresponding sampling law.
These observations are used to construct the empirical measure
$\mathbb P_n$, where $n$ is the sample size.
All stochastic-order and weak-convergence statements in this section
are understood under $\mathbb P_{\mu^\star}$. For a fixed $k > 0$, 
$$
g_k(x) = \left(\frac{e^{rT}}{k}\right)^{\frac{1}{1-\alpha}}x^{\frac{1}{1-\alpha}},
$$
and $k^*$ represents the optimal Lagrangian multiplier. We still make a sub-Gaussian assumption on $B$. Many computational details are similar to those in Section \ref{Th3} (only the sub-Gaussian parameters are different), and thus we omit some proofs for simplicity. 
\begin{assumption}\label{subGaussianCLT}
Let $\beta:=\frac{\alpha}{1-\alpha},$
and define
\[
M_{\mathrm{mid}}(\beta)
:=
\begin{cases}
\dfrac{2\beta}{(2-\beta)^2},
& 1<\beta<2,\\[2mm]
0,
& \text{otherwise}.
\end{cases}
\]
Suppose there exists $\gamma_0>0$ such that
\[
\E_{\mu^\star}
\left[
\exp\!\left(\gamma^2\|B\|_2^2\right)
\right]
<\infty
\qquad
\text{for every }\gamma<\gamma_0,
\]
and
\[
\frac{\gamma_0^2}{\|\sigma^{-1}\|_F^2}
>
T\max\left\{
4\beta^2-2\beta,\,
M_{\mathrm{mid}}(\beta),\,
16,\,
8\beta+8
\right\}.
\]
\end{assumption}
\begin{theorem}\label{CLT}

For each fixed $y\in\mathbb R^d$, define $F^\star(y)
:=
\E_{\mu^\star}
\left[
L_T(B,y)
\right].$
Also define $$c^\star
:=
\E_{\mu^\star}
\left[
\left\|
\int_{\mathbb R^d}
g_{k^\star}'\!\left(F^\star(y)\right)
\nabla_bL_T(B,y)
\varphi_T(y)\,dy
\right\|_2^2
\right],$$
and suppose that $c^\star>0.$
Then, under Assumption~\ref{subGaussianCLT}, as $n\to\infty$, we have the
asymptotic stochastic upper bound
\[
nR_n(k^\star)
\lesssim_D
L
:=
\frac{\tau Z^2}{c^\star},
\]
where $Z\sim\mathcal N\!\left(0,h(\mu^\star)\right)$,
and
\[
h(\mu^\star)
:=
\int_{\mathbb R^d}\int_{\mathbb R^d}
g_{k^\star}'\!\left(F^\star(y_1)\right)
g_{k^\star}'\!\left(F^\star(y_2)\right)
\operatorname{Cov}_{\mu^\star}\!\left(L_T(B,y_1),L_T(B,y_2)\right)
\varphi_T(y_1)\varphi_T(y_2)\,dy_1\,dy_2
<\infty.
\]
\end{theorem}
The proof of Theorem~\ref{CLT} consists of six parts. First, we use
Taylor expansions to separate the terms in $R_n(k^\star)$ and obtain
Eq.~\eqref{eq:Taylor_final}. Second, we estimate the remainder terms
$R_1$ and $R_2$ in Eq.~\eqref{eq:Taylor_final}. Third, we construct a
bounded deterministic correction direction and compute its linearized
form in Eq.~\eqref{eq:correction}. Fourth, we establish the convergence
results needed for the main statement. Fifth, we solve the exact
correction equation and evaluate the corresponding exponential
transport cost. Sixth, we combine the preceding estimates, let the
truncation level tend to infinity, and conclude the asymptotic
stochastic upper bound.

\subsection{Part I}
To begin with, the uniqueness of $k^*$ is easy to see. In order to prove Theorem~\ref{CLT}, for notational convenience, write $g:=g_{k^\star},$
and define, for a drift prior $\mu\in\mathcal P(\mathbb R^d)$,
\[
J(\mu)
:=
\int_{\mathbb R^d}
g\!\left(
\E_\mu[L_T(B,y)]
\right)
\varphi_T(y)\,dy.
\]
    
In particular, 
$$
g'(x) = \frac{1}{1-\alpha}\left(\frac{e^{rT}}{k^*}\right)^{\frac{1}{1-\alpha}}x^{\frac{\alpha}{1-\alpha}},
$$
and 
$$
g''(x) = \frac{\alpha}{(1-\alpha)^2}\left(\frac{e^{rT}}{k^*}\right)^{\frac{1}{1-\alpha}}x^{\frac{2\alpha - 1}{1-\alpha}}.
$$    
We notice that $L_T$ and $g$ are both twice continuously differentiable functions in each argument. Therefore, the RWP function becomes 
    \begin{align*}
R_n(k^\star)
&:=
\inf_{\mu \in \mathcal{M}_{k^\star}}
D_c(\mathbb{P}_n,\mu)
\\
&=
\inf\left\{
D_c(\mu,\mathbb{P}_n):
\int_{\mathbb{R}^d}
I\left(
\frac{k^\star e^{-rT}}
{\mathbb{E}_{\mu}\left[L_T(B,y)\right]}
\right)
\varphi_T(y)\,dy
=
x_0e^{rT}
\right\}
\\
&=
\inf\left\{
D_c(\mu,\mathbb{P}_n):
J(\mu)=J(\mu^\star)
\right\}.
\end{align*}
{Given the empirical measure $\mathbb P_n$, define $\mathbb P_n^\Delta$ by
moving each atom $B^{(i)}$ to $B^{(i)}+n^{-1/2}\Delta_i$, where
$\Delta_i\in\mathbb R^d$. We show that, for every fixed $\varepsilon>0$ and
all sufficiently large $n$, with probability at least $1-\varepsilon$ there
exists a correction $\Delta$ such that $\|\Delta\|_n=O_p(1)$ and
$J(\mathbb P_n^\Delta)=J(\mu^\star)$. For $\gamma\in\mathbb Z$, write
$\|\Delta\|_n^\gamma:=n^{-1}\sum_{i=1}^n\|\Delta_i\|^\gamma$.}

For each $y$ and $i$, set $h_i:=\Delta_i/\sqrt n$.
The Taylor theorem in the $b$--variable gives
\[
L_T(B^{(i)}+h_i,y)
= L_T(B^{(i)},y)+\nabla_b L_T(B^{(i)},y)\!\cdot\!h_i
+\int_0^1(1-t)\,h_i^\top \nabla_b^2 L_T(B^{(i)}+t h_i,y)\,h_i\,dt.
\]
Defining $m_\mu(y):=\E_\mu[L_T(B,y)].$ for a probability measure $\mathbb{P}$ yields
\[
m_{\mathbb{P}_n^\Delta}(y)
= m_{\mathbb{P}_n}(y)\;+\;A_1(y)\;+\;A_2(y),
\]
with
\[
A_1(y)=\frac{1}{\sqrt n}\cdot\frac1n\sum_{i=1}^n \nabla_b L_T(B^{(i)},y)\cdot \Delta_i
\]
and
\[
A_2(y)=\frac1n\sum_{i=1}^n\int_0^1(1-t)\,\frac{\Delta_i^\top}{\sqrt n}\,
\nabla_b^2 L_T\!\Big(B^{(i)}+t\tfrac{\Delta_i}{\sqrt n},y\Big)\,\frac{\Delta_i}{\sqrt n}\,dt.
\]
Another Taylor expansion gives 
\begin{align*}
    g(m_{\mathbb{P}_n^\Delta}(y))& = g(m_{\mathbb{P}_n}(y)) + g'(m_{\mathbb{P}_n}(y))(A_1(y) + A_2(y))\\ 
    &+ \int_0^1 (1-t)g''(m_{\mathbb{P}_n}(y) + t(A_1(y) + A_2(y)))(A_1(y) + A_2(y))^2dt ,
\end{align*}which implies that 
\begin{equation}\label{eq:Taylor_final}
J(\mathbb P_n^\Delta)-J(\mu^\star)
=
J(\mathbb P_n)-J(\mu^\star)
+
\int_{\mathbb R^d}
g'\!\left(m_{\mathbb P_n}(y)\right)
A_1(y)\varphi_T(y)\,dy
+
R_1+R_2,
\end{equation}
where $$R_1 = \int g'(m_{\mathbb{P}_n}(y))A_2(y)\varphi_T(y)dy,$$
and $$R_2 = \int_0^1 (1-t)\int g''(m_{\mathbb{P}_n}(y) + t(A_1(y) + A_2(y)))(A_1(y) + A_2(y))^2\varphi_T(y)dydt.$$
If we define $$C_i(n) = \int g'(m_{\mathbb{P}_n}(y))\nabla_bL_T(B^{(i)},y)\varphi_T(y)dy,$$
then $$\int g'(m_{\mathbb{P}_n}(y))A_1(y)\varphi_T(y)dy = \frac{1}{n^{1/2}}\frac{1}{n}\sum_{i = 1}^n C_i(n) \cdot \Delta_i.$$

\subsection{Part II}
In this part, Lemmas~\ref{negalphaR1} and~\ref{boundR2} provide
general bounds for the remainder terms $R_1$ and $R_2$. We then
specialize these bounds to bounded truncations of the population
gradient, which will be used in Parts III--V. The sub-Gaussian
assumptions are used to control the required Gaussian integrals. We
give the proof details of the most involved estimate,
Lemma~\ref{boundR2}.

\subsubsection{Part II-i}
We propose the following sub-Gaussian assumption. 
\begin{assumption}\label{Lemma12ass}
Suppose there exists $\gamma_0>0$ such that
\[
\E_{\mu^\star}
\left[
\exp\!\left(\gamma^2\|B\|_2^2\right)
\right]
<\infty
\qquad
\text{for every }\gamma<\gamma_0,
\]
and
\[
\frac{\gamma_0^2}{\|\sigma^{-1}\|_F^2}
>
T\max\left\{
2(2\beta^2-\beta),\,4
\right\},
\qquad
\beta=\frac{\alpha}{1-\alpha}.
\]
\end{assumption}

\begin{lemma}\label{negalphaR1}
In the context of the proof of Theorem \ref{CLT}, under Assumption \ref{Lemma12ass}, when $\alpha < 1$ and $\alpha \neq 0$, there exists a constant $C > 0$ such that
\[
|R_1|
= O_p\!\left(
\frac{1}{n}\Big(\frac{1}{n}\sum_{i=1}^n \|\Delta_i\|^4\Big)^{\!1/2}\,
\Big[1+\frac{\|\Delta\|_n}{\sqrt{n}}+\frac{1}{n}\Big(\frac{1}{n}\sum_{i=1}^n \|\Delta_i\|^4\Big)^{\!1/2}\Big]\,
\exp\!\Big(C'\,\max_{1\le i\le n}\frac{\|\Delta_i\|^2}{n}\Big)
\right).
\]
\end{lemma}

\begin{assumption}\label{Lemma13ass}
Suppose there exists $\gamma_0>0$ such that
\[
\E_{\mu^\star}
\left[
\exp\!\left(\gamma^2\|B\|_2^2\right)
\right]
<\infty
\qquad
\text{for every }\gamma<\gamma_0,
\]
and
\[
\frac{\gamma_0^2}{\|\sigma^{-1}\|_F^2}
>
T\max\left\{
M_{\mathrm{mid}}(\beta),\,
6,\,
2p^2-p
\right\},
\qquad
p:=\beta-1,
\qquad
\beta:=\frac{\alpha}{1-\alpha}.
\]
\end{assumption}

\begin{lemma}\label{boundR2}
In the context of proof of Theorem \ref{CLT}, under Assumption \ref{Lemma13ass}, when $\alpha < 1$ and $\alpha \neq 0$, there exists a constant $C > 0$ such that
$$|R_2| = O_p\left(\frac{\left\Vert \Delta\right\Vert_n^2}{n}\exp\left(C\max_{1\ \leq i \leq n}\frac{\left\Vert \Delta_i\right\Vert^2}{n}\right)\right) + O_p\!\left(
\frac{\|\Delta\|_n^2}{n}\cdot \frac{\max_{1\le i\le n}\|\Delta_i\|^2}{n}\;
\exp\!\Big(C\,\max_{1\le i\le n}\frac{\|\Delta_i\|^2}{n}\Big)
\right).$$
\end{lemma}

\begin{proof}
Bounding $R_2$ is equivalent to bounding these two terms:
$$|R_2| \leq I_1 + I_2,$$
where 
$$I_1 = 2\int_0^1 (1-t)\int \left|g''(m_{\mathbb{P}_n}(y) + t(A_1(y) + A_2(y)))\right|A_1(y)^2\varphi_T(y)dydt$$and 
$$I_2 = 2\int_0^1 (1-t)\int\left| g''(m_{\mathbb{P}_n}(y) + t(A_1(y) + A_2(y)))\right|A_2(y)^2\varphi_T(y)dydt.$$  

Recall that $\beta = \frac{\alpha}{1-\alpha}$ and $g''(x)=C_\alpha\,x^{\beta-1}$ for a constant $C_{\alpha}$. The {Cauchy--Schwarz inequality} yields
\begin{align*}
|A_1(y)|
&\le \frac{1}{\sqrt n}\cdot\frac{1}{n}\Big(\sum_{i=1}^n \|\nabla_b L_T(B^{(i)},y)\|_2^2\Big)^{\!1/2}\!
\Big(\sum_{i=1}^n \|\Delta_i\|_2^2\Big)^{\!1/2}
= \sqrt{\frac{\|\Delta\|_n^{\,2}}{n}}\;\Big(\overline{\Xi^2}(y)\Big)^{\!1/2},
\end{align*}
with
\[
\overline{\Xi^2}(y):=\frac{1}{n}\sum_{i=1}^n \|\nabla_b L_T(B^{(i)},y)\|_2^2
\ \le\ C_1\Big(\|y\|_2^2\,\overline{L^2}(y)\;+\;T^2\,\overline{A^2L^2}(y)\Big),
\]
where
\[
\overline{L^2}(y):=\frac1n\sum_{i=1}^n L_T(B^{(i)},y)^2,\qquad
\overline{A^2L^2}(y):=\frac1n\sum_{i=1}^n \|a(B^{(i)})\|_2^2\,L_T(B^{(i)},y)^2,
\]
and therefore, 
$$|A_1(y)|^2 \;\le\; \frac{\|\Delta\|_n^{\,2}}{n}\;C_1\Big(\|y\|_2^2\,\overline{L^2}(y)\;+\;T^2\,\overline{A^2L^2}(y)\Big).$$

Recall that $h_i:=\Delta_i/\sqrt n$. Define the shifted mixture
\[
m_{\mathbb{P}_n}^{(h)}(y):=\frac1n\sum_{i=1}^n L_T(B^{(i)}+h_i,y).
\]
Therefore, 
\[
x_t(y):=(1-t)\,m_{\mathbb{P}_n}(y)+t\,m_{\mathbb{P}_n}^{(h)}(y)
=m_{\mathbb{P}_n}(y)+t\big(A_1(y)+A_2(y)\big).
\]

Let $p:=\beta-1 = \frac{\alpha}{1-\alpha} - 1$. 
We write $m:=m_{\mathbb{P}_n}(y)>0$ and $m^{(h)}:=m_{\mathbb{P}_n}^{(h)}(y)>0$ to ease notation, so
$x_t=(1-t)m+t m^{(h)}$.

When $p < 0$ or $p \geq 1$,
the function $f(x)=x^{p}$ on $(0,\infty)$ is convex. Therefore
\[
x_t^{p}=f\big((1-t)m+t m^{(h)}\big)
\;\le\; (1-t)f(m)+t f(m^{(h)})=(1-t)m^{p}+t (m^{(h)})^{p}
\;\le\; m^{p}+(m^{(h)})^{p}.
\]

When $p \in (0,1)$, the map $f(x)=x^{p}$ is increasing and concave. Since $x_t\le m+m^{(h)}$ and $f$ is increasing, we have
\[
x_t^{p}\ \le\ (m+m^{(h)})^{p}\ \le\ m^{p}+(m^{(h)})^{p},
\]
where the last step uses the subadditivity $(a+b)^{p}\le a^{p}+b^{p}$ for $a,b\ge0$ and $0<p\le1$.

Therefore, $$|g''(x_t)|=|C_\alpha|\,x_t^p \le |C_\alpha|\big(m^{\beta-1}+(m^{(h)})^{\beta-1}\big).$$
Hence, 
\begin{align*}
I_1
&=2\!\int_0^1(1-t)\!\int |g''(x_t)|\,A_1^2(y)\,\varphi_T(y)dy
\ \le\ \frac{C}{n}\,\|\Delta\|_n^{\,2}\,\Big(J_1+J_2+J_1^{(h)}+J_2^{(h)}\Big),
\end{align*}
where
\[
J_1:=\int \|y\|^2 m_{\mathbb{P}_n}^{\beta-1}\,\overline{L^2}\,\varphi_T,\qquad
J_2:=\int m_{\mathbb{P}_n}^{\beta-1}\,\overline{A^2L^2}\,\varphi_T,
\]
\[
J_1^{(h)}:=\int \|y\|^2 (m_{\mathbb{P}_n}^{(h)})^{\beta-1}\,\overline{L^2}\,\varphi_T,\qquad
J_2^{(h)}:=\int (m_{\mathbb{P}_n}^{(h)})^{\beta-1}\,\overline{A^2L^2}\,\varphi_T.
\]

We first focus on the case when $p:= \beta - 1 < 0$, by the Jensen inequality,
\[
m_{\mathbb{P}_n}(y) \;=\; \E_{\mathbb{P}_n}\big[e^{\langle a(B),y\rangle-\frac{T}{2}\|a(B)\|_2^2}\big]
\;\ge\;
\exp\!\Big( \langle \E_{\mathbb{P}_n}[a(B)],y\rangle - \tfrac{T}{2}\,\E_{\mathbb{P}_n}\left[\|a(B)\|_2^2\right] \Big).
\]
Raising to the negative power $p<0$ reverses the inequality, giving
$$
m_{\mathbb{P}_n}(y)^{p}\;\le\;C_a\,\exp\!\big(p\,\langle v,y\rangle\big),
\qquad
C_a:=\exp\!\Big(-\tfrac{p T}{2}\,\E_{\mathbb{P}_n}\left[\|a(B)\|_2^2\right]\Big),\ \ 
v:=\E_{\mathbb{P}_n}[a(B)].
$$
We denote $c_i:=\sigma^{-1}h_i$. Then we have 
\[
m_{\mathbb{P}_n}^{(h)}(y)=\frac1n\sum_{i=1}^n L_T(B^{(i)}+h_i,y)
=\frac1n\sum_{i=1}^n \exp\!\Big(\,\big\langle a(B^{(i)})+c_i,\,y\big\rangle-\tfrac{T}{2}\|a(B^{(i)})+c_i\|_2^2\Big).
\]
Introduce the notations
\[
v_h:=\frac1n\sum_{i=1}^n \big(a(B^{(i)})+c_i\big), 
\qquad
s_{2,h}:=\frac1n\sum_{i=1}^n \|a(B^{(i)})+c_i\|_2^2,
\qquad
H_n:=\max_{1\le i\le n}\|c_i\|_2.
\]
By {Jensen's inequality},
\[
m_{\mathbb{P}_n}^{(h)}(y)\ \ge\ \exp\!\Big(\,\langle v_h,y\rangle-\tfrac{T}{2}\,s_{2,h}\Big).
\]
Since $p<0$, raising both sides to the power $p$ reverses the inequality:
\[
\big(m_{\mathbb{P}_n}^{(h)}(y)\big)^{p}
\ \le\ \exp\!\Big(\,p\,\langle v_h,y\rangle-\tfrac{pT}{2}\,s_{2,h}\Big).
\]
Set
\[
s_2:=\E_{\mathbb{P}_n}\!\big[\|a(B)\|_2^2\big],\qquad
C_a:=\exp\!\Big(-\tfrac{pT}{2}\,s_2\Big).
\]
Hence, 
\[
J_1
\ \le\ \frac{C_a}{n}\sum_{i=1}^n
\int \|y\|_2^2\,e^{\,p\langle v,y\rangle}\,L_T(B^{(i)},y)^2\,\varphi_T(y)\,dy.
\]
Write $a_i:=a(B^{(i)})$ and note $L_T(B^{(i)},y)^2=\exp(2\langle a_i,y\rangle - T\|a_i\|_2^2)$.
Let $Y\sim N(0,TI_d)$ so that $\varphi_T$ is its density. Then, with
\[
\lambda_i:=2a_i+p\,v\in\R^d,
\]
we have 
\[
\int \|y\|_2^2\,e^{\,p\langle v,y\rangle}\,L_T(B^{(i)},y)^2\,\varphi_T(y)\,dy
= e^{-T\|a_i\|_2^2}\,(Td+T^2\|\lambda_i\|_2^2)\,e^{\frac{T}{2}\|\lambda_i\|_2^2}.
\]
Therefore,
$$ 
J_1
\ \le\ \frac{C_a}{n}\sum_{i=1}^n
\Big(Td+T^2\|\lambda_i\|_2^2\Big)\,
\exp\!\Big(\tfrac{T}{2}\|\lambda_i\|_2^2 - T\|a_i\|_2^2\Big).$$
By the inequality $(\|u+z\|_2^2\le 2\|u\|_2^2+2\|z\|_2^2)$ with $u=2a_i$, $z=pv$:
\[
\|\lambda_i\|_2^2=\|2a_i+p v\|_2^2 \ \le\ 8\|a_i\|_2^2 + 2p^2\|v\|_2^2.
\]
Hence, 
\[
\frac{T}{2}\|\lambda_i\|_2^2 - T\|a_i\|_2^2
\ \le\ 3T\|a_i\|_2^2 + T p^2\|v\|_2^2.
\]
Also, $\|\lambda_i\|_2^2\le 8\|a_i\|_2^2+2p^2\|v\|_2^2$ implies
\(
Td+T^2\|\lambda_i\|_2^2
\le C\big(1+\|a_i\|_2^2+\|v\|_2^2\big)
\)
for a constant $C=C(T,p,d)$.
Therefore,
\begin{align*}
  J_1
\ &\le\ C\,C_a\,e^{Tp^2\|v\|_2^2}\,
\frac{1}{n}\sum_{i=1}^n \big(1+\|a_i\|_2^2+\|v\|_2^2\big)\,e^{3T\|a_i\|_2^2}.
\end{align*}
Set
\[
A_n:=\frac1n\sum_{i=1}^n e^{3T\|a_i\|_2^2},\qquad
B_n:=\frac1n\sum_{i=1}^n \|a_i\|_2^2\,e^{3T\|a_i\|_2^2}.
\]
Then
\[
J_1\ \le\ C\,C_a\,e^{Tp^2\|v\|_2^2}\,\Big[(1+\|v\|_2^2)\,A_n + B_n\Big].
\]
Since $a(B)=\sigma^{-1}B-m$, we have
\(
\|a(B)\|_2^2 \le 2\|\sigma^{-1}B\|_2^2 + 2\|m\|_2^2
\),
and thus
\[
\E\big[e^{3T\|a(B)\|_2^2}\big]
\ \le\ e^{6T\|m\|_2^2}\,\E\big[e^{6T\|\sigma^{-1}B\|_2^2}\big]
\ \le\ e^{6T\|m\|_2^2}\,\E\big[e^{6T\|\sigma^{-1}\|_F^2\|B\|_2^2}\big].
\]
Hence, $\E[e^{3T\|a(B)\|_2^2}]<\infty$ if the sub-Gaussian parameter satisfies $\frac{\gamma_0^2}{\|\sigma^{-1}\|_F^2}\;>\;6T.$ The same argument with a polynomial prefactor yields
\(
\E[\|a(B)\|_2^2 e^{3T\|a(B)\|_2^2}]<\infty
\).
By the law of large numbers, $A_n=O_p(1)$ and $B_n=O_p(1)$. Also, by Jensen's inequality,
\(
\|v\|_2^2=\big\|\tfrac1n\sum a_i\big\|_2^2 \le \tfrac1n\sum \|a_i\|_2^2
\),
so $\|v\|_2^2=O_p(1)$ and therefore $C_a e^{Tp^2\|v\|_2^2}=O_p(1)$.
Thus, $J_1=O_p(1)$. The bound of $J_2 = O_p(1)$ is almost verbatim since there is no additional term on the exponential.

Similarly, it suffices to upper bound $J_1^{(h)}$. 
We introduce new analogous notation
\[
C_a^{(h)}:=\exp\!\Big(-\tfrac{pT}{2}\,s_{2,h}\Big).
\]
Therefore, 
\[
J_1^{(h)}
\ \le\ \frac{C_a^{(h)}}{n}\sum_{i=1}^n
\int \|y\|_2^2\,e^{\,p\langle v_h,y\rangle}\,L_T(B^{(i)},y)^2\,\varphi_T(y)\,dy.
\]
An analogous computation of the case of $J_1$ shows that
\[
J_1^{(h)}
\ \le\ C\,C_a^{(h)}\,e^{Tp^2\|v_h\|_2^2}\,
\frac{1}{n}\sum_{i=1}^n \big(1+\|a_i\|_2^2+\|v_h\|_2^2\big)\,e^{3T\|a_i\|_2^2}.
\]
Set
\[
A_n:=\frac1n\sum_{i=1}^n e^{3T\|a_i\|_2^2},\qquad
B_n:=\frac1n\sum_{i=1}^n \|a_i\|_2^2\,e^{3T\|a_i\|_2^2}.
\]
Then
\[
J_1^{(h)}
\ \le\ C\,C_a^{(h)}\,e^{Tp^2\|v_h\|_2^2}\,\Big[(1+\|v_h\|_2^2)\,A_n + B_n\Big].
\]
Note
\(
\|v_h\|_2\le \|v\|_2+\bar c,\ \bar c:=\tfrac1n\sum_i\|c_i\|_2\le H_n,
\)
so \(\|v_h\|_2^2\le 2\|v\|_2^2+2H_n^2\). Also, 
\(
s_{2,h}\le 2\,s_2+2\,H_n^2,\ \ s_2:=\tfrac1n\sum_i\|a_i\|_2^2.
\)
Hence, 
\[
C_a^{(h)}\,e^{Tp^2\|v_h\|_2^2}\ \le\ e^{C\,H_n^2}\,\underbrace{\exp\!\Big(-\tfrac{pT}{2}\,s_2\Big)}_{=:C_a}\,\exp\!\big(Tp^2\|v\|_2^2\big).
\]
Under the same sub-Gaussian condition as for $ J_1$,
\[
A_n=O_p(1),\quad B_n=O_p(1),\quad \|v\|_2=O_p(1),
\]
so for a constant $C>0$,
\[
\ J_1^{(h)}\ =\ O_p\!\big(e^{C H_n^2}\big) = O_p\left(\exp\left(C\max_{1\leq i \leq n}\frac{\left\Vert \Delta\right\Vert^2}{n}\right)\right). 
\]
Therefore, when $p < 0$,
$$I_1 = O_p\left(\frac{\left\Vert \Delta\right\Vert^2_n}{n}\exp\left(C\max_{1\leq i \leq n}\frac{\left\Vert \Delta\right\Vert^2}{n}\right)\right).$$

Next, we bound $I_1$ in the case when $p \geq 0$. It suffices to consider $J_1$ and $J_1^{(h)}$ here again.
First we assume $p \in (0,1)$, then
{by H\"older's inequality with exponents $r=\frac1p$ and $s=\frac1{1-p}$,}
\[
J_1=\int \|y\|^2\,m_{\mathbb{P}_n}(y)^p\,\overline{L^2}(y)\,\varphi_T(y)\,dy
\ \le\ \Big(\int m_{\mathbb{P}_n}(y)\,\varphi_T(y)\,dy\Big)^{\!p}
\Big(\int \|y\|^{2s}\,\overline{L^2}(y)^{\,s}\,\varphi_T(y)\,dy\Big)^{\!1/s}.
\]
Since $\int L_T(b,y)\varphi_T(y)dy=1$ for all $b$, Fubini's theorem gives $\int m_{\mathbb{P}_n}\varphi_T=1$,
so the first factor is $1$. For the second term, by convexity of $x\mapsto x^{s}$ ($s>1$),
\begin{align*}
J_1\ &\le\ \Big(\tfrac1n\sum_{i=1}^n \int \|y\|^{2s} L_T(B^{(i)},y)^{2s}\varphi_T(y)dy\Big)^{\!1/s}\\
&\le\ C\,\Big(\tfrac1n\sum_{i=1}^n \big(1+\|a_i\|^{2s}\big)\,
e^{\,T s(2s-1)\,\|a_i\|^2}\Big)^{\!1/s}.
\end{align*}
Therefore, a sufficient sub-Gaussian condition ensuring $J_1=O_p(1)$ is
\[
 \frac{\gamma_0^2}{\|\sigma^{-1}\|_F^2} > 2\,T\,s(2s-1).
\]
Similarly,
\[
J_1^{(h)}=\int \|y\|^2\,(m_{\mathbb{P}_n}^{(h)}(y))^{p}\,\overline{L^2}(y)\,\varphi_T(y)\,dy
\ \le\ \Big(\int m_{\mathbb{P}_n}^{(h)}(y)\,\varphi_T(y)\,dy\Big)^{\!p}
\Big(\int \|y\|^{2s}\,\overline{L^2}(y)^{\,s}\,\varphi_T(y)\,dy\Big)^{\!1/s}.
\]
Since $\int L_T(b,y)\varphi_T(y)dy=1$ for every $b$, we have
\[
\int m_{\mathbb{P}_n}^{(h)}(y)\,\varphi_T(y)\,dy
=\frac1n\sum_{j=1}^n \int L_T(B^{(j)}+h_j,y)\varphi_T(y)\,dy
=1,
\]
so with the same sub-Gaussian assumption $J_1^{(h)}=O_p(1)$, when $p \in (0,1)$,
$$I_1 = O_p\left(\frac{\left\Vert \Delta\right\Vert^2_n}{n}\right).$$

When $p \geq 1$, Jensen's inequality gives \(m_{\mathbb{P}_n}(y)^p \le \E_{\mathbb{P}_n}[L_T(B,y)^{p}]\). Thus
\[
J_1\ \le\ \frac1n\sum_{i=1}^n \E_{\mathbb{P}_n}\!\left[\int \|y\|^2\,L_T(B,y)^{p}\,L_T(B^{(i)},y)^{2}\,\varphi_T(y)\,dy\right].
\]
For fixed $(B,B^{(i)})$, write $\lambda:=p\,a(B)+2\,a_i$.
Completing the square (as in the $J_1$ computation) yields
\[
\int \|y\|^2\,L_T(B,y)^{p}L_T(B^{(i)},y)^{2}\varphi_T(y)dy
= (Td+T^2\|\lambda\|^2)\,\exp\!\Big\{\tfrac{T}{2}\|\lambda\|^2-\tfrac{T}{2}p\|a(B)\|^2-T\|a_i\|^2\Big\}.
\]
Using $\|\lambda\|^2\le 2p^2\|a(B)\|^2+8\|a_i\|^2$, we have 
\[
\frac{T}{2}\|\lambda\|^2-\tfrac{T}{2}p\|a(B)\|^2-T\|a_i\|^2
\ \le\ \tfrac{T}{2}\,(2p^2-p)\,\|a(B)\|^2\ +\ 3T\,\|a_i\|^2.
\]
Taking expectation in $B$, we have
\[
J_1\ \le\ C\,\frac1n\sum_{i=1}^n \big(1+\|a_i\|^2\big)\,e^{3T\|a_i\|^2}\;
\E_{\mathbb{P}_n}\!\left[e^{\frac{T}{2}(2p^2-p)\,\|a(B)\|^2}\right].
\]
Hence, a sufficient sub-Gaussian condition (thus $J_1=O_p(1)$) is
$$\frac{\gamma_0^2}{\|\sigma^{-1}\|_F^2} > T\max\left\{\,(2p^2-p),6\right\}.$$
A similar computation gives
\[
\quad
J_1^{(h)}\ \le\ C\,e^{C p\,H_n^2}\,
\Big(\frac1n\sum_{i=1}^n (1+\|a_i\|^2)\,e^{3T\|a_i\|^2}\Big)\,
\Big(\frac1n\sum_{j=1}^n e^{\,[\,T(2p^2-p)+\varepsilon\,]\|a_j\|^2}\Big),
\quad
\]where the bound is finite under 
the same sub-Gaussian assumption
\[
\ \frac{\gamma_0^2}{\|\sigma^{-1}\|_F^2}\ >\ \max\!\big\{6T,\;T(2p^2-p)\big\}.
\]
Therefore, when $p \geq 1$,
$$I_1 = O_p\left(\frac{\left\Vert \Delta\right\Vert^2_n}{n}\exp\left(C\max_{1\leq i \leq n}\frac{\left\Vert \Delta\right\Vert^2}{n}\right)\right).$$

Next, we focus on the bound of $I_2$ and begin with a bound for $A_2(y)^2$.
To begin with, direct computations show that there exists a constant $C_F=C_F(d,T,\|\sigma^{-1}\|_F) > 0$ such that
$$\big\|\nabla_b^2 L_T(b,y)\big\|_F
\;\le\; C_F\,L_T(b,y)\,\Big(1+\|y\|_2^2+T^2\|a(b)\|_2^2\Big). $$

{
Recall that $c_i:=\sigma^{-1}h_i$ and choose
$\eta_y>0$ and $\eta_a>0$ sufficiently large as specified below.
Then
}
\[
L_T(B^{(i)}+t h_i,y)
= L_T(B^{(i)},y)\,\exp\!\Big\{t\langle c_i,y\rangle - Tt\langle a_i,c_i\rangle - \tfrac{T}{2}t^2\|c_i\|_2^2\Big\}.
\]
By the Young inequalities,
\[
t\langle c_i,y\rangle \le \frac{\|y\|_2^2}{4\eta_y} + \eta_y t^2\|c_i\|_2^2,
\qquad
-\,Tt\langle a_i,c_i\rangle
\le \frac{T^2 t^2}{4\eta_a}\|a_i\|_2^2+\eta_a\|c_i\|_2^2.
\]
{
Since $t\in[0,1]$, we have $t^2\le1$, and the negative
$-\tfrac{T}{2}t^2\|c_i\|_2^2$ can be dropped. Thus, for a finite
constant $C=C(T,\eta_y,\eta_a)$, we have
}
$$L_T(B^{(i)}+t h_i,y)
\ \le\ \exp\!\Big\{\tfrac{\|y\|_2^2}{4\eta_y}\Big\}\,
\exp\!\Big\{C\,\|c_i\|_2^2\Big\}\,
\exp\!\Big\{\tfrac{T^2}{4\eta_a}\,\|a_i\|_2^2\Big\}\,
L_T(B^{(i)},y).$$
Also, $\|a(B^{(i)}+t h_i)\|_2 \le \|a_i\|_2+\|c_i\|_2$. Hence,
$$
\|a(B^{(i)}+t h_i)\|_2^2 \le 2\|a_i\|_2^2+2\|c_i\|_2^2.
$$
Therefore,
\[
|A_2(y)|
=\frac1n\sum_{i=1}^n\int_0^1 (1-t)\,\big|h_i^\top\nabla_b^2 L_T(B^{(i)}+t h_i,y)\,h_i\big|\,dt
\le \frac1n\sum_{i=1}^n \|h_i\|_2^2\int_0^1 \big\|\nabla_b^2 L_T(B^{(i)}+t h_i,y)\big\|_F\,dt.
\]
{
Set
\[
\kappa:=\frac{T^2}{4\eta_a},\qquad
S_n:=\frac1n\sum_{i=1}^n\|c_i\|_2^2,\qquad
H_n:=\max_{1\le i\le n}\|c_i\|_2,
\]
and define the weighted empirical mixtures
\[
\overline{L^2}_{\kappa}(y)
:=
\frac1n\sum_{i=1}^n
e^{2\kappa\|a_i\|_2^2}L_T(B^{(i)},y)^2,
\]
and
\[
\overline{A^4L^2}_{\kappa}(y)
:=
\frac1n\sum_{i=1}^n
\|a_i\|_2^4
e^{2\kappa\|a_i\|_2^2}L_T(B^{(i)},y)^2.
\]

Applying Cauchy--Schwarz to the empirical sum and using $\frac1n\sum_{i=1}^n\|c_i\|_2^4
\le H_n^2S_n,$
we obtain
\[
|A_2(y)|
\le
C\|\sigma\|_F^2(S_nH_n^2)^{1/2}
\exp\!\Big\{\frac{\|y\|_2^2}{4\eta_y}\Big\}
e^{CH_n^2}
\Big[
(1+\|y\|_2^2+T^2H_n^2)
\overline{L^2}_{\kappa}(y)^{1/2}
+
T^2\overline{A^4L^2}_{\kappa}(y)^{1/2}
\Big].
\]
Consequently,
\[
A_2(y)^2
\le
C\|\sigma\|_F^4S_nH_n^2
\exp\!\Big\{\frac{\|y\|_2^2}{2\eta_y}\Big\}
e^{CH_n^2}
\Big[
(1+\|y\|_2^2+T^2H_n^2)^2
\overline{L^2}_{\kappa}(y)
+
T^4\overline{A^4L^2}_{\kappa}(y)
\Big].
\]
}

Similarly as in the bound for $I_1$, we have
\begin{align*}
I_2
&=
2\!\int_0^1(1-t)
\int |g''(x_t)|A_2(y)^2\varphi_T(y)\,dy\,dt\\
&\le
C\|\sigma\|_F^4S_nH_n^2e^{CH_n^2}
\Big(
\widetilde J_1+\widetilde J_2
+\widetilde J_1^{(h)}+\widetilde J_2^{(h)}
\Big),
\end{align*}
where
\begin{align*}
\widetilde J_1
&:=
\int
(1+\|y\|_2^2+T^2H_n^2)^2
e^{\|y\|_2^2/(2\eta_y)}
m_{\mathbb P_n}(y)^{\beta-1}
\overline{L^2}_{\kappa}(y)
\varphi_T(y)\,dy,\\
\widetilde J_2
&:=
T^4\int
e^{\|y\|_2^2/(2\eta_y)}
m_{\mathbb P_n}(y)^{\beta-1}
\overline{A^4L^2}_{\kappa}(y)
\varphi_T(y)\,dy,\\
\widetilde J_1^{(h)}
&:=
\int
(1+\|y\|_2^2+T^2H_n^2)^2
e^{\|y\|_2^2/(2\eta_y)}
\big(m_{\mathbb P_n}^{(h)}(y)\big)^{\beta-1}
\overline{L^2}_{\kappa}(y)
\varphi_T(y)\,dy,\\
\widetilde J_2^{(h)}
&:=
T^4\int
e^{\|y\|_2^2/(2\eta_y)}
\big(m_{\mathbb P_n}^{(h)}(y)\big)^{\beta-1}
\overline{A^4L^2}_{\kappa}(y)
\varphi_T(y)\,dy.
\end{align*}

{
The four $\widetilde J$-terms are controlled by the same
completion-of-squares estimates as their $I_1$ analogues.
The only changes are the additional weight
$e^{2\kappa\|a_i\|_2^2}$ and the Gaussian factor
$e^{\|y\|_2^2/(2\eta_y)}$.

Since all inequalities in Assumption~\ref{Lemma13ass} are strict,
we may first choose $\eta_a$ sufficiently large so that
$\kappa=T^2/(4\eta_a)$ is small enough to fit within the same
sub-Gaussian moment slack, and then choose $\eta_y$ sufficiently
large so that the corresponding Gaussian completion-of-squares
coefficients also remain within that slack.

When $0<p<1$, with $s=(1-p)^{-1}$, choose in addition
$\eta_y>sT$. After H\"older's inequality, $e^{s\|y\|_2^2/(2\eta_y)}\varphi_T(y)$
is proportional to a centered Gaussian density with variance
\[
T_{\eta_y,s}
:=
\frac{T\eta_y}{\eta_y-sT}
\longrightarrow T
\qquad\text{as }\eta_y\to\infty.
\]
Hence the strict condition $\frac{\gamma_0^2}{\|\sigma^{-1}\|_F^2}
>
2Ts(2s-1)$
continues to give the required $O_p(1)$ bounds after
$\eta_a,\eta_y$ are chosen sufficiently large.

For $p<0$ and $p\ge1$, the same argument applies with the
corresponding strict thresholds $6T$ and $\max\{6T,T(2p^2-p)\}.$
The additional polynomial factor $\|a_i\|_2^4$ does not change
the exponential-moment threshold. Moreover, all polynomial
dependence on $H_n$ can be absorbed into $e^{CH_n^2}$, since
for every $\varepsilon>0$, $1+H_n^4\le C_\varepsilon e^{\varepsilon H_n^2}.$
Thus,
\[
\widetilde J_1+\widetilde J_2+
\widetilde J_1^{(h)}+\widetilde J_2^{(h)}
=
O_p(e^{CH_n^2}),
\]
where the exponential factor is absorbed by enlarging $C$.
}

Therefore,
\[
I_2
=
O_p\!\left(
\|\sigma\|_F^4S_nH_n^2e^{CH_n^2}
\right).
\]
Since
\[
S_n
\le
\|\sigma^{-1}\|_F^2
\frac{\|\Delta\|_n^2}{n},
\qquad
H_n^2
\le
\|\sigma^{-1}\|_F^2
\frac{\max_{1\le i\le n}\|\Delta_i\|_2^2}{n},
\]
absorbing $\|\sigma\|_F$ and $\|\sigma^{-1}\|_F$ into $C$ gives
\[
I_2
=
O_p\!\left(
\frac{\|\Delta\|_n^2}{n}
\cdot
\frac{\max_{1\le i\le n}\|\Delta_i\|_2^2}{n}
\exp\!\Big(
C\max_{1\le i\le n}
\frac{\|\Delta_i\|_2^2}{n}
\Big)
\right).
\]
Therefore, 
$$|R_2| = O_p\left(\frac{\left\Vert \Delta\right\Vert_n^2}{n}\exp\left(C\max_{1\ \leq i \leq n}\frac{\left\Vert \Delta_i\right\Vert^2}{n}\right)\right) + O_p\!\left(
\frac{\|\Delta\|_n^2}{n}\cdot \frac{\max_{1\le i\le n}\|\Delta_i\|^2}{n}\;
\exp\!\Big(C\,\max_{1\le i\le n}\frac{\|\Delta_i\|^2}{n}\Big)
\right).$$
\end{proof}
By Lemmas \ref{negalphaR1} and \ref{boundR2}, we have the bounds for the remainder terms with a constant $C, C^* >0$:
\begin{align*}
    &|R_1| + |R_2|\\
    &= O_p\!\left(
\frac{1}{n}\Big(\frac{1}{n}\sum_{i=1}^n \|\Delta_i\|^4\Big)^{\!1/2}\,
\Big[1+\frac{\|\Delta\|_n}{\sqrt{n}}+\frac{1}{n}\Big(\frac{1}{n}\sum_{i=1}^n \|\Delta_i\|^4\Big)^{\!1/2}\Big]\,
\exp\!\Big(C'\,\max_{1\le i\le n}\frac{\|\Delta_i\|^2}{n}\Big)
\right)\\
&+ O_p\left(\frac{\left\Vert \Delta\right\Vert_n^2}{n}\exp\left(C\max_{1\ \leq i \leq n}\frac{\left\Vert \Delta_i\right\Vert^2}{n}\right)\right) + O_p\!\left(
\frac{\|\Delta\|_n^2}{n}\cdot \frac{\max_{1\le i\le n}\|\Delta_i\|^2}{n}\;
\exp\!\Big(C\,\max_{1\le i\le n}\frac{\|\Delta_i\|^2}{n}\Big)
\right).
\end{align*} 
\subsubsection{Part II-ii}

Define the population gradient $G^\star(b)
:=
\int_{\mathbb R^d}
g'\!\left(F^\star(y)\right)
\nabla_bL_T(b,y)
\varphi_T(y)\,dy.$
For $M>0$, define its bounded truncation $G_M^\star(b)
:=
G^\star(b)
\mathbf 1_{\{\|G^\star(b)\|_2\leq M\}}.$
For $\lambda\in\mathbb R$, define the correction direction
\[
\Delta_i^{M}(\lambda)
:=
\lambda G_M^\star(B^{(i)}),
\qquad
i=1,\ldots,n.
\]

\begin{lemma}\label{lem:Ci-moment-only}
Under Assumption~\ref{subGaussianCLT}, the vector field $G^\star$ is
well-defined $\mu^\star$-almost surely and satisfies $\E_{\mu^\star}
\left[
\|G^\star(B)\|_2^2
\right]
<\infty.$
Moreover, for every fixed $M,K<\infty$, uniformly over
$|\lambda|\leq K$, $\frac1n
\sum_{i=1}^n
\|\Delta_i^{M}(\lambda)\|_2^2
\leq
K^2M^2,$ $\frac1n
\sum_{i=1}^n
\|\Delta_i^{M}(\lambda)\|_2^4
\leq
K^4M^4,$
and
\[
\max_{1\leq i\leq n}
\frac{\|\Delta_i^{M}(\lambda)\|_2^2}{n}
\leq
\frac{K^2M^2}{n}
\longrightarrow0.
\]
\end{lemma}

\begin{proof}
Since $g'\!\left(F^\star(y)\right)
=
\frac{1}{1-\alpha}
\left(
\frac{e^{rT}}{k^\star}
\right)^{\frac{1}{1-\alpha}}
F^\star(y)^\beta,$
the square-integrability of $G^\star$ follows from the same
sub-Gaussian Gaussian-integral estimates used in the proof of
Lemma~\ref{lem:L22-fro}, with $\mathbb P_0$ replaced by $\mu^\star$.
The quantitative conditions needed for those estimates are included
in Assumption~\ref{subGaussianCLT}.

The remaining inequalities follow immediately from $\|G_M^\star(b)\|_2\leq M$ and $\Delta_i^M(\lambda)
=
\lambda G_M^\star(B^{(i)}).$
\end{proof}

For fixed $M>0$, define $\mathbb P_{n,\lambda,M}
:=
\mathbb P_n^{\Delta^M(\lambda)}$, $F_{n,M}(\lambda)
:=
J(\mathbb P_{n,\lambda,M})-J(\mu^\star)$, and $c_{n,M}
:=
\frac1n
\sum_{i=1}^n
C_i(n)\cdot G_M^\star(B^{(i)}).$

\begin{lemma}\label{cor:clean-rate-explicit}
Under Assumption~\ref{subGaussianCLT}, for every fixed
$M,K<\infty$,
\[
F_{n,M}(\lambda)
=
F_{n,M}(0)
+
\frac{\lambda}{\sqrt n}c_{n,M}
+
R_{n,M}(\lambda),
\qquad
|\lambda|\leq K,
\]
where $\sup_{|\lambda|\leq K}
|R_{n,M}(\lambda)|
=
O_p\!\left(\frac1n\right).$
Moreover, there exists a random variable
$K_{n,M,K}=O_p(1)$ such that, with probability tending to one, for all $|\lambda_1|,|\lambda_2|\leq K$
\[
\left|
R_{n,M}(\lambda_1)
-
R_{n,M}(\lambda_2)
\right|
\leq
\frac{K_{n,M,K}}{n}
|\lambda_1-\lambda_2|.
\]
\end{lemma}

\begin{proof}
The displayed expansion follows from
Eq.~\eqref{eq:Taylor_final}, evaluated at $\Delta_i
=
\Delta_i^M(\lambda)
=
\lambda G_M^\star(B^{(i)}).$
Indeed,
\[
\frac1{\sqrt n}\frac1n
\sum_{i=1}^n
C_i(n)\cdot\Delta_i^M(\lambda)
=
\frac{\lambda}{\sqrt n}c_{n,M}.
\]

By Lemma~\ref{lem:Ci-moment-only}, uniformly over
$|\lambda|\leq K$, $\frac1n\sum_{i=1}^n
\|\Delta_i^M(\lambda)\|_2^4
\leq
K^4M^4,$
and $\max_{1\leq i\leq n}
\frac{\|\Delta_i^M(\lambda)\|_2^2}{n}
\leq
\frac{K^2M^2}{n}.$
Substituting these bounds into
Lemmas~\ref{negalphaR1} and~\ref{boundR2} gives
\[
\sup_{|\lambda|\leq K}
|R_{n,M}(\lambda)|
=
O_p\!\left(\frac1n\right).
\]

It remains to establish the Lipschitz bound. Write $v_i:=G_M^\star(B^{(i)})$ 
and define $m_{n,\lambda}(y)
:=
\frac1n
\sum_{i=1}^n
L_T\!\left(
B^{(i)}+\frac{\lambda v_i}{\sqrt n},
y
\right).$
Then
\[
\dot m_{n,\lambda}(y)
=
\frac1{\sqrt n}\frac1n
\sum_{i=1}^n
\nabla_bL_T\!\left(
B^{(i)}+\frac{\lambda v_i}{\sqrt n},
y
\right)\cdot v_i,
\]
and
\[
\ddot m_{n,\lambda}(y)
=
\frac1{n^2}
\sum_{i=1}^n
v_i^\top
\nabla_b^2L_T\!\left(
B^{(i)}+\frac{\lambda v_i}{\sqrt n},
y
\right)
v_i.
\]
Consequently,
\[
F_{n,M}''(\lambda)
=
\int_{\mathbb R^d}
\left[
g''\!\left(m_{n,\lambda}(y)\right)
\dot m_{n,\lambda}(y)^2
+
g'\!\left(m_{n,\lambda}(y)\right)
\ddot m_{n,\lambda}(y)
\right]
\varphi_T(y)\,dy.
\]
The two terms on the right are of the same form as the terms bounded
in Lemmas~\ref{negalphaR1} and~\ref{boundR2}. Since $\|v_i\|_2\leq M$
and $|\lambda|\leq K$, the same sub-Gaussian estimates yield
\[
\sup_{|\lambda|\leq K}
n|F_{n,M}''(\lambda)|
=
O_p(1).
\]
Taylor's integral remainder formula therefore gives
\[
\sup_{|\lambda|\leq K}
|R_{n,M}(\lambda)|
=
O_p\!\left(\frac1n\right),
\]
and the mean-value theorem gives a random
$K_{n,M,K}=O_p(1)$ such that
\[
|R_{n,M}(\lambda_1)-R_{n,M}(\lambda_2)|
\leq
\frac{K_{n,M,K}}{n}
|\lambda_1-\lambda_2|.
\]
\end{proof}

\subsection{Part III}

Since $\mathbb P_{n,0,M}=\mathbb P_n,$
we have $F_{n,M}(0)
=
J(\mathbb P_n)-J(\mu^\star).$
Define
\[
d_{n,M}
:=
\frac1n
\sum_{i=1}^n
\left\|
G_M^\star(B^{(i)})
\right\|_2^2
\]
and
\[
c_M
:=
\E_{\mu^\star}
\left[
G^\star(B)\cdot G_M^\star(B)
\right].
\]
By the definition of $G_M^\star$,
\[
c_M
=
\E_{\mu^\star}
\left[
\|G^\star(B)\|_2^2
\mathbf 1_{\{\|G^\star(B)\|_2\leq M\}}
\right].
\]
Therefore,
\[
c_M\uparrow c^\star
\qquad
\text{as }M\to\infty.
\]
Since $c^\star>0$, we may fix $M$ sufficiently large such that $c_M>0.$
In Part IV, we will show that
\[
c_{n,M}\xrightarrow{p}c_M
\qquad\text{and}\qquad
d_{n,M}\xrightarrow{p}c_M.
\]
In particular, $\mathbb P_{\mu^\star}(c_{n,M}>0)\longrightarrow1.$

Ignoring the remainder term temporarily, the linearized correction
equation is
\[
F_{n,M}(0)
+
\frac{\lambda}{\sqrt n}c_{n,M}
=
0.
\]
On the event $\{c_{n,M}>0\}$, its solution is $\lambda_{n,M}^{\mathrm{lin}}
:=
-
\frac{\sqrt n\,F_{n,M}(0)}
{c_{n,M}}.$
The corresponding bounded correction direction is
\begin{equation}\label{eq:correction}
\Delta_{i,n,M}^{\mathrm{lin}}
:=
\lambda_{n,M}^{\mathrm{lin}}
G_M^\star(B^{(i)})
=
-
\frac{\sqrt n\,F_{n,M}(0)}
{c_{n,M}}
G_M^\star(B^{(i)}),
\qquad
i=1,\ldots,n.
\end{equation}
By construction,
\[
F_{n,M}(0)
+
\frac1{\sqrt n}\frac1n
\sum_{i=1}^n
C_i(n)\cdot
\Delta_{i,n,M}^{\mathrm{lin}}
=
0.
\]

The squared average norm of the linearized correction is
\[
\left\|\Delta_{n,M}^{\mathrm{lin}}\right\|_n^2
=
\frac{1}{n}\sum_{i=1}^n
\left\|\Delta_{i,n,M}^{\mathrm{lin}}\right\|_2^2
=
\left(\lambda_{n,M}^{\mathrm{lin}}\right)^2 d_{n,M}
=
nF_{n,M}(0)^2\frac{d_{n,M}}{c_{n,M}^2}.
\]
The correction in \eqref{eq:correction} solves the linearized equation.
In Part V, using the Lipschitz remainder estimate from
Lemma~\ref{cor:clean-rate-explicit}, we will construct an exact
correction $\lambda_{n,M}^\star$ satisfying $F_{n,M}(\lambda_{n,M}^\star)=0$
and show that $\lambda_{n,M}^\star
-
\lambda_{n,M}^{\mathrm{lin}}
=
O_p(n^{-1/2}).$

\subsection{Part IV}

In this part, we present
Lemmas~\ref{convergetoc}, \ref{JCLT}, and \ref{Frechet}
to establish the convergence results used in Parts III and V.

\begin{lemma}\label{convergetoc}
For every fixed $M>0$, under Assumption~\ref{subGaussianCLT},
\[
c_{n,M}\xrightarrow{p}c_M
\qquad\text{and}\qquad
d_{n,M}\xrightarrow{p}c_M.
\]
\end{lemma}

\begin{proof}
Since
\[
G_M^\star(B)
=
G^\star(B)
\mathbf 1_{\{\|G^\star(B)\|_2\leq M\}},
\]
we have
\[
G^\star(B)\cdot G_M^\star(B)
=
\|G_M^\star(B)\|_2^2.
\]
Therefore,
\[
c_M
=
\E_{\mu^\star}
\left[
\|G_M^\star(B)\|_2^2
\right].
\]
Since
\[
\|G_M^\star(B)\|_2^2\leq M^2,
\]
the law of large numbers gives
\[
d_{n,M}
=
\frac1n
\sum_{i=1}^n
\|G_M^\star(B^{(i)})\|_2^2
\xrightarrow{p}
c_M.
\]

It remains to compare $c_{n,M}$ and $d_{n,M}$. Define
\[
\Xi_{n,M}(y)
:=
\frac1n
\sum_{i=1}^n
\nabla_bL_T(B^{(i)},y)
\cdot
G_M^\star(B^{(i)}).
\]
By the definitions of $C_i(n)$ and $G^\star$,
\[
c_{n,M}
=
\int_{\mathbb R^d}
g'\!\left(m_{\mathbb P_n}(y)\right)
\Xi_{n,M}(y)
\varphi_T(y)\,dy,
\]
whereas
\[
d_{n,M}
=
\int_{\mathbb R^d}
g'\!\left(F^\star(y)\right)
\Xi_{n,M}(y)
\varphi_T(y)\,dy.
\]
Consequently, by the Cauchy--Schwarz inequality,
\begin{equation}\label{eq:cnM-dnM-bound}
|c_{n,M}-d_{n,M}|
\leq
\left(
\int_{\mathbb R^d}
\left|
g'\!\left(m_{\mathbb P_n}(y)\right)
-
g'\!\left(F^\star(y)\right)
\right|^2
\varphi_T(y)\,dy
\right)^{1/2}
\left(
\int_{\mathbb R^d}
|\Xi_{n,M}(y)|^2
\varphi_T(y)\,dy
\right)^{1/2}.
\end{equation}

We first show that
\[
\|\Xi_{n,M}\|_{L^2(\varphi_T)}
=
O_p(1).
\]
By Jensen's inequality and
$\|G_M^\star(B^{(i)})\|_2\leq M$,
\begin{align*}
\int_{\mathbb R^d}
|\Xi_{n,M}(y)|^2
\varphi_T(y)\,dy
&\leq
\frac1n
\sum_{i=1}^n
\int_{\mathbb R^d}
\left|
\nabla_bL_T(B^{(i)},y)
\cdot
G_M^\star(B^{(i)})
\right|^2
\varphi_T(y)\,dy
\\
&\leq
\frac{M^2}{n}
\sum_{i=1}^n
\int_{\mathbb R^d}
\|\nabla_bL_T(B^{(i)},y)\|_2^2
\varphi_T(y)\,dy.
\end{align*}
The Gaussian completion-of-squares calculation gives
\[
\int_{\mathbb R^d}
\|\nabla_bL_T(b,y)\|_2^2
\varphi_T(y)\,dy
\leq
C
\left(
1+\|a(b)\|_2^2
\right)
\exp\!\left(
T\|a(b)\|_2^2
\right).
\]
The expectation of the right-hand side under $\mu^\star$ is finite
under Assumption~\ref{subGaussianCLT}. Hence, by the law of large
numbers,
\[
\int_{\mathbb R^d}
|\Xi_{n,M}(y)|^2
\varphi_T(y)\,dy
=
O_p(1).
\]

We next prove
\begin{equation}\label{eq:gprime-mixture-convergence}
\left\|
g'\!\left(m_{\mathbb P_n}\right)
-
g'\!\left(F^\star\right)
\right\|_{L^2(\varphi_T)}
\xrightarrow{p}
0.
\end{equation}
Suppose first that $\beta\geq0$. For every fixed
$y\in\mathbb R^d$, the law of large numbers gives $m_{\mathbb P_n}(y)
\xrightarrow{\mathrm{a.s.}}
F^\star(y).$
By Fubini's theorem, there exists an event of
$\mathbb P_{\mu^\star}$-probability one on which this convergence
holds for $\varphi_T(y)\,dy$-almost every $y$.

When $\beta=0$, $g'$ is constant, so
\eqref{eq:gprime-mixture-convergence} is immediate. Assume henceforth
that $\beta>0$. Choose $\eta>0$ sufficiently small and set $\rho:=2\beta(1+\eta).$
By the strict inequalities in
Assumption~\ref{subGaussianCLT}, $\eta$ may be chosen so that, whenever 
$\rho\geq1$,
\[
T(\rho^2-\rho)\|\sigma^{-1}\|_F^2
<
\gamma_0^2.
\]
For $\rho\geq1$, Jensen's inequality gives
\[
m_{\mathbb P_n}(y)^\rho
\leq
\frac1n
\sum_{i=1}^n
L_T(B^{(i)},y)^\rho,
\]
and hence
\[
\E_{\mathbb P_{\mu^\star}}
\left[
m_{\mathbb P_n}(y)^\rho
\right]
\leq
\E_{\mu^\star}
\left[
L_T(B,y)^\rho
\right].
\]
Moreover,
\[
F^\star(y)^\rho
\leq
\E_{\mu^\star}
\left[
L_T(B,y)^\rho
\right].
\]
Using
\[
\int_{\mathbb R^d}
L_T(b,y)^\rho
\varphi_T(y)\,dy
=
\exp\!\left(
\frac{T}{2}
(\rho^2-\rho)
\|a(b)\|_2^2
\right),
\]
the sub-Gaussian assumption yields
\[
\sup_{n\geq1}
\E_{\mathbb P_{\mu^\star}}
\left[
\int_{\mathbb R^d}
m_{\mathbb P_n}(y)^\rho
\varphi_T(y)\,dy
\right]
<\infty
\]
and
\[
\int_{\mathbb R^d}
F^\star(y)^\rho
\varphi_T(y)\,dy
<\infty.
\]

If $0<\rho<1$, concavity gives
\[
\E_{\mathbb P_{\mu^\star}}
\left[
m_{\mathbb P_n}(y)^\rho
\right]
\leq
F^\star(y)^\rho,
\]
and, since $\varphi_T(y)\,dy$ is a probability measure,
\[
\int_{\mathbb R^d}
F^\star(y)^\rho
\varphi_T(y)\,dy
\leq
\left(
\int_{\mathbb R^d}
F^\star(y)
\varphi_T(y)\,dy
\right)^\rho
=
1.
\]
Thus in either case,
\[
\sup_{n\geq1}
\E_{\mathbb P_{\mu^\star}}
\left[
\int_{\mathbb R^d}
\left|
g'\!\left(m_{\mathbb P_n}(y)\right)
-
g'\!\left(F^\star(y)\right)
\right|^{2(1+\eta)}
\varphi_T(y)\,dy
\right]
<\infty.
\]
Therefore, $\left\{
\left|
g'\!\left(m_{\mathbb P_n}(y)\right)
-
g'\!\left(F^\star(y)\right)
\right|^2
\right\}_{n\geq1}$
is uniformly integrable on the product space with measure $\mathbb P_{\mu^\star}
\otimes
\bigl(\varphi_T(y)\,dy\bigr).$
Since the integrands converge to zero almost everywhere on this product
space, Vitali's theorem gives
\[
\E_{\mathbb P_{\mu^\star}}
\left[
\left\|
g'\!\left(m_{\mathbb P_n}\right)
-
g'\!\left(F^\star\right)
\right\|_{L^2(\varphi_T)}^2
\right]
\longrightarrow0.
\]
It follows that
\[
\left\|
g'\!\left(m_{\mathbb P_n}\right)
-
g'\!\left(F^\star\right)
\right\|_{L^2(\varphi_T)}
\xrightarrow{p}0.
\]

Now suppose that $\beta<0$. Choose $R>0$ such that
\[
p_R
:=
\mu^\star
\left(
\|a(B)\|_2\leq R
\right)
>0,
\]
and define
\[
\widehat p_{n,R}
:=
\frac1n
\sum_{i=1}^n
\mathbf 1_{\{\|a(B^{(i)})\|_2\leq R\}}.
\]
By the law of large numbers,
\[
\widehat p_{n,R}
\xrightarrow{\mathrm{a.s.}}
p_R.
\]
On the event $\widehat p_{n,R}\geq\frac{p_R}{2},$
which eventually holds almost surely,
\[
m_{\mathbb P_n}(y)
\geq
\frac{p_R}{2}
\exp\!\left(
-R\|y\|_2-\frac{T}{2}R^2
\right).
\]
Similarly,
\[
F^\star(y)
\geq
p_R
\exp\!\left(
-R\|y\|_2-\frac{T}{2}R^2
\right).
\]
Since $\beta<0$, these lower bounds imply
\[
\left|
g'\!\left(m_{\mathbb P_n}(y)\right)
-
g'\!\left(F^\star(y)\right)
\right|^2
\leq
C_R
\exp\!\left(
2|\beta|R\|y\|_2
\right)
\]
for all sufficiently large $n$, almost surely. The right-hand side is
integrable against the Gaussian density $\varphi_T$. Since $m_{\mathbb P_n}(y)
\longrightarrow
F^\star(y)$
for $\varphi_T(y)\,dy$-almost every $y$, dominated convergence yields
\[
\left\|
g'\!\left(m_{\mathbb P_n}\right)
-
g'\!\left(F^\star\right)
\right\|_{L^2(\varphi_T)}
\longrightarrow0
\]
almost surely, and hence in probability. This proves
\eqref{eq:gprime-mixture-convergence} in all cases.

Combining \eqref{eq:cnM-dnM-bound} and
\eqref{eq:gprime-mixture-convergence}, we obtain $c_{n,M}-d_{n,M}
\xrightarrow{p}0.$
Since $d_{n,M}\xrightarrow{p}c_M,$
then $c_{n,M}\xrightarrow{p}c_M.$
\end{proof}

\begin{lemma}\label{JCLT}
Under Assumption~\ref{subGaussianCLT},
\[
J(\mathbb P_n)-J(\mu^\star)
=
O_p(n^{-1/2}),
\]
and
\[
\sqrt n
\left(
J(\mathbb P_n)-J(\mu^\star)
\right)
\Rightarrow
\mathcal N\!\left(
0,h(\mu^\star)
\right),
\]
where
\[
h(\mu^\star)
=
\int_{\mathbb R^d}\int_{\mathbb R^d}
g'\!\left(F^\star(y_1)\right)
g'\!\left(F^\star(y_2)\right)
\operatorname{Cov}_{\mu^\star}
\left(
L_T(B,y_1),
L_T(B,y_2)
\right)
\varphi_T(y_1)\varphi_T(y_2)
\,dy_1\,dy_2
<\infty.
\]
\end{lemma}

\begin{proof}
Define the centered influence function
\[
\psi^\star(b)
:=
\int_{\mathbb R^d}
g'\!\left(F^\star(y)\right)
\left(
L_T(b,y)-F^\star(y)
\right)
\varphi_T(y)\,dy.
\]
By Fubini's theorem, $\E_{\mu^\star}
\left[
\psi^\star(B)
\right]
=0.$
The sub-Gaussian Gaussian-integral estimates in
Assumption~\ref{subGaussianCLT} imply $\E_{\mu^\star}
\left[
\psi^\star(B)^2
\right]
<\infty.$
Moreover, $\operatorname{Var}_{\mu^\star}
\left(
\psi^\star(B)
\right)
=
h(\mu^\star).$
By Lemma~\ref{Frechet},
\[
J(\mathbb P_n)-J(\mu^\star)
=
\frac1n
\sum_{i=1}^n
\psi^\star(B^{(i)})
+
o_p(n^{-1/2}).
\]
The classical central limit theorem therefore gives
\[
\sqrt n
\left(
J(\mathbb P_n)-J(\mu^\star)
\right)
\Rightarrow
\mathcal N
\left(
0,h(\mu^\star)
\right).
\]
In particular,
\[
J(\mathbb P_n)-J(\mu^\star)
=
O_p(n^{-1/2}).
\]
\end{proof}

\begin{lemma}\label{Frechet}
Let
\[
D_n(y)
:=
m_{\mathbb P_n}(y)-F^\star(y).
\]
Under Assumption~\ref{subGaussianCLT},
\[
J(\mathbb P_n)-J(\mu^\star)
=
\int_{\mathbb R^d}
g'\!\left(F^\star(y)\right)
D_n(y)
\varphi_T(y)\,dy
+
o_p(n^{-1/2}).
\]
Equivalently,
\[
J(\mathbb P_n)-J(\mu^\star)
=
\frac1n
\sum_{i=1}^n
\psi^\star(B^{(i)})
+
o_p(n^{-1/2}),
\]
where
\[
\psi^\star(b)
=
\int_{\mathbb R^d}
g'\!\left(F^\star(y)\right)
\left(
L_T(b,y)-F^\star(y)
\right)
\varphi_T(y)\,dy.
\]
\end{lemma}

\begin{proof}
Since both $F^\star(y)$ and $m_{\mathbb P_n}(y)$ are strictly
positive, Taylor's integral formula gives
\begin{align}
&
J(\mathbb P_n)-J(\mu^\star)
-
\int_{\mathbb R^d}
g'\!\left(F^\star(y)\right)
D_n(y)
\varphi_T(y)\,dy
\nonumber\\
&\quad=
\int_0^1(1-t)
\int_{\mathbb R^d}
g''\!\left(
F^\star(y)+tD_n(y)
\right)
D_n(y)^2
\varphi_T(y)\,dy\,dt.
\label{eq:J-von-Mises-remainder}
\end{align}

We show that the right-hand side of
\eqref{eq:J-von-Mises-remainder} is $o_p(n^{-1/2})$.

Suppose first that $\beta<1$. Choose $R>0$ such that
\[
p_R
:=
\mu^\star\!\left(
\|a(B)\|_2\leq R
\right)
>0
\]
and define
\[
\widehat p_{n,R}
:=
\frac1n
\sum_{i=1}^n
\mathbf 1_{\{\|a(B^{(i)})\|_2\leq R\}}.
\]
On the event $\widehat p_{n,R}\geq\frac{p_R}{2},$
whose probability tends to one, we have, for every
$t\in[0,1]$,
\[
F^\star(y)+tD_n(y)
=
(1-t)F^\star(y)+tm_{\mathbb P_n}(y)
\geq
c_R e^{-R\|y\|_2}
\]
for a deterministic constant $c_R>0$. Since
\[
g''(x)
=
\frac{\alpha}{(1-\alpha)^2}
\left(
\frac{e^{rT}}{k^\star}
\right)^{\frac{1}{1-\alpha}}
x^{\beta-1},
\]
it follows that
\[
\left|
g''\!\left(
F^\star(y)+tD_n(y)
\right)
\right|
\leq
C_R
e^{(1-\beta)R\|y\|_2}.
\]
Therefore, on this event, the absolute value of the remainder in
\eqref{eq:J-von-Mises-remainder} is bounded by
\[
C_R
\int_{\mathbb R^d}
e^{(1-\beta)R\|y\|_2}
D_n(y)^2
\varphi_T(y)\,dy.
\]
Moreover,
\[
\E_{\mathbb P_{\mu^\star}}
\left[
D_n(y)^2
\right]
=
\frac1n
\operatorname{Var}_{\mu^\star}
\left(
L_T(B,y)
\right).
\]
Hence,
\begin{align*}
&
\E_{\mathbb P_{\mu^\star}}
\left[
\int_{\mathbb R^d}
e^{(1-\beta)R\|y\|_2}
D_n(y)^2
\varphi_T(y)\,dy
\right]
\\
&\qquad=
\frac1n
\int_{\mathbb R^d}
e^{(1-\beta)R\|y\|_2}
\operatorname{Var}_{\mu^\star}
\left(
L_T(B,y)
\right)
\varphi_T(y)\,dy
\\
&\qquad\leq
\frac{C}{n},
\end{align*}
where the last integral is finite under
Assumption~\ref{subGaussianCLT}. Thus, when $\beta<1$, the
remainder is $O_p(n^{-1})$ and therefore
$o_p(n^{-1/2})$.

Now suppose that $\beta\geq1$. For $x,z>0$ and $t\in[0,1]$,
\[
\left(
(1-t)x+tz
\right)^{\beta-1}
\leq
C_\beta
\left(
x^{\beta-1}+z^{\beta-1}
\right).
\]
Consequently,
\[
\left|
g''\!\left(
F^\star(y)+tD_n(y)
\right)
\right|
\leq
C_\beta
\left(
F^\star(y)^{\beta-1}
+
m_{\mathbb P_n}(y)^{\beta-1}
\right).
\]
It therefore suffices to control
\[
Q_n
:=
\int_{\mathbb R^d}
\left(
F^\star(y)^{\beta-1}
+
m_{\mathbb P_n}(y)^{\beta-1}
\right)
D_n(y)^2
\varphi_T(y)\,dy.
\]

For every fixed $y$, the standard fourth-moment bound for empirical
averages gives
\[
\E_{\mathbb P_{\mu^\star}}
\left[
|D_n(y)|^4
\right]
\leq
\frac{C}{n^2}
\E_{\mu^\star}
\left[
\left|
L_T(B,y)-F^\star(y)
\right|^4
\right].
\]
By Cauchy--Schwarz, $\E_{\mathbb P_{\mu^\star}}
\left[
F^\star(y)^{\beta-1}
D_n(y)^2
\right] \leq
F^\star(y)^{\beta-1}
\left(
\E_{\mathbb P_{\mu^\star}}
|D_n(y)|^4
\right)^{1/2},$
and $$\E_{\mathbb P_{\mu^\star}}
\left[
m_{\mathbb P_n}(y)^{\beta-1}
D_n(y)^2
\right] \leq
\left(
\E_{\mathbb P_{\mu^\star}}
m_{\mathbb P_n}(y)^{2(\beta-1)}
\right)^{1/2}
\left(
\E_{\mathbb P_{\mu^\star}}
|D_n(y)|^4
\right)^{1/2}.$$
The Gaussian-integral and sub-Gaussian moment bounds in
Assumption~\ref{subGaussianCLT} imply that the resulting
$y$-integrals are finite uniformly in $n$. Hence, $\E_{\mathbb P_{\mu^\star}}[Q_n]
\leq
\frac{C}{n}.$
Thus, $Q_n=O_p(n^{-1}),$
and the remainder is again $o_p(n^{-1/2})$.

Combining the two cases proves
\[
J(\mathbb P_n)-J(\mu^\star)
=
\int_{\mathbb R^d}
g'\!\left(F^\star(y)\right)
D_n(y)
\varphi_T(y)\,dy
+
o_p(n^{-1/2}).
\]
Finally,
\[
D_n(y)
=
\frac1n
\sum_{i=1}^n
\left(
L_T(B^{(i)},y)-F^\star(y)
\right),
\]
so Fubini's theorem gives
\[
\int_{\mathbb R^d}
g'\!\left(F^\star(y)\right)
D_n(y)
\varphi_T(y)\,dy
=
\frac1n
\sum_{i=1}^n
\psi^\star(B^{(i)}).
\]
\end{proof}

By Lemmas~\ref{convergetoc} and~\ref{JCLT}, for every fixed
sufficiently large $M$ such that $c_M>0$,
\[
\sqrt n\,F_{n,M}(0)=O_p(1),
\qquad
c_{n,M}\xrightarrow{p}c_M,
\qquad
d_{n,M}\xrightarrow{p}c_M.
\]
Consequently,
\[
\lambda_{n,M}^{\mathrm{lin}}
=
-
\frac{\sqrt n\,F_{n,M}(0)}
{c_{n,M}}
=
O_p(1),
\]
and
\[
\left\|
\Delta_{n,M}^{\mathrm{lin}}
\right\|_n^2
=
nF_{n,M}(0)^2
\frac{d_{n,M}}{c_{n,M}^2}
=
O_p(1).
\]

\subsection{Part V}

Fix $M>0$ sufficiently large that $c_M>0.$
For every fixed $K<\infty$, define $A_{n,M,K}
:=
n\sup_{|\lambda|\leq K}
|R_{n,M}(\lambda)|.$
Then $A_{n,M,K}=O_p(1).$
Moreover, by Lemma~\ref{cor:clean-rate-explicit}, there exists
$B_{n,M,K}=O_p(1)$ such that, with probability tending to one,
\begin{align}
\sup_{|\lambda|\leq K}
|R_{n,M}(\lambda)|
&\leq
\frac{A_{n,M,K}}{n},
\label{eq:Rn-mag}
\\
|R_{n,M}(\lambda_1)-R_{n,M}(\lambda_2)|
&\leq
\frac{B_{n,M,K}}{n}
|\lambda_1-\lambda_2|,
\qquad
|\lambda_1|,|\lambda_2|\leq K.
\label{eq:Rn-lip}
\end{align}

Define
\[
\mathcal T_{n,M}(\lambda)
:=
-
\frac{
\sqrt n\,F_{n,M}(0)
+
\sqrt n\,R_{n,M}(\lambda)
}{
c_{n,M}
}.
\]
A fixed point of $\mathcal T_{n,M}$ is an exact solution of $F_{n,M}(\lambda)=0.$

We next verify the contraction conditions. Fix $\varepsilon>0$ and 
choose $c_0\in(0,c_M).$
By Lemmas~\ref{convergetoc} and~\ref{JCLT},
\[
c_{n,M}\xrightarrow{p}c_M
\qquad\text{and}\qquad
\sqrt n\,F_{n,M}(0)=O_p(1).
\]
Hence, we may choose $K<\infty$ sufficiently large and constants
$A_0,B_0<\infty$ such that, for all sufficiently large $n$, with
$\mathbb P_{\mu^\star}$-probability at least $1-\varepsilon$,
\[
c_{n,M}\geq c_0,
\qquad
\left|
\sqrt n\,F_{n,M}(0)
\right|
\leq
\frac{Kc_0}{2},
\qquad
A_{n,M,K}\leq A_0,
\qquad
B_{n,M,K}\leq B_0.
\]
Choose $n$ further sufficiently large that
\[
\frac{A_0}{c_0\sqrt n}
\leq
\frac K2
\qquad\text{and}\qquad
\frac{B_0}{c_0\sqrt n}
\leq
\frac12.
\]

For $|\lambda|\leq K$, Eq.~\eqref{eq:Rn-mag} gives $$|\mathcal T_{n,M}(\lambda)| \leq \frac{
\left|
\sqrt n\,F_{n,M}(0)
\right|
}{
c_{n,M}
}
+
\frac{
\sqrt n\,|R_{n,M}(\lambda)|
}{
c_{n,M}
} \leq \frac K2
+
\frac{A_0}{c_0\sqrt n} \leq K$$
Thus, $\mathcal T_{n,M}([-K,K])
\subseteq
[-K,K].$

Furthermore, for
$|\lambda_1|,|\lambda_2|\leq K$, Eq.~\eqref{eq:Rn-lip} gives
\[
|\mathcal T_{n,M}(\lambda_1)-\mathcal T_{n,M}(\lambda_2)|
=
\frac{\sqrt n}{c_{n,M}}
|R_{n,M}(\lambda_1)-R_{n,M}(\lambda_2)|
\leq
\frac{B_0}{c_0\sqrt n}
|\lambda_1-\lambda_2|
\leq
\frac12|\lambda_1-\lambda_2|.
\]
Therefore, with probability at least $1-\varepsilon$ for all
sufficiently large $n$, $\mathcal T_{n,M}$ is a contraction from
$[-K,K]$ into itself. By the Banach fixed-point theorem, there exists
a unique $\lambda_{n,M}^\star\in[-K,K]$ such that $\mathcal T_{n,M}(\lambda_{n,M}^\star)
=
\lambda_{n,M}^\star.$
Equivalently, $F_{n,M}(\lambda_{n,M}^\star)=0.$

Recall the linearized correction from Part III: $\lambda_{n,M}^{\mathrm{lin}}
=
-
\frac{
\sqrt n\,F_{n,M}(0)
}{
c_{n,M}
}.$
The fixed-point equation gives
\[
\lambda_{n,M}^\star
-
\lambda_{n,M}^{\mathrm{lin}}
=
-
\frac{
\sqrt n\,R_{n,M}(\lambda_{n,M}^\star)
}{
c_{n,M}
}.
\]
Since $\lambda_{n,M}^\star\in[-K,K]$, Eq.~\eqref{eq:Rn-mag} implies $\lambda_{n,M}^\star
-
\lambda_{n,M}^{\mathrm{lin}}
=
O_p(n^{-1/2}).$

Define
\[
\Delta_{i,n,M}^\star
:=
\lambda_{n,M}^\star
G_M^\star(B^{(i)}),
\qquad
i=1,\ldots,n.
\]
A standard diagonal selection over $\varepsilon\downarrow0$ yields
events $E_{n,M}$ satisfying $\mathbb P_{\mu^\star}(E_{n,M})\longrightarrow1$
and a measurable choice of $\lambda_{n,M}^\star$ such that, on
$E_{n,M}$,
\begin{equation}\label{eq:existence}
J\!\left(
\mathbb P_{n,\lambda_{n,M}^\star,M}
\right)
=
J(\mu^\star),
\qquad
\left\|
\Delta_{n,M}^\star
\right\|_n
=
O_p(1),
\qquad
\lambda_{n,M}^\star
-
\lambda_{n,M}^{\mathrm{lin}}
=
O_p(n^{-1/2}).
\end{equation}

We now evaluate the exponential transport cost. On $E_{n,M}$, the
atomwise coupling that sends
\[
B^{(i)}
\quad\text{to}\quad
B^{(i)}
+
\frac{
\lambda_{n,M}^\star
G_M^\star(B^{(i)})
}{
\sqrt n
}
\]
gives
\[
\begin{aligned}
D_c\!\left(
\mathbb P_n,
\mathbb P_{n,\lambda_{n,M}^\star,M}
\right)
&\leq
\frac1n
\sum_{i=1}^n
\left[
\exp\!\left(
\frac{
\tau(\lambda_{n,M}^\star)^2
\|G_M^\star(B^{(i)})\|_2^2
}{n}
\right)
-1
\right].
\end{aligned}
\]
Since
\[
\|G_M^\star(B^{(i)})\|_2\leq M
\qquad\text{and}\qquad
\lambda_{n,M}^\star=O_p(1),
\]
we have
\[
\frac{
\tau(\lambda_{n,M}^\star)^2M^2
}{n}
=
o_p(1).
\]
Using
\[
e^x-1\leq xe^x,
\qquad x\geq0,
\]
we obtain
\[
\begin{aligned}
nD_c\!\left(
\mathbb P_n,
\mathbb P_{n,\lambda_{n,M}^\star,M}
\right)
&\leq
\tau(\lambda_{n,M}^\star)^2
d_{n,M}
\exp\!\left(
\frac{
\tau(\lambda_{n,M}^\star)^2M^2
}{n}
\right)
\\
&=
\tau(\lambda_{n,M}^\star)^2d_{n,M}
+
o_p(1).
\end{aligned}
\]

Because
\[
\lambda_{n,M}^\star
-
\lambda_{n,M}^{\mathrm{lin}}
=
O_p(n^{-1/2}),
\qquad
\lambda_{n,M}^{\mathrm{lin}}=O_p(1),
\qquad
d_{n,M}=O_p(1),
\]
we have
\[
(\lambda_{n,M}^\star)^2d_{n,M}
=
(\lambda_{n,M}^{\mathrm{lin}})^2d_{n,M}
+
o_p(1).
\]
By the definition of $\lambda_{n,M}^{\mathrm{lin}}$,
\[
(\lambda_{n,M}^{\mathrm{lin}})^2d_{n,M}
=
nF_{n,M}(0)^2
\frac{d_{n,M}}{c_{n,M}^2}.
\]
Consequently, there exists
\[
\varepsilon_{n,M}=o_p(1)
\]
such that, on $E_{n,M}$,
\begin{equation}\label{eq:fixed-M-upper-bound}
nR_n(k^\star)
\leq
\tau nF_{n,M}(0)^2
\frac{d_{n,M}}{c_{n,M}^2}
+
\varepsilon_{n,M}.
\end{equation}

\subsection{Part VI}

We first verify that the single sub-Gaussian condition in
Assumption~\ref{subGaussianCLT} implies all auxiliary moment
conditions used in the preceding parts.

\begin{lemma}\label{lem:section5-assumptions}
Assumption~\ref{subGaussianCLT} implies
Assumptions~\ref{Lemma12ass} and~\ref{Lemma13ass}, as well as the
Gaussian-moment conditions used in
Lemmas~\ref{lem:Ci-moment-only},
\ref{cor:clean-rate-explicit},
\ref{convergetoc},
\ref{JCLT}, and~\ref{Frechet}.
\end{lemma}

\begin{proof}
Assumption~\ref{Lemma12ass} requires
\[
\frac{\gamma_0^2}{\|\sigma^{-1}\|_F^2}
>
T\max\left\{
2(2\beta^2-\beta),\,4
\right\}
=
T\max\left\{
4\beta^2-2\beta,\,4
\right\},
\]
which follows immediately from
Assumption~\ref{subGaussianCLT}.
Next, Assumption~\ref{Lemma13ass} requires
\[
\frac{\gamma_0^2}{\|\sigma^{-1}\|_F^2}
>
T\max\left\{
M_{\mathrm{mid}}(\beta),\,
6,\,
2p^2-p
\right\},
\qquad
p:=\beta-1.
\]
We have
\[
2p^2-p
=
2(\beta-1)^2-(\beta-1)
=
2\beta^2-5\beta+3.
\]
Moreover,
\[
(4\beta^2-2\beta)
-
(2\beta^2-5\beta+3)
=
2\beta^2+3\beta-3.
\]
Let
\[
\beta_0
:=
\frac{-3+\sqrt{33}}{4}.
\]
For $\beta\geq\beta_0$,
\[
2\beta^2-5\beta+3
\leq
4\beta^2-2\beta.
\]
For $-1<\beta\leq\beta_0$, the function
\[
2\beta^2-5\beta+3
\]
is bounded above by its value at $\beta=-1$.
Therefore,
\[
2p^2-p
\leq
\max\left\{
4\beta^2-2\beta,\,
16
\right\}.
\]
Also we notice that
Assumption~\ref{subGaussianCLT} implies
Assumption~\ref{Lemma13ass}.

The remaining estimates require Gaussian likelihood moments of the
orders appearing in the terms
\[
4\beta^2-2\beta,
\qquad
M_{\mathrm{mid}}(\beta),
\qquad
16,
\qquad
8\beta+8.
\]
These four quantities are included explicitly in
Assumption~\ref{subGaussianCLT}. Therefore, all Gaussian-integral and
uniform-integrability conditions invoked in the preceding lemmas are
satisfied.
\end{proof}

We now conclude the proof of Theorem~\ref{CLT}. Fix $M$ sufficiently
large that $c_M>0.$
By Lemmas~\ref{JCLT} and~\ref{convergetoc}, $\sqrt n\,F_{n,M}(0)
=
\sqrt n
\left(
J(\mathbb P_n)-J(\mu^\star)
\right)
\Rightarrow
Z,$
where $Z\sim\mathcal N\!\left(0,h(\mu^\star)\right),$ $c_{n,M}\xrightarrow{p}c_M,$ and $d_{n,M}\xrightarrow{p}c_M.$
Hence, by Slutsky's theorem,
\[
\tau nF_{n,M}(0)^2
\frac{d_{n,M}}{c_{n,M}^2}
+
\varepsilon_{n,M}
\Rightarrow
\frac{\tau Z^2}{c_M}.
\]

On the events $E_{n,M}$ from Part V,
Eq.~\eqref{eq:fixed-M-upper-bound} gives
\[
0
\leq
nR_n(k^\star)
\leq
\tau nF_{n,M}(0)^2
\frac{d_{n,M}}{c_{n,M}^2}
+
\varepsilon_{n,M}.
\]
Therefore, for every $t\geq0$,
\[
\begin{aligned}
\mathbb P_{\mu^\star}
\left(
nR_n(k^\star)>t
\right)
&\leq
\mathbb P_{\mu^\star}
\left(
\tau nF_{n,M}(0)^2
\frac{d_{n,M}}{c_{n,M}^2}
+
\varepsilon_{n,M}
\geq t
\right)
+
\mathbb P_{\mu^\star}(E_{n,M}^c).
\end{aligned}
\]
Since $\mathbb P_{\mu^\star}(E_{n,M}^c)\longrightarrow0,$
the Portmanteau theorem yields
\[
\limsup_{n\to\infty}
\mathbb P_{\mu^\star}
\left(
nR_n(k^\star)>t
\right)
\leq
\Pr\!\left(
\frac{\tau Z^2}{c_M}
\geq t
\right).
\]

Finally, by monotone convergence,
\[
c_M
=
\E_{\mu^\star}
\left[
\|G^\star(B)\|_2^2
\mathbf 1_{\{\|G^\star(B)\|_2\leq M\}}
\right]
\uparrow
\E_{\mu^\star}
\left[
\|G^\star(B)\|_2^2
\right]
=
c^\star.
\]
Consequently,
\[
\frac{\tau Z^2}{c_M}
\downarrow
\frac{\tau Z^2}{c^\star}
\qquad
\text{almost surely}.
\]
Letting $M\to\infty$, we obtain, for every continuity point $t$ of
the distribution of $L
:=
\frac{\tau Z^2}{c^\star},$
that
\[
\limsup_{n\to\infty}
\mathbb P_{\mu^\star}
\left(
nR_n(k^\star)>t
\right)
\leq
\Pr(L>t).
\]
This completes the proof of Theorem~\ref{CLT}.

\end{document}